\documentclass[11pt]{amsart}
\usepackage{graphicx}
\usepackage{amsmath}
\usepackage{amsfonts}
\usepackage{amssymb}
\usepackage{a4}
\usepackage{color}
\def\oh{\frac{1}{2}}

\def\oh{\frac{1}{2}}

\begin{document}
\title[High-order well-balanced finite-volume schemes. Draft \today]
      {High-order well-balanced finite-volume schemes for barotropic flows.\\
	Development and numerical comparisons.}

\thanks{
  This joint work was supported by the EU financed network no.
  HPRN-CT-2002-00282 ("Hyke"). The work of N.P. was funded by German
  Science Foundation grant Graduiertenkolleg 775 and that of J.N. by
  the BeMatA program 139144/431 of the Research Council of Norway.}

\author[Pankratz]{Normann Pankratz}
\address[Normann Pankratz]{\newline
  Institut f\"ur Geometrie und Praktische Mathematik, \newline
  RWTH Aachen, 52056 Aachen, Germany}
\email{pankratz@igpm.rwth-aachen.de}
\author[Natvig]{Jostein~R.~Natvig}
\address[Jostein Natvig]{\newline
  SINTEF ICT,  Deptartement of Applied Mathematics, \newline 
  P.O.Box 124. Blindern, N-0314 Oslo, Norway.}
\email{Jostein.R.Natvig@sintef.no}
\author[Gjevik]{Bj\o rn Gjevik}
\address[Bj\o rn Gjevik]{\newline
  Institute of Mathematics, University of Oslo,\newline
  P.O.Box 1053 Blindern, N-0316 Oslo, Norway.}
\email{bjorng@math.uio.no}
\urladdr{http://www.math.uio.no/~bjorng/}
\author[Noelle, DRAFT \today]{Sebastian Noelle}
\address[Sebastian Noelle]{\newline
  Institut f\"ur Geometrie und Praktische Mathematik, \newline
  RWTH Aachen, 52056 Aachen, Germany}
\email{noelle@igpm.rwth-aachen.de}
\urladdr{http://www.igpm.rwth-aachen.de/\~{}noelle/}
%

%\date{to be submitted to Ocean Modelling}

\date{Draft \today}

%=========================================================================
\begin{abstract}
%=========================================================================
  In this paper we compare a classical finite-difference and a high
  order finite-volume scheme for barotropic ocean flows. 
  We compare the schemes with respect to their accuracy, stability, and study various
  outflow and inflow boundary conditions. 
  We apply the schemes to the problem of eddy formation in shelf slope jets along 
the Ormen Lange section of the Norwegian
  shelf. Our results strongly confirm the development of
  mesoscale eddies caused by instability of the flows.
\end{abstract}

%\subjclass{65M06, 35L65, 35L67 !!!!!!!!!!!!!!!!!!!!!!!!!!!!!!!!!!!!!!!!!!!!!}

%\keywords{hyperbolic conservation laws, source terms, shallow water equations,
%  well-balanced schemes, finite-volume schemes, WENO reconstruction}

\maketitle

\tableofcontents
%
%
%
%===============================================================================
\section{Introduction}
%===============================================================================
%
%
%
In medium-scale geophysical fluid flow, with length scales of hundreds
of kilometres, the geometry of the earth, its rotation and curvature
are of great importance.  The modelling of flow phenomena at these
scales involves complex nonlinear equations with extra terms
accounting for the geometry and the rotating frame of reference.

Many geophysical flow problems are shallow in the sense that the waves
length of horizontal motion greatly exceeds the scale of changes in
the vertical direction.  In many cases, this justifies a
simplification of the governing equations for the vertical motion.
The shallow-water equations is one such system where the dependent
variables are depth-averaged and only first-order differential terms
are retained.  In this paper we consider numerical solutions of the
shallow-water system, written as a system of first-order hyperbolic
conservation laws with source terms modelling the effects of variable
bottom and a rotating frame of reference,
\begin{equation}
\label{eq:sw}
\begin{bmatrix}
\eta \\[0.15cm] U \\[0.15cm] V 
\end{bmatrix}_t +
\begin{bmatrix}
U \\[0.15cm] \frac{U^2}{H} + \frac {1}{2}gH^2 \\[0.15cm] \frac{UV}{H}
\end{bmatrix}_x +
\begin{bmatrix}
V \\[0.15cm] \frac{UV}{H} \\[0.15cm] \frac{V^2}{H} + \frac{1}{2}gH^2
\end{bmatrix}_y 
=
\begin{bmatrix}
0 \\[0.15cm]-gH z_x + fV \\[0.15cm] -gH z_y - fU
\end{bmatrix}.
\end{equation}
Here subscripts denote differentiation, $\eta$ is the surface
elevation, $z$ is the bottom topography and $H=\eta-z$ is the total
water depth.  The components of the volume-flux per unit length in 
the $x$- and $y$-direction are $U$ and $V$, respectively.  The source
terms in \eqref{eq:sw} model two different physical effects: the
rotation and the variable bottom topography.  The rotating frame of
reference introduces a Coriolis force $[0, fV, -fU]^T$ acting
transversely and proportionally to the volume-flux.  The other source
term $[0, -gHz_x, -gHz_y]^T$ accounts for the variations in the bottom
topography $z$.  In applications, this barotropic model is used to
study weather systems, mean currents and transport and wave phenomena
in coastal zones, rivers and lakes, in cases where the density
stratification has negligible influence on the flow.

Classically, i.e. at least since the 1940s, such initial value
problems have been solved by finite-difference methods
\cite{Neumann_Richtmyer_1950,Richtmyer_Morton_1967}. To this day, such
methods are the working horse of many models. They are
easy to implement, fast, and for smooth flows they give accurate
results. On the other hand, for non-smooth solutions they suffer from
dispersive oscillations which need to be damped by adding artificial
viscosity.

These stability problems led (roughly from 1950s into the 1990s) to
the development of more robust finite-difference, finite-volume, ENO
and WENO schemes
\cite{Lax_1954,Godunov_1959,vanLeer_1979,Harten_1983,LeVeque_1992,
Shu_1998}. For geophysical flows it was important to develop schemes
which maintain fundamental equilibrium solutions on the discrete
level, the so-called well-balanced schemes (see e.g.
\cite{Audusse_etal_2004,Noelle_etal_2006,Xing_Shu_2005} and the
references therein).  Recently, Bouchut et al.\cite{Bouchut_etal_2004}
have described a technique to obtain a well-balanced discretisation of
the Coriolis terms in the one-dimensional case.  The well-balanced
discretisation preserves geostrophically balanced states exactly at
the discrete level. This technique may be generalised to
two-dimensional jets which are aligned with a Cartesian grid. With
these extensions, well-balanced finite-volume schemes are a very
stable and -- if equipped with high-order reconstructions -- highly
accurate alternative for the computation of depth-averaged geophysical
flows, which may contain shock-, or bore-waves. An advantage of these
schemes is that the solution is damped only in region where damping is
needed.

The present paper reports on the joint work of a researcher, who has
over many years developed and used a finite-difference ocean models
\cite{Gjevik_etal_2002}, an engineer and two numerical analysts who
have developed a high order well-balanced finite-volume scheme
\cite{Noelle_etal_2006}. Our goal is to study and, if possible,
quantify the advantages of either code. We hope that other researchers
will draw some useful conclusions from our results, when they decide
which type of code they should use.

As test-case we study a class of jets along the Norwegian shelf. Such
shelf slope jets have been studied extensively (see
\cite{Gjevik_etal_2002, Thiem_etal_2006} and the references
therein). A series of numerical examples indicates that these currents
can become unstable, in the sense that an initially almost laminar
flow generates strong eddies and oscillations. Linear stability
analysis \cite{Gjevik_etal_2002, Thiem_etal_2006} confirms the
existence of unstable modes. This provides us with a challenging test
problem within a relatively simple topography. Other test problems are
used to study numerical convergence and accuracy.

It will come as no surprise that the setup of analytical and numerical
in- and outflow boundary conditions was one of the main difficulties
in this study.

The outline of the paper is as follows: In
Section~\ref{sec:TheBgridScheme} we give an overview of the
finite-difference method used in \cite{Gjevik_etal_2002}. We rearrange
the temporal update to assure second-order accuracy. In
Sections~\ref{sec:FV_Scheme} we review the high-order well-balanced
finite-volume scheme derived recently in \cite{Noelle_etal_2006}.  In
particular, the source term treatment is described in
Section~\ref{sec:FV_SourceTerms}. The entire Section
\ref{sec:BoundaryConditions} is devoted to boundary conditions,
particularly inflow and absorbing outflow boundary conditions for the
finite-volume scheme. These have a strong impact upon the accuracy and
the flow features computed by our schemes. In
Section~\ref{sec:Comparison}, we evaluate the accuracy, order of
convergence and resolution for various test problems. Then we focus on
the formation of eddies in shelf slope jets. Here we study and
thereby rule out several possible numerical sources of the
instability. 
%Among the boundary conditions we find the free-slip
%balanced inflow boundary condition derived in
%Section~\ref{sec:FV_MomBalInf_BC} particularly convincing. 
We conclude the paper in Section~\ref{sec:Conclusion} by discussing in
detail the advantages of the finite-difference and finite-volume
solvers, the boundary conditions and the eddy formation in the along
shelf current.\\ \textbf{Acknowledgement:} We would like to thank
Roland Sch\"afer for lively and stimulating discussions. Also we would
like to thank Frank Knoben and Markus J\"urgens for there unresting
support for the parallelisation of the scheme.
%
%
%
%=========================================================================================
\section{Discretisation}\label{sec:Discretisation}
%=========================================================================================
%
%
%
In this section we give an overview of the two numerical schemes.  The
philosophies underlying these schemes are quite different.  In the
finite-difference scheme, the solution is approximated by point values
on a grid.  To advance the solution, the derivative of the flux terms
are computed using central differencing and averaging operators.  A
special staggering of variables, called a B-grid, is used, where the
volume-flux is approximated on the mesh $(ih, jh, n\Delta t)$, and the
surface elevation is approximated on a mesh shifted by $h/2$ in each
spatial direction and $\Delta t/2$ in time.  The original scheme of
\cite{Gjevik_etal_2002} is second-order accurate in space, but only
first-order accurate in time.  A simple extension yields a fully
second-order finite-difference scheme.  In actual computations, this
scheme performs very well for smooth solutions.  However, we do
observe spurious oscillations when shocks appear in the solution.  In
Section~\ref{sec:TheBgridScheme} we give a complete description of
this scheme.

In the finite-volume scheme, the solution is approximated in terms of
cell-averages.  These cell averages are advanced in time by computing
fluxes across cell interfaces. To evaluate the fluxes, accurate
point-values of each variable must be reconstructed from cell
averages.  We use a fifth-order WENO procedure \cite{Shu_1988,
Shu_1998} for the reconstruction combined and Roe's approximate
Riemann solver \cite{Roe_1981} for the interface flux.  A standard
fourth-order Runge-Kutta scheme is used as temporal discretisation.
In addition, the scheme is equipped with a high-order well-balanced
discretisation of the geometrical source term \cite{Noelle_etal_2006}.
This scheme has proven to be highly accurate both for smooth and
non-smooth solutions.  In Section~\ref{sec:FV_Scheme} we give an
overview of this scheme, and refer to \cite{Noelle_etal_2006} for a
full description.
%
%==============================================================================================================
\subsection{The Finite-Difference Scheme and a fully second-order accurate extension}\label{sec:TheBgridScheme}
%==============================================================================================================
%
\begin{figure}
 \centering \includegraphics{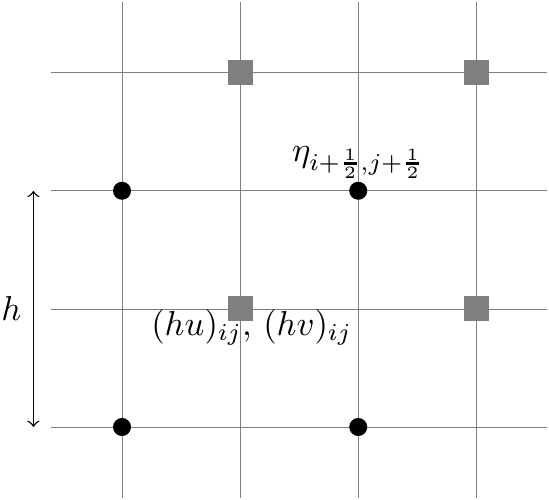}
\caption{\small\label{fig:bgrid} The layout of a B-grid.  The surface
elevation $\eta$ is approximated in the black circles, and the
volume-fluxes $U$ and $V$ in the grey squares.}
\end{figure}

\begin{figure}\label{fig:stencil_fd}
  \includegraphics[width=0.49\linewidth]{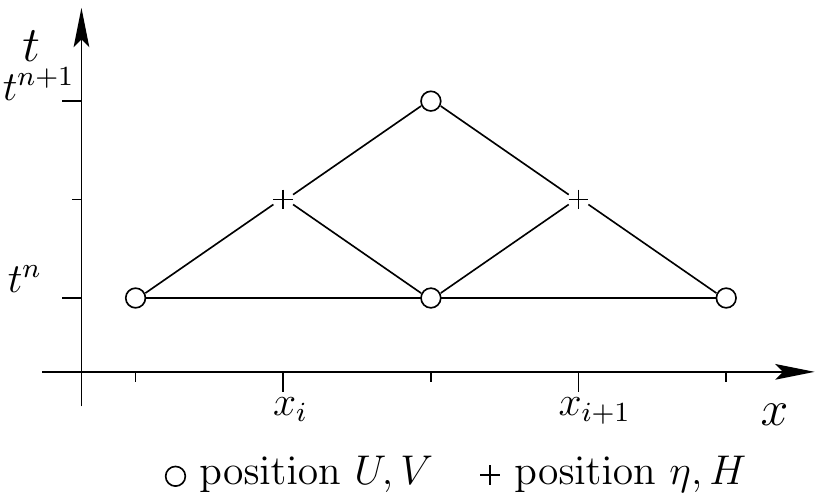}
  \includegraphics[width=0.49\linewidth]{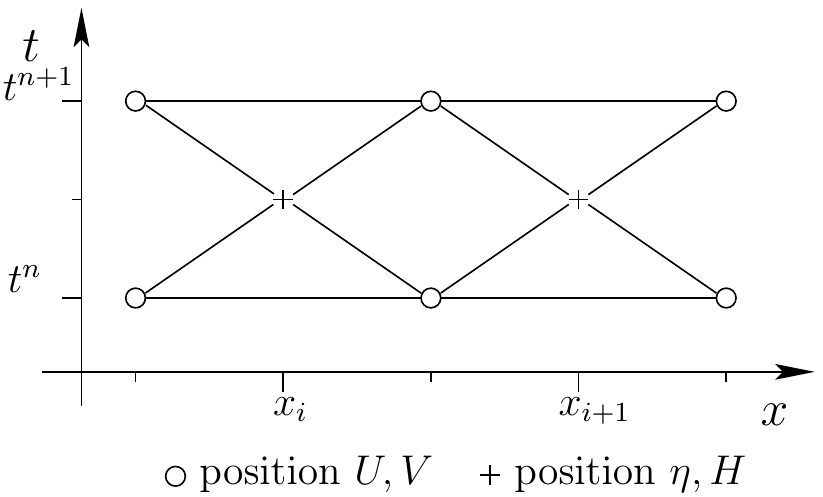}
  \caption{\small\label{fig:stag_stencil} Stencils of volume-flux update
  for finite-difference schemes. Left: first order scheme, not
  symmetric with respect to time. Right: time symmetry recovered.}
\end{figure}

The original B-grid scheme of \cite{Gjevik_etal_2002} is based on a
staggering of unknowns, where the volume-fluxes $U$ and $V$ are
approximated in the grid points $(ih, jh, n\Delta t)$ and the surface
elevation $\eta$ is approximated in shifted grid points $((i+\oh)h,
(j+\oh)h, (n+\oh)\Delta t)$ as shown in Figure~\ref{fig:bgrid}.  To
ease the presentation we introduce the following standard differencing
and averaging operators:
\begin{equation*}
\begin{array}{l@{\quad\quad\quad}r}
  \delta_x()_{\cdot,\cdot} = 
  \frac{1}{h}\left[()_{\cdot+\oh,\cdot} -()_{\cdot-\oh,\cdot}\right], & 
  \mu_x()_{\cdot,\cdot} = 
  \oh\left[()_{\cdot+\oh,\cdot} +()_{\cdot-\oh,\cdot}\right],\\[0.3cm]
  \delta_y()_{\cdot,\cdot} = 
  \frac{1}{h}\left[()_{\cdot,\cdot+\oh} -()_{\cdot,\cdot-\oh}\right], & 
  \mu_y()_{\cdot,\cdot} = 
  \oh\left[()_{\cdot,\cdot+\oh} +()_{\cdot,\cdot-\oh}\right].
\end{array}
\end{equation*}
Note that it is implied that the result of these operations is shifted
by $\frac{h}{2}$ relative to the argument. For simplicity in notation,
we omit $i$ and $j$ indices in the following scheme.  The meaning
should be clear from the aforementioned shift and the position of the
point-wise approximations.
%, i.e.,
%\begin{align*}
%  &\eta^{n+\oh}_{i+\oh,j+\oh} \approx\eta((i+\oh)h,(j+\oh)h, (n+\oh)\Delta t), \\
%  &U^n_{ij} \approx U(ih,jh, n\Delta t),\\
%  &V^n_{ij} \approx V(ih,jh, n\Delta t)
%\end{align*}
For instance, the approximation of $\partial_x (U^2/H)$ in the point
$(ih, jh)$ is given by
\begin{align*}
  \left(\delta_x\frac{\big(\mu_x U\big)^2}{\mu_y H}\right)_{ij}&=
   \frac{1}{2h}\left(\frac{\big(U_{i+1,j}+U_{ij}\big)^2}{H_{i+\oh, j-\oh}+H_{i+\oh, j+\oh}}
  -\frac{\big(U_{ij}+U_{i-1,j}\big)^2}{H_{i-\oh, j-\oh}+H_{i-\oh, j+\oh}}\right)
\end{align*}
With this notation, we can  write the finite-difference scheme for \eqref{eq:sw}
as
\begin{align}
  \label{eq:fd1_eta}
  \eta^{n+\oh} &= \eta^{n-\oh}
  -\Delta t[\delta_x \mu_y U^n +\delta_y \mu_x V^n]\\
  \label{eq:fd1_U}
  U^{n+1} &= U^n
  -\Delta t\left[\delta_x \frac{\big(\mu_x U^n\big)^2}{\mu_y H^{n+\oh}} 
    +\delta_y \frac{\big(\mu_y U^n\big)\;\big(\mu_y V^n\big)}{\mu_x H^{n+\oh}}
    +\big(g\mu_x\mu_yH^{n+\oh}\big) \delta_x\mu_y\eta^{n+\oh}-fV^n\right],\\
  \label{eq:fd1_V}
  V^{n+1} &= V^n
  -\Delta t\left[\delta_x \frac{\big(\mu_x U^n\big)\;\big(\mu_x V^n\big)}{\mu_y H^{n+\oh}}
    +\delta_y \frac{\big(\mu_y V^n\big)^2}{\mu_x H^{n+\oh}} 
    +\big(g\mu_x\mu_yH^{n+\oh}\big)\delta_y\mu_x\eta^{n+\oh}+fU^n\right],
\end{align}
where $\Delta t$ is the time step and $H^{n+\oh} = \eta^{n+\oh}-z$.
The use of a B-grid yields a quite compact second-order discretisation
of the flux and source terms.  We would like to point out that, due to
the central differencing, the scheme
\eqref{eq:fd1_eta}--\eqref{eq:fd1_V} is second order accurate in
space, but not in time.  This can be seen most easily from left part
of Figure \ref{fig:stag_stencil}, which shows that the stencil of the
volume-flux update is not symmetric with respect to time.

In order to correct the asymmetry, we introduce the following shorthands:
\begin{align*}
\Phi_x(\eta) &:=\big(g\mu_x\mu_yH\big)\delta_x\mu_y\eta,\\
\Phi_y(\eta) &:=\big(g\mu_x\mu_yH\big)\delta_y\mu_x\eta,\\
\Omega(U,V,\eta) &:= \delta_x \frac{\big(\mu_x U\big)^2}{\mu_y H} 
  +\delta_y \frac{\big(\mu_y U\big)\;\big(\mu_y V\big)}{\mu_x H} - fV,\\
\Psi(U,V,\eta) &:= \delta_x \frac{\big(\mu_x U\big)\;\big(\mu_x V\big)}{\mu_y H}
  +\delta_y \frac{\big(\mu_y V\big)^2}{\mu_x H} + fU.
\end{align*}
where $H=H(\eta,z)$. With this notation \eqref{eq:fd1_U} -- \eqref{eq:fd1_V} read
\begin{align}
  \label{eq:fd11_U}
  U^{n+1} &= U^n -\Delta t\left[\Omega(U^n,V^n,\eta^{n+\oh}) + \Phi_x(\eta^{n+\oh}) \right],\\
  \label{eq:fd11_V}
  V^{n+1} &= V^n -\Delta t\left[\Psi(U^n,V^n,\eta^{n+\oh}) + \Phi_y(\eta^{n+\oh}) \right].
\end{align}

Let us now introduce the correction which assures second order
accuracy in time.  For this we denote the volume-flux update in
\eqref{eq:fd11_U} -- \eqref{eq:fd11_V} by $(U,V)^{n+1}_*$ and centre
the terms $\Omega$ and $\Psi$, i.e the flux differences and the
coriolis term, with respect to time. This gives
\begin{align}
\label{eq:fd2_U}
U^{n+1} &= U^{n+1}_* + \frac{\Delta t}2 \left[\Omega(U^n,V^n,\eta^{n+\oh}) - \Omega(U^{n+1}_*,V^{n+1}_*,\eta^{n+\oh}) \right],\\
\label{eq:fd2_V}
V^{n+1} &= V^{n+1}_* + \frac{\Delta t}2 \left[\Psi(U^n,V^n,\eta^{n+\oh}) - \Psi(U^{n+1}_*,V^{n+1}_*,\eta^{n+\oh}) \right],
\end{align}
which is the symmetric stencil shown in right part of Figure \ref{fig:stag_stencil}. 
%\begin{center}
%\end{center}
In Tables \ref{table:convergence2D_enlag} and \ref{table:convergence2D_improved_enlag} one
can clearly observe the gain in accuracy.

An elementary calculation shows that both the first-order version and
the second-order version of this scheme are well-balanced for the
stationary state of water at rest $U=V=0$ and $\eta - z = Const$.  For
smooth solutions driven by inflow boundary conditions, both scheme
yields quite sharp results with moderate numerical diffusion.  For
non-smooth solutions, both versions of the scheme experience
instabilities in the form of oscillations.
%
%================================================================================================
\subsection{The High-Order Finite-Volume Scheme}\label{sec:FV_Scheme}
%================================================================================================
%
\begin{figure}
  \centering
  \includegraphics{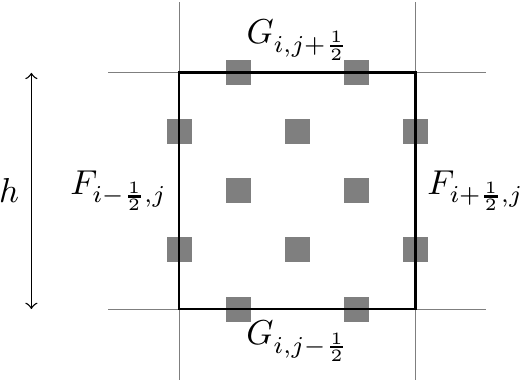}
  \caption{\small\label{fig:cvolumes} A finite-volume, with
    corresponding fluxes and integration points on cell interfaces and
    in the interior. The integration points in the interior are needed
    for the well-balanced integration rule for the source source
    term.}
\end{figure}
To   simplify  the   presentation of the finite-volume scheme
somewhat,   we  rewrite
\eqref{eq:sw} as
\begin{align}\label{eq:sw2} 
  Q_t + F(Q)_x+ G(Q)_y = B(Q, Q_x, Q_y)+C(Q),
\end{align}
where subscript denotes differentiation, $Q=[\eta, U, V]^T$ is the
vector of unknown functions and $F$ and $G$ are vector-valued
functions.  The source terms $B=-gH [0, z_x, z_y]^T$ and $C=f[0,V,
-U]^T$ are the geometrical source accounting for variable bottom and
the Coriolis force term, respectively.

The discretisation of the homogeneous part of \eqref{eq:sw2} is
straightforward.  As opposed to the scheme in the previous section,
the high order finite-volume scheme is based on computing cell
averages over grid cells $I_{ij}=\left[\left(i-\oh\right)h,
\left(i+\oh\right)h\right]\times
\left[\left(j-\oh\right)h,\left(j+\oh\right)h\right]$ of a uniform
Cartesian mesh.  In each grid cell we approximate
\begin{align*}
  \bar {Q}^n_{ij} = \int\!\!\!\int_{I_{ij}}   Q(x,y,t^n) dxdy.
\end{align*}
To compute the evolution of the cell averages $\eta_{ij}$, $U_{ij}$
and $V_{ij}$, we must approximate the flux over the cell interfaces.
This is accomplished by \emph{reconstructing} the solution within each
grid cell using (nonlinear) interpolation.  To obtain a stable and
accurate solution, the reconstruction procedure must ensure that no
spurious oscillations are introduced, even when the solution is
non-smooth. This may be done, for example, by the WENO technique (see
\cite{Shu_1998} and the references therein). The reconstruction yields
one-sided approximations of point values at cell interfaces as well as
in the interior of the cell, see Figure~\ref{fig:cvolumes}.  Thus, at
each cell interface we obtain two one-sided approximations.  To
compute consistent interface fluxes we integrate Roe's approximate
flux function $\hat{F}(Q^+, Q^-,n)$ \cite{Roe_1981} over each cell
interface.  In computations, these integrals are approximated using
Gaussian quadrature. For the $x$-direction, $F_{i+\oh, j}$ is computed
as,
\begin{align*}
  F_{i+\oh,j} &= \int_{(j-\oh)h}^{(j+\oh)h} \hat{F}\Big(Q\big((i+\oh)h, y\big)^-,Q\big((i+\oh)h, y\big)^+,n_x\Big)dy, \\
  &\approx h\sum_\alpha \omega_\alpha \hat{F}\Big(Q\big((i+\oh)h, y_\alpha\big)^-,Q\big((i+\oh)h, y_\alpha\big)^+,n_x\Big),
\end{align*}
where $y_\alpha$ and $\omega_\alpha$ are the quadrature points and
weights. The fluxes in the $y$-direction $G_{i,j+\oh}$ are computed in
the same manner.  To complete the spatial discretisation we must
compute the averaged source terms $B_{ij}$ and $C_{ij}$ over each grid
cell.  The final evolution of the cell averages is computed by solving
the semi-discrete equation
\begin{align}
  \frac{d}{dt}\bar{Q}_{ij}  = - (F_{i+\oh,j}-F_{i-\oh},j)/h  -  
  (G_{i,j+\oh}-G_{i,j-\oh})/h  +  B_{ij} + C_{ij},
\end{align}
using a standard fourth-order Runge-Kutta scheme.  The details of the
source term discretisation are given in
Section~\ref{sec:FV_SourceTerms}, and the boundary conditions are
discussed in Section~\ref{sec:FV_BoundaryConditions}.
%
%===============================================================================================
\subsection{Well-balanced Finite-Volume Treatment of the Source Terms}\label{sec:FV_SourceTerms}
%===============================================================================================
%
The shallow water system has stationary solutions where source terms
and flux terms are in equilibrium.  If we discretise the fluxes and
source terms naively, this may lead to spurious oscillations near
equilibria. To accurately resolve small perturbations of such
equilibria, we must ensure that the discrete fluxes and sources
exactly balance at equilibrium.  Ideally, the truncation error of the
scheme should vanish at equilibrium states.  This is often called a
well-balanced treatment of source terms.

Flat, stationary water of variable depth is an equilibrium where it is
possible to construct such a discretisation.  In
\cite{Noelle_etal_2006}, an arbitrary-order well-balanced
discretisation of the geometrical source term was constructed.
A well-balanced second-order discretisation of the geometrical source
term is extended to arbitrary order of accuracy using an asymptotic
expansion.  In a single grid cell in one spatial dimension, the
integral of the source term can be approximated by the fourth-order
rule
\begin{align}
  \label{eq:wbint}
  S & = \frac{g}{6}\left[4(\eta_{l}+\eta_{c})(z_{l}-z_{c})
    +4(\eta_{c}+\eta_{r})(z_{c}-z_{r})-(\eta_{l}+\eta_{r})(z_{l}-z_{r})\right]
\end{align}
where $\eta_l$, $\eta_c$ and $\eta_r$ are reconstructed point-values
in the left-, centre- and right endpoints of the cell.  The
reconstruction of the central point is an additional cost associated
with this source term.  This integration rule can be extended to two
spatial dimensions using Gaussian quadrature. In the equation for
$x$-volume-flux, \eqref{eq:wbint} is applied in the $x$-direction and the
Gaussian rule in the $y$-direction.  For the $y$-volume-flux, the order
is reversed.  The integration points used for the source terms are
shown in Figure~\ref{fig:cvolumes}.

The Coriolis term is approximated by $C_{ij}=f[0, V_{ij}, -U_{ij}]^T$.
As in \cite{Bouchut_etal_2004}, this is well-balanced for grid-aligned
geostrophic jets.
%
%
%
%=========================================================================================
\section{Treatment of boundary conditions}\label{sec:BoundaryConditions}
%=========================================================================================
%
%
%
For the experiments presented in this paper we need three types of
boundary conditions, reflective, outflow and inflow. These are
presented below. While reflective boundaries are rather
straightforward, out- and inflow conditions have to be translated
carefully from the B-grid finite-difference setting to the
finite-volume setting. Moreover, for the finite-volume scheme we
discovered a subtle perturbation introduced by a no-slip inflow
boundary condition.  In Sections~\ref{sec:FV_InflowBC} and
\ref{sec:FV_MomBalInf_BC}, we introduce free-slip, Neumann-type
boundary conditions which give smoother inflow.
%
%=========================================================================================
\subsection{Finite-Difference Boundary Conditions}\label{sec:BGridBoundaryConditions}
%=========================================================================================
%
Reflective boundary condition are treated with one-sided differences
and normal volume-flux equal zero. At {\em outflow} boundaries the
normal volume-flux $V$ is defined by $\eta\sqrt{gH}$ and the transverse
volume-flux $U$ is set to zero. This boundary condition is called Flather
condition in mechanics, and it coincides e.g. with the first order
absorbing boundary condition given by Eng\-quist and Majda
\cite{Engquist_Majda_1977}. The normal velocity $v$ on the {\em
inflow} boundary is given by a time dependent velocity profile
function $v_{\textnormal{jet}}(x,y,t)$ (see \eqref{eq:vjet1} and
\eqref{eq:vjet2}) and the tangential velocity $u$ is set to zero
(no-slip). To compute the volume-flux on the inflow boundary the height
is extrapolated from the interior. These boundary conditions are
straightforward to implement in the finite-difference scheme.
%
%===========================================================================================
\subsection{Finite-Volume Boundary Conditions}\label{sec:FV_BoundaryConditions}
%===========================================================================================
%

%===========================================================================================
\subsubsection{Reflective boundary conditions}\label{sec:FV_ReflectiveBC}
%===========================================================================================
To treat reflective boundary conditions we are using ghost-cells and
solve the Riemann problem on the reflective boundary, where the ghost
cell contains the same data as the interior cell, but with reflected
volume-fluxes.

%================================================================================
\subsubsection{Absorbing Outflow Boundary Condition}\label{sec:FV_AbsorbingBC}
%================================================================================
Here we adopt a technique developed by Eng\-quist and Majda
\cite{Engquist_Majda_1977} to derive a so-called first order absorbing
boundary conditions for the outflow boundary. In particular, we follow
Kr\"oner's adaptation \cite{Kroener_1991} of the Eng\-quist-Majda
absorbing boundary condition, who has computed the relevant
decomposition into normal and tangential waves for the linearised
Euler equations.

For the shelf flows which we would like to compute, there are two
relevant cases.  Due to the very large speed of long gravity waves in
the ocean $c_0=\sqrt{gH} \gg v$ we are always in a subcritical flow
and either one or two characteristics are leaving the domain. Carrying
over the results of \cite{Engquist_Majda_1977,Kroener_1991} to the
linearised shallow water equations, we obtain the following analytical
boundary conditions:

If two characteristics are leaving the domain, the first order
absorbing boundary conditions specifies the normal volume-flux
$V=\eta\sqrt{g H}$ at the open boundary.

In the case of one outgoing characteristic we obtain as before
$V=\eta\sqrt{g H}$ for the normal volume-flux. In addition, we obtain a
no-slip condition for the tangential volume-flux, namely $U=0$.

Now we translate these analytical boundary conditions to obtain data
for the Riemann solver at the absorbing boundary. Let $\xi$ be the
outward pointing coordinate normal to the boundary and let $Q_L$ be
the approximation at the interior point $\xi=0-$. We want to determine
$Q_R$ at $\xi=0+$ such that an appropriate discretisation of the
absorbing boundary condition is fulfilled at the boundary. We
linearise the system around the interior state $Q_L$.  In the
subcritical case, there are two possibilities (see
Figure~\ref{fig:riemann_inflow}):
\renewcommand{\labelenumi}{(\alph{enumi})}
\begin{enumerate}
\item $v\ge0$ at  $\xi=0$ and  only $\lambda_1<0$  or,
\item $v<0$  at $\xi=0$  and both $\lambda_1<0$ and $\lambda_2 <0$.
\end{enumerate}
In case (a) we need to specify $Q_A$, since only one characteristic
enters the domain, and we can only specify one condition on the
absorbing boundary. For this we choose
\begin{align}
  \label{eq:absorbing1}
  V_{A}=V_{\eta}:=\eta_{L}\sqrt{gH_{L}},
\end{align}
corresponding to the absorbing condition. This condition is equivalent
to the radiation condition of Flather \cite{Flather_1976} which is
used by Gjevik et. al. \cite{Gjevik_etal_2002}. 

In case (b), we want to prescribe the state $Q_B$, since two
characteristics enter the domain and we have to specify two conditions
at the boundary.  In this case the first order absorbing boundary
conditions are 
\begin{align}  
\label{eq:absorbing2a}
  V_{B} &=V_{\eta} \\
 \label{eq:absorbing2b}
   U_B &= 0.
\end{align}
In both cases, this yields a well-posed problem.  In case (a), the
states $Q_L$ and $Q_A$ are separated by a simple wave,
\begin{align*}
  Q_A :=  Q_L + \alpha_1 r_1.
\end{align*}
These are three equations for the unknowns $\alpha_1,H_A$ and $U_A$.
An elementary calculations gives
\begin{align}  
  \label{eq:H_A}
  H_{A} &= H_L (1+\beta_L) \\
  \label{eq:U_A} 
  U_A &= U_L (1+\beta_L)
\end{align}
where
\begin{align}
\beta_L:=  \frac{V_{\eta}-V_L}{V_L - H_L \sqrt{g{H_L}}}.
\end{align}
Dividing \eqref{eq:U_A} by \eqref{eq:H_A} one obtains that
\begin{align}
u_A = u_L,
\end{align}
so there is no jump in tangential velocity!

In case (b), the states $Q_L$ and
$Q_B$ are connected by two waves separated by the intermediate state
$Q_A$,
\begin{align*}
Q_B &:=  Q_L + \alpha_1 r_1 + \alpha_2 r_2
\end{align*}
From \eqref{eq:absorbing2a}--\eqref{eq:absorbing2b}, we have
$V_B=V_{\eta}:=\eta_L\sqrt{gH_L}$ and $U_B=0$. A straightforward
computation yields
\begin{align}
H_B &= H_L (1+\beta_L).
\end{align}
The numerical flux at the northern boundary is simply
\begin{align}
G_{i,j_{max}+\oh} := \left\{\begin{array}{cc}
G(Q_A) & \quad\mathrm{for\;case\;(a)} \\
G(Q_B) & \quad\mathrm{for\;case\;(b)}
\end{array}\right.,
\end{align}
where $(i,j_{max})$ is the index of the northern cells adjacent to the
boundary.

%===========================================================================================
\subsubsection{Free-slip inflow boundary conditions}\label{sec:FV_InflowBC}
%===========================================================================================
For the finite-volume scheme, we have implemented two types of inflow
boundary conditions. The first is the no-slip inflow boundary
condition which we have already described for the finite-difference
scheme. As is documented in Table~\ref{table:conv2D_nonsmooth_FV} for
our finite-volume scheme, this leads to a loss of accuracy even for
smooth incoming jets. Since the water depth is computed from values
downstream, small inaccuracies in the specification of the boundary
condition can lead to large numerical errors or even instabilities.

This may be explained as follows: fixing the tangential velocity at
the boundary to be zero leads to a jump in tangential velocity when
the internal flow develops vortices near the boundary. This admittedly
small discontinuity can cause loss of accuracy, and must be removed to
get the expected rate of convergence.

Therefore we would like to propose a second type of boundary
condition, which we call {\em free-slip}.  We will show that this
leads to smoother solutions.
\iffalse
%
\begin{figure}
  \centering
  \includegraphics[width=0.49\linewidth]{figures/FV_NoSlip_Linear_vel_E0_T120h}
  \includegraphics[width=0.49\linewidth]{figures/FV_FreeSlip_Linear_vel_E0_T120h}
  \caption{\small\label{fig:FreeSlipNoSlip}
    Velocity plots at time
  120h, on the left side no-slip inflow is used and on right side free-slip.
  Eddies develop close to the boundary when free-slip inflow is used,
  details of this experiment will be given in Section~{sec:Comparison}.}
\end{figure}
%
\fi
To analyse the possible inflow boundary conditions, we linearise the
system at the inflow boundary around a state
$\hat{Q}:=(\hat{H},\hat{U},\hat{V})^T$, with $U=Hu$, $V=Hv$ and
$H=\eta-z$. Assuming $z_\xi =0$, we obtain
\begin{align*}
  Q_t + \hat{A}Q_{\xi}=0, 
\end{align*}
where $\xi$ is the coordinate normal to the inflow boundary, $\hat{A}$
is the Jacobian of the flux function in the $\xi$-direction, $V$ is
the volume-flux in the $\xi$-direction and $U$ is the volume-flux parallel
to the boundary.

\begin{figure}
  \centering
  \includegraphics[width=0.3\linewidth]{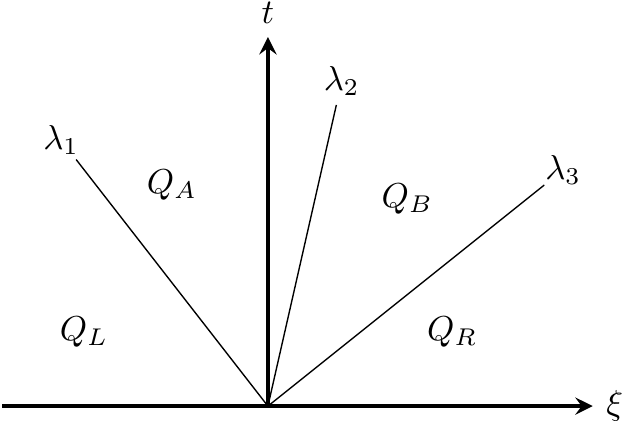} 
  \hspace{2cm}
  \includegraphics[width=0.3\linewidth]{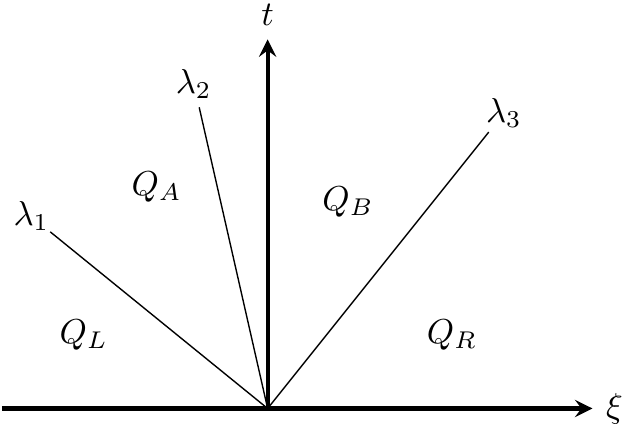} 
  \caption{\small\label{fig:riemann_inflow}Solution of linear  Riemann problem at inflow
    boundary, $Q_L$ exterior, $Q_R$ interior domain.}
\end{figure}

The general Riemann solution consists of four states $Q_L$, $Q_A$,
$Q_B$ and $Q_R$ as shown in Figure~\ref{fig:riemann_inflow}. They are
connected by three waves travelling with speeds
\begin{align*}
  \lambda_1  &= \hat{v}  -  \sqrt{g\hat{H}},\qquad  \lambda_2 =  \hat{v},\qquad
  \lambda_3 = \hat{v} + \sqrt{g\hat{H}},
\end{align*}
where $\hat{v}$ is the component of the velocity in the
$\xi$-direction. The corresponding eigenvectors are denoted by
$r_1,r_2,r_3$.

For our boundary value problem, we have $\hat v > 0$, since we assume
that we are at an inflow boundary. Typical velocities $\hat v$ will
not exceed one meter per second. But the typical speed of long gravity
waves $\sqrt{gH}$ will be of the order of 30 meters per second to 140
meters per second for water depths of 100 to 2000 meters. Thus, the
inflow velocity $\hat v$ is much smaller that the speed of long
gravity waves, and we have subcritical flow. Therefore, the
eigenvalues satisfy
\begin{align*}
  \lambda_1\ll 0 \leq \lambda_2\ll\lambda_3.
\end{align*}
As a result, the numerical boundary data we are looking for are given
by $Q_A$. This state $Q_A$ will be connected by waves of the second
and the third families to the state $Q_R$,
\begin{align}
  \label{eq:linearRiemannProblem}
  \begin{bmatrix}
    H_A \\ U_A \\ V_A 
  \end{bmatrix} -
  \begin{bmatrix}
    H_R \\ U_R \\ V_R
  \end{bmatrix} 
  =
  \alpha_2
  \begin{bmatrix}
    0 \\ 1 \\ 0
  \end{bmatrix}+
  \alpha_3
  \begin{bmatrix}
    1 \\ \frac{\hat{U}}{\hat{H}} \\ \frac{\hat{V}}{\hat{H}}+\sqrt{g\hat{H}} 
  \end{bmatrix}.
\end{align}
These are three equation for the five unknowns $H_A, U_A, V_A,
\alpha_1$ and $\alpha_2$.  To obtain a uniquely solvable system we
need to specify two of the unknowns. This corresponds to the fact that
exactly two characteristic are entering the domain.
Oliger and Sundstr\o m  \cite{Oliger_Sundstroem_1978} showed that the  initial boundary
value problem \eqref{eq:sw2} is well posed under the boundary condition
\begin{align} \label{eq:inflowbc}
  \frac{\partial u}{\partial \xi}=0\quad\mbox{ and }\quad
  v=v_{jet},\qquad \mbox{ when }\xi=0.
\end{align}
We translate condition \eqref{eq:inflowbc} to our inflow Riemann
problem by requiring that
\begin{align} \label{eq:inflowbcA}
  u_A=u_R \quad\mbox{ and }\quad
  v_A=v_{jet}.
\end{align}
Now we have only three unknowns left. Plugging $u_A$ and $v_A$ into
\eqref{eq:linearRiemannProblem} gives
\begin{align}
  \label{eq:linearRiemannProblem2}
  \begin{bmatrix}
    H_A \\ H_A\;u_A \\ H_A\;v_{\textnormal jet}  
  \end{bmatrix} -
  \begin{bmatrix}
    H_R \\ U_R \\ V_R
  \end{bmatrix} 
  =
  \alpha_2
  \begin{bmatrix}
    0 \\ 1 \\ 0
  \end{bmatrix}+
  \alpha_3
  \begin{bmatrix}
    1 \\ \frac{\hat{U}}{\hat{H}} \\ \frac{\hat{V}}{\hat{H}}+\sqrt{g\hat{H}} 
  \end{bmatrix},
\end{align}
which yields
\begin{align*}
  \alpha_2 &= H_A\;u_A - H_R\;u_R -\alpha_3\frac{\hat{U}}{\hat{H}},  \qquad
  \alpha_3 = \frac{V_A-V_R}{\frac{\hat{V}}{\hat{H}}+\sqrt{g \hat{H}}}.
\end{align*}
Choosing $\hat{Q}  = Q_R$, this  leads to the  following formula for  the state
$Q_A$ at the boundary,
\begin{align}\label{eq:FV_freeslipBC}
  H_A &= \frac{H_R\sqrt{gH_R}}{v_R+\sqrt{gH_R}-v_{jet}},\qquad U_A = H_Au_R,\qquad
  V_A = H_Av_{jet}.
\end{align}
The numerical flux at the boundary is simply
\begin{align}
G_{i,\oh} := G(Q_A).
\end{align}
%========================================================================================
\subsubsection{Balanced Inflow Boundary Condition}\label{sec:FV_MomBalInf_BC}
%========================================================================================
In Section~\ref{sec:jetexample} we will apply another variant of the
jet inflow boundary condition. In order to motivate it, let us
consider once more the free-slip boundary condition derived in the
previous section. As can be seen from Figure~\ref{fig:riemann_inflow}
and equation \eqref{eq:linearRiemannProblem}, the jet inflow data were
assigned to the intermediate state $Q_A$ via $v_A=v_{jet}$. Then the
state $Q_A$ was connected to the inner state by two waves. This
defined $u_A=u_R$ and $H_A$ implicitly. We would like to point out
that the third wave, a long gravity wave leaving the domain, is
effectively suppressed, and no wave can leave the domain. Indeed in
Section~\ref{sec:BarostrophicJetSetupE2} we show that this may lead to
an increase of the overall water height.

Now we modify the boundary condition to include the outgoing wave. For
this, we apply the jet inflow condition to the outer state $Q_L$
instead of $Q_A$. Since we now have three waves to connect the inner
state $Q_R$ with the jet, there is one more degree of freedom.  We
determine this by the following reasoning: at the jet, we already know
the normal velocity $v(x,0-,t)=v_L(x,t)=v_{jet}(x,t)$. By the
free/slip condition, we also know the tangential velocity
$u(x,0-,t)=u_L(x,t)=u_R(x,t)$. It remains to determine
$H(x,0-,t)=H_L(x,t)$.  Now we request that these values
$(H,u,v)(x,0-,t)$ are compatible with the shallow water equations.  To
determine $H_L(x,t)$ it is sufficient to use the balance of tangential
volume-flux. Taking into account that $u_y(x,0,t)=0$ by the the free-slip
condition, we obtain
\begin{align}\label{eq:MomentumBalancePV}
u_t + uu_x =  - g (H + z)_x + f v,
\end{align}
or
\begin{align}\label{eq:MomentumBalancePV2}
H(x) + z(x) = H(x_0) + z(x_0) - \frac{1}{2g}(u(x)^2-u(x_0)^2)
-\frac{1}{g}\int^x_{x_0}(u_t - f v).
\end{align}
Note that the geostrophic balance 
\begin{align}\label{eq:GeostrophicBalance}
\eta(x)-\eta(x_0) = \frac{f}{g}\int^x_{x_0}v
\end{align}
is a special case of the volume-flux balance
Equation~\eqref{eq:MomentumBalancePV2} when $u \equiv 0$ (remember that
$\eta = H + z$).

For the cells $(i,1)$ at the southern boundary, the in-flowing flux is
given by the Riemann solver $\hat F$ via
\begin{align}
G_{i,\oh} := \hat F(Q_L,Q_R,n_y),
\end{align}
where $n_y=(0,1)$ is the inward unit vector normal to the southern boundary.
Note that $y_{\oh}$ is the position of the boundary edge. 

In Section~\ref{sec:jetexample} we will see that this boundary
condition is transparent, i.e. it admits both in- and outflow.
%
%
%
%===============================================================================
\section{Comparison of the Schemes}\label{sec:Comparison}
%===============================================================================
%
%
%
In this section we present comparisons of the staggered scheme and the
high-order finite-volume scheme on different challenging test
problems.
%
%===============================================================================
\subsection{Order of Accuracy}\label{sec:TestingOrder}
%===============================================================================
%
To compute the numerical order of accuracy of the finite-volume scheme
we use a slight modification of an experiment of Xing and Shu
(\cite{Xing_Shu_2005}, see also \cite{Lukacova_etal_2006}.  On the
unit square $[0,1]\times [0,1]$ the bottom topography, initial surface
elevation, and initial volume-flux are given by the smooth functions
\begin{align*}
  z(x,y)      &= \sin(2 \pi x)+\cos(2 \pi y),\\
  \eta(x,y,0) &= 10+e^{\sin(2 \pi x)}\cos(2\pi y),\\
  U(x,y,0)   &= \sin(\cos(2 \pi x))\;\sin(2 \pi y),\\
  V(x,y,0)   &= \cos(2\pi x)\;\cos(\sin(2 \pi y)).
\end{align*}
We compute the solution up to time $T=0.05$ with CFL-number $0.5$.
The physical parameters are $g=9.812$ and $f=10.0$. 
The reference solution is computed with the finite-volume
scheme on a grid with $1600 \times 1600$ cells.
%
% first order enlag %
%
\begin{table}  
\begin{tabular}{c|c c|c c|c c}\hline 
 $N$    & \multicolumn{2}{c|}{$H$}&\multicolumn{2}{c|}{$U$}&\multicolumn{2}{c}{$V$} \\
        & $L^1$-error & rate & $L^1$-error & rate& $L^1$-error & rate \\ \hline
  25 & 4.56E-02 &      & 1.70E-01 &      & 4.37E-01 &      \\ 
  50 & 1.69E-02 & 1.43 & 7.44E-02 & 1.19 & 1.76E-01 & 1.31 \\ 
 100 & 7.19E-03 & 1.23 & 3.33E-02 & 1.16 & 7.58E-02 & 1.22 \\ 
 200 & 3.35E-03 & 1.10 & 1.58E-02 & 1.08 & 3.48E-02 & 1.12 \\ 
 400 & 1.63E-03 & 1.04 & 7.68E-03 & 1.04 & 1.66E-02 & 1.07 \\ 
 800 & 8.02E-04 & 1.02 & 3.80E-03 & 1.02 & 8.12E-03 & 1.03 \\ 
\hline
\end{tabular}\\[0.3cm]
 \caption{\small\label{table:convergence2D_enlag} The $L^1$-errors and
convergence rates for each of the components in the convergence test
of Section~\ref{sec:TestingOrder}, computed with the first-order
finite-difference scheme of Section~\ref{sec:TheBgridScheme}.  The
reference solution is computed with the high-order finite-volume
scheme on a $1600 \times 1600$ grid.}
\end{table}
%
% second order enlag %
%
\begin{table}  
\begin{tabular}{c|c c|c c|c c}\hline 
 $N$ & \multicolumn{2}{c|}{$H$}&\multicolumn{2}{c|}{$U$}&\multicolumn{2}{c}{$V$} \\
     & $L^1$-error & rate & $L^1$-error & rate& $L^1$-error & rate \\ \hline
  25  & 3.27E-02 &      & 1.19E-01 &      & 2.41E-01 &      \\
  50  & 8.45E-03 & 1.96 & 3.30E-02 & 1.85 & 6.27E-02 & 1.94 \\
 100  & 2.10E-03 & 2.01 & 8.45E-03 & 1.96 & 1.60E-02 & 1.97 \\
 200  & 5.26E-04 & 2.00 & 2.12E-03 & 2.00 & 4.01E-03 & 2.00 \\
 400  & 1.32E-04 & 2.00 & 5.31E-04 & 2.00 & 1.01E-03 & 2.00 \\
 800  & 3.29E-05 & 2.00 & 1.33E-04 & 2.00 & 2.52E-04 & 2.00 \\
\hline 
\end{tabular}\\[0.3cm]
 \caption{\small\label{table:convergence2D_improved_enlag} The
$L^1$-errors and convergence rates for each of the components in the
convergence test of Section~\ref{sec:TestingOrder}, computed with the
improved second-order finite-difference scheme of
Section~\ref{sec:TheBgridScheme}.  The reference solution is computed
with the high-order finite-volume scheme on a $1600 \times 1600$
grid.}
\end{table}
%
% fourth order finite-volume %
% 
\begin{table}  
\begin{tabular}{c|c c|c c|c c}\hline 
   $N$  & \multicolumn{2}{c|}{$H$}&\multicolumn{2}{c|}{$U$}&\multicolumn{2}{c}{$V$} \\
        & $L^1$-error&   rate   & $L^1$-error&   rate   & $L^1$-error& rate     \\ \hline
 25     & 6.70E-03 &      & 2.06E-02 &      & 5.34E-02 &      \\  
 50     & 8.46E-04 & 2.99 & 1.60E-03 & 3.69 & 7.30E-03 & 2.87 \\ 
 100    & 6.84E-05 & 3.63 & 9.19E-05 & 4.13 & 5.57E-04 & 3.71 \\ 
 200    & 3.06E-06 & 4.48 & 3.70E-06 & 4.64 & 2.48E-05 & 4.49 \\ 
 400    & 1.10E-07 & 4.79 & 1.32E-07 & 4.81 & 9.03E-07 & 4.78 \\ 
 800    & 3.66E-09 & 4.91 & 4.38E-09 & 4.91 & 3.04E-08 & 4.90 \\ 
 \hline 
\end{tabular}\\[0.3cm]
\caption{\small\label{table:convergence2D} The $L^1$-errors and convergence rate
for each  component in  the convergence test  of Section~\ref{sec:TestingOrder},
computed      with     the      high-order      finite-volume     scheme      of
Section~\ref{sec:FV_Scheme}.   Each grid  is  $N\times N$  and the  reference
solution is computed on a $1600 \times 1600$ grid.}
\end{table}

According to the discussion in Section~\ref{sec:TheBgridScheme}, we
expect the original finite-difference scheme
\eqref{eq:fd1_eta}--\eqref{eq:fd1_V} to be first order accurate. This
is confirmed by the results in Table~\ref{table:convergence2D_enlag}.
The improved enlag scheme \eqref{eq:fd2_U} and \eqref{eq:fd2_V} is
indeed second-order accurate, see
Table~\ref{table:convergence2D_improved_enlag}.

For the finite-volume scheme we expect fourth-order accuracy (indeed
the Runge-Kutta scheme for time integration, the Gaussian rules for
integrating the numerical fluxes and the cell centred source term are
all formally fourth-order accurate, and the spatial WENO
reconstruction procedure is even fifth-order accurate).

Table~\ref{table:convergence2D} reports the $L^1$-errors together with
convergence rates for the finite-volume scheme.  For this test case,
we get the expected fourth-order accuracy (in fact almost fifth-order)
in all components.
%
%
%
%===============================================================================
\subsection{Large Eddies in a Doubly Periodic Domain.}\label{sec:LargeEddiesDoublePeriodicDomain}
%===============================================================================
%
%
%
\begin{figure}
  \centering
  \includegraphics[width=0.6\linewidth]{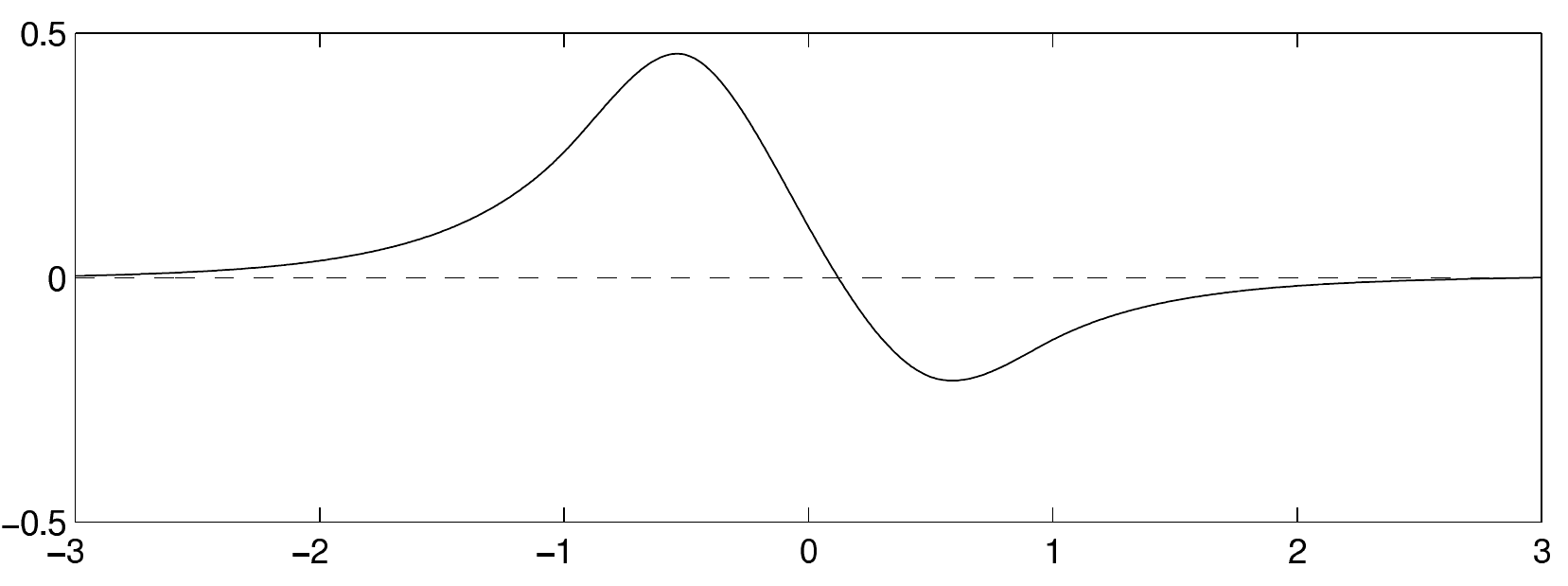}\\
  \includegraphics[width=0.6\linewidth]{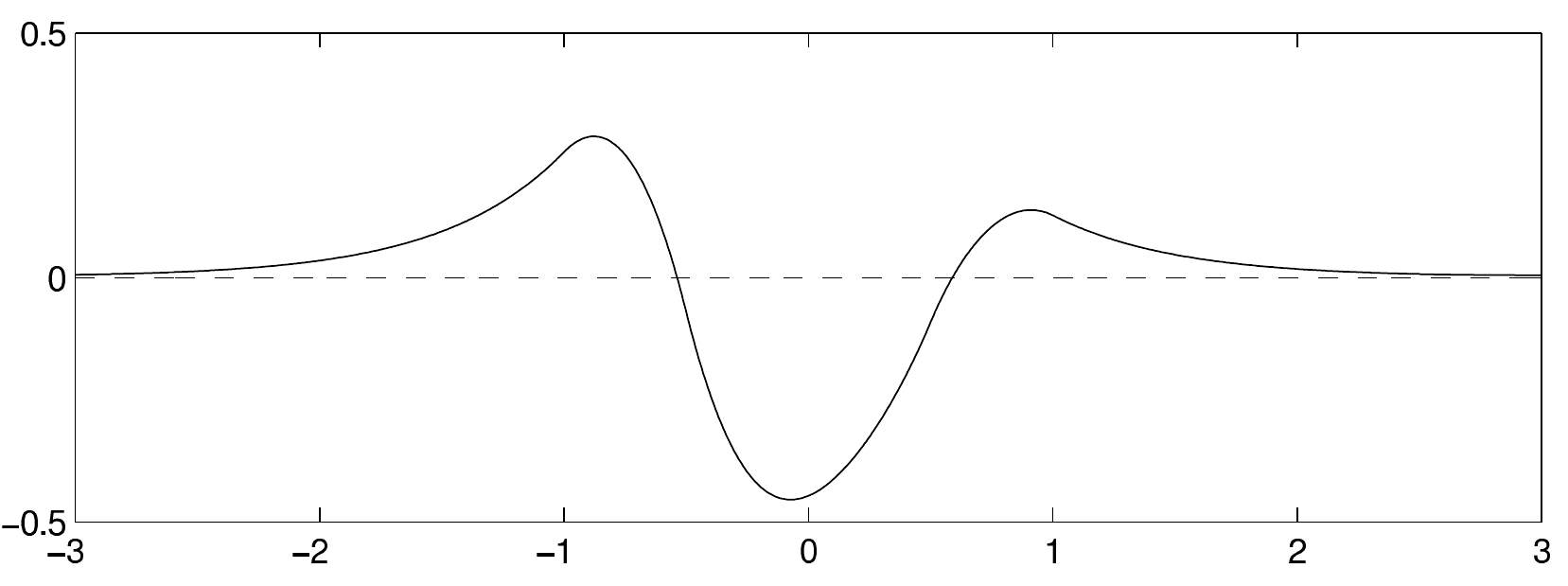} \\
  \includegraphics[width=0.6\linewidth]{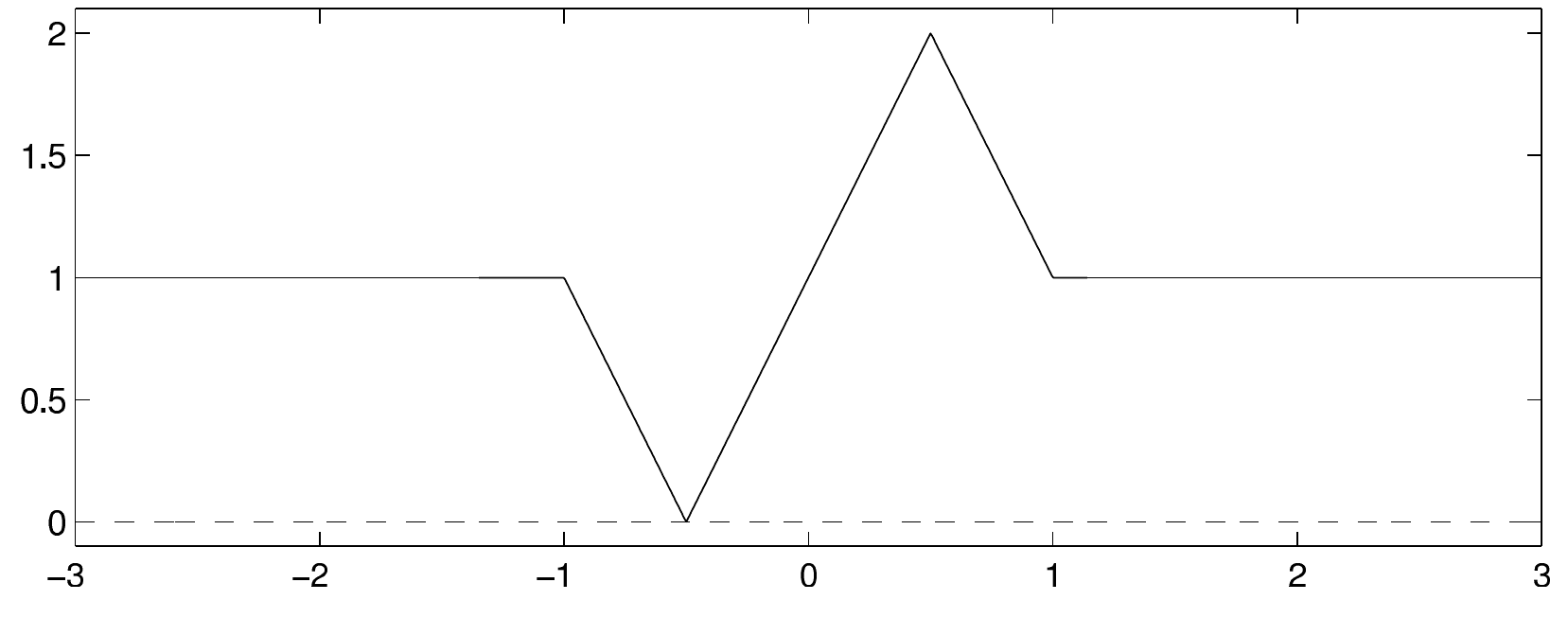} \\
  \caption{\small\label{fig:balance}The figure shows a cross-section
    of the initial data of Example~2: (top) surface elevation,
    (middle) $x$-component of the velocity field, and (bottom)
    potential vorticity.}
\end{figure}

To illustrate visually the difference in performance of the
finite-difference and the finite-volume schemes, we compute the
evolution of potential vorticity (PV) in a very hard test-case taken
from \cite{Dritschel_etal}.  The PV is a conserved quantity that is
advected with the flow and is a good test of the effect of the
numerical diffusion on complex smooth solutions of the rotating
shallow water equations.  Since this test-case is doubly periodic,
which can easily be implemented in both schemes, the comparison does
not involve the complications of boundary conditions.

Consider a doubly periodic domain $(-\pi, \pi)^2$ with flat bottom
topography.  Let $\mathbf{u} := (u,v)$ be the velocity field.  The
potential vorticity is given by
\begin{align}\label{eq:pot_vort}
  q := \frac{\nabla\times\mathbf{u} +f}{H}.
\end{align}
Assume  that  the  flow  is geostrophically  balanced initially, i.e.,  that  the
gravitational forces  exactly balance the Coriolis force. Using \eqref{eq:sw},
this balance can be written as 
\begin{align}
  \label{eq:balance}
  g\nabla H+f\mathbf{u}^\perp&=0.
\end{align}
If the potential vorticity is known throughout the domain, the balance
condition \eqref{eq:balance} specifies the state of the shallow water
system completely.  At this state, the surface elevation solves the
following equation,
\begin{align}
  H_{xx}+ H_{yy} + \frac{fq}{g}H &= \frac{f^2}{g},
\end{align}
with evanescent boundary conditions.

In this example, we use the initial potential vorticity of \cite{Dritschel_etal},
\begin{align*}
  q(x,y,0)&=
  \begin{cases}
    \bar{q} + Q\mbox{ sign}(\hat{y})(a-\big| |\hat{y}|-a\big|),& |\hat{y}|<2a,\\
    \bar{q}, & \mbox{otherwise},
  \end{cases},\\
  \hat{y} &=y+c_m\sin mx + c_n\sin nx, 
\end{align*}
where $\bar{q}$ is the mean potential vorticity, $\bar{q} \pm Qa$ is
the maximum/minimum of the potential vorticity, and $2a$ is the width
of the jet.  As in \cite{Dritschel_etal}, we use the scalings
$\bar{h}=1$, $L_R^2=gH/f^2=0.25$, $a=0.5$, $\bar{h}Q/f=2$ and
$f=4\pi$, with one unit time corresponding to one day.  The parameters
of the perturbation are $m=2$, $n=3$, $c_2 =-0.1$ and $c_3=0.1$.

By solving the balance condition \eqref{eq:balance}, the potential
vorticity field yields a balanced double jet flow.  Cross-sections of
the initial surface elevation, velocity field and potential vorticity
are shown in Figure~\ref{fig:balance}.  The actual solution of the
balance condition is computed using a simple central finite-difference
scheme on a $512\times 512$ grid.

\begin{figure}
  \centering
  \begin{tabular}{cc}
    Finite-difference scheme, day 4 & Finite-volume scheme, day 4 \\
    \includegraphics[width=0.495\linewidth]{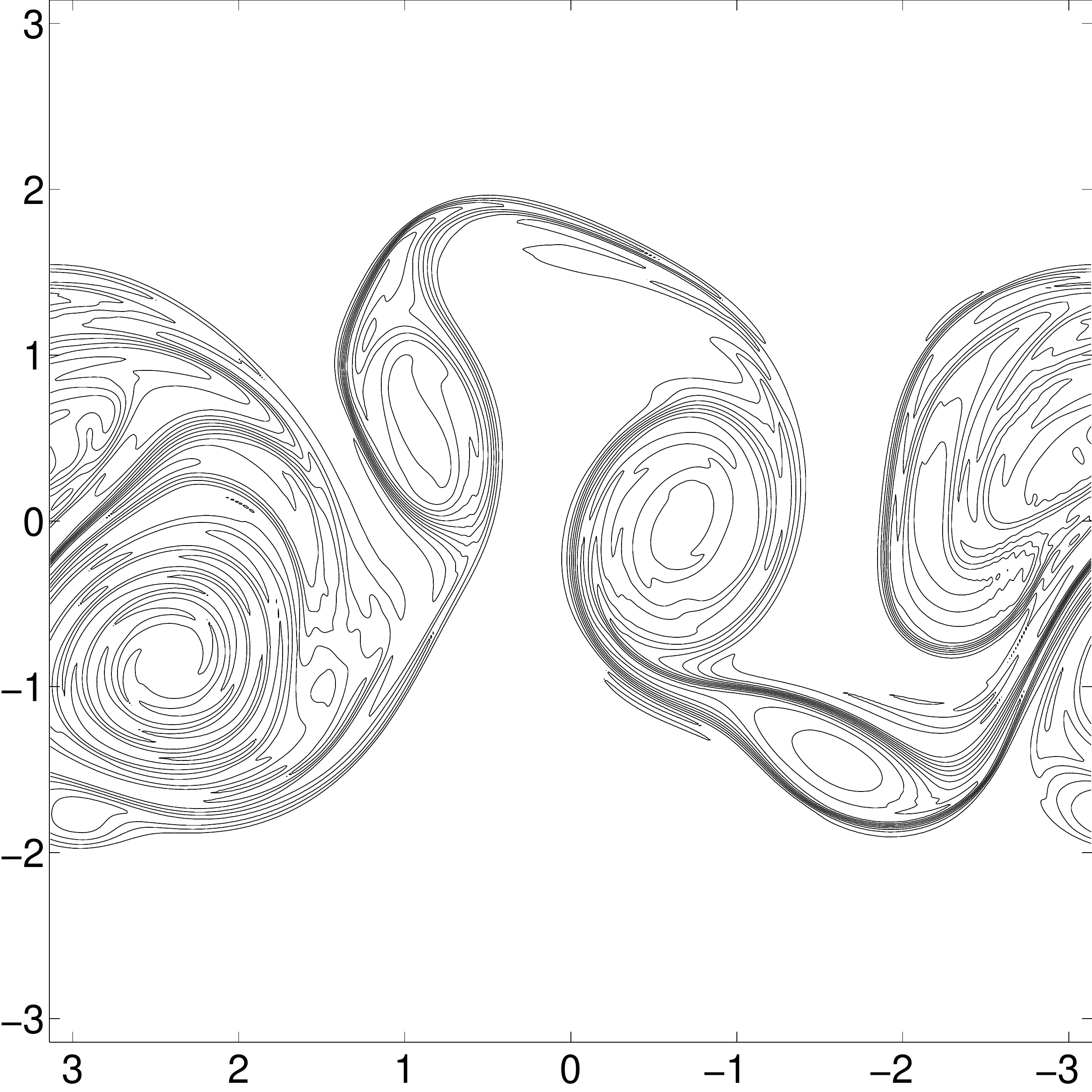}&
    \includegraphics[width=0.495\linewidth]{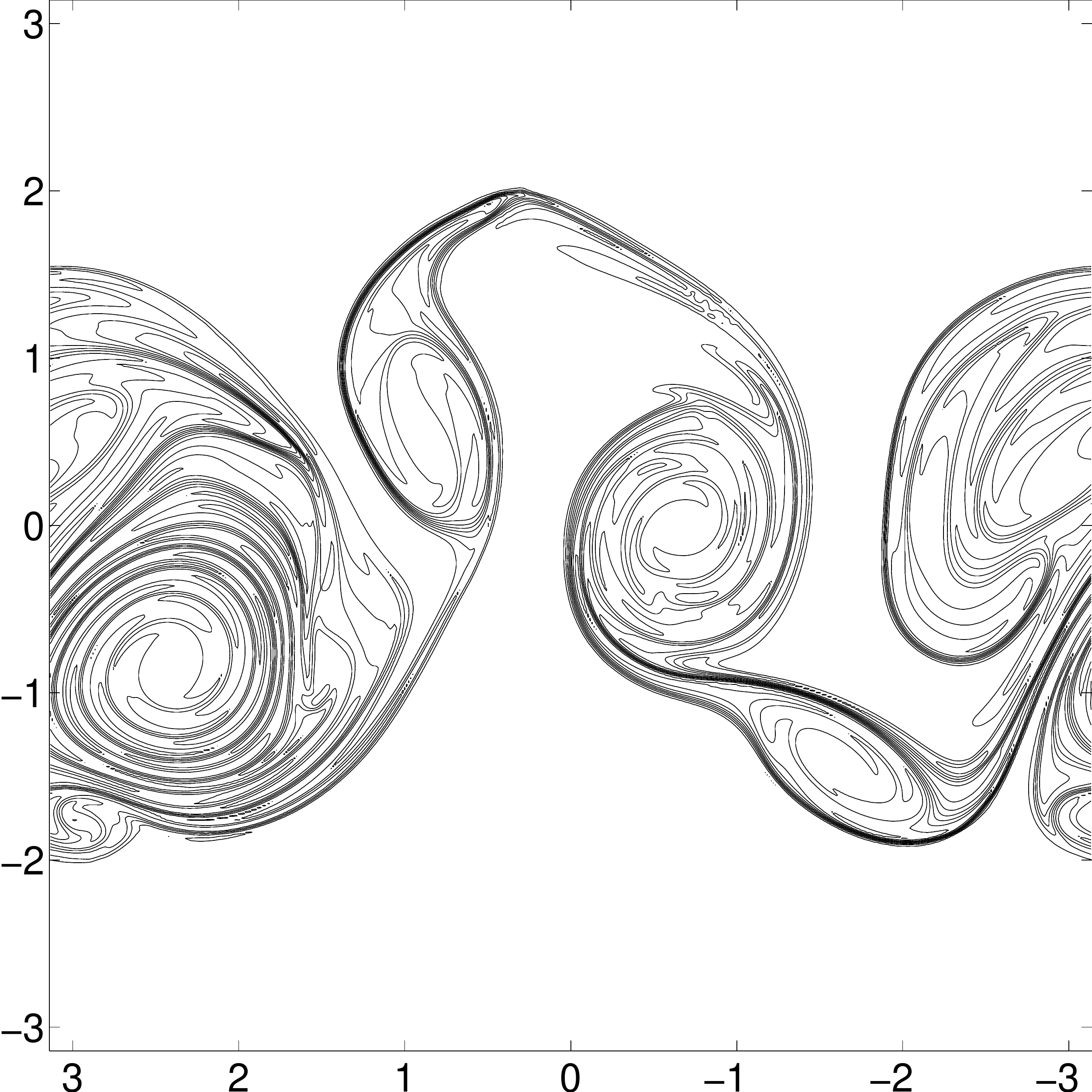}\\
    Finite-difference scheme, day 8 & Finite-volume scheme, day 8 \\
    \includegraphics[width=0.495\linewidth]{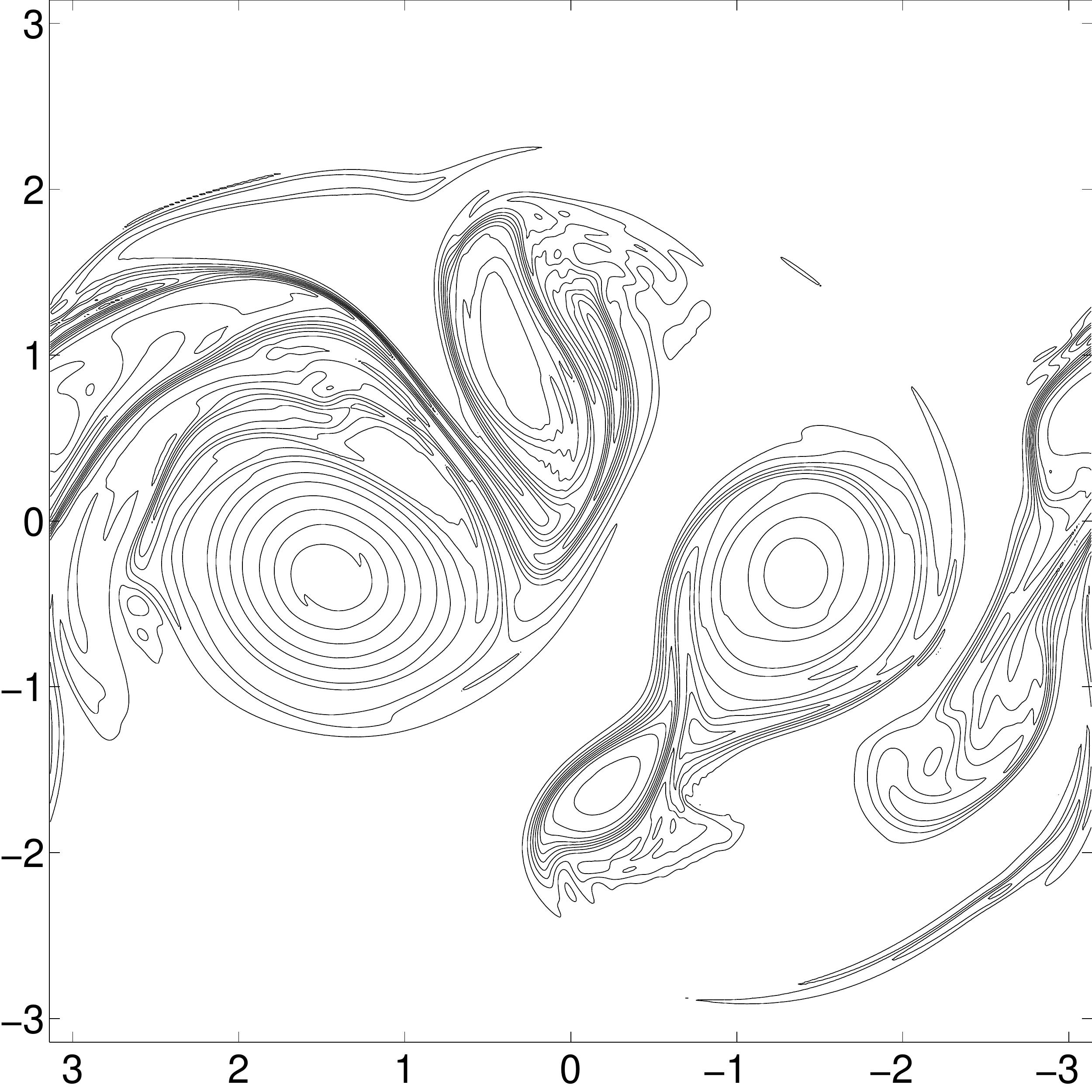}&
    \includegraphics[width=0.495\linewidth]{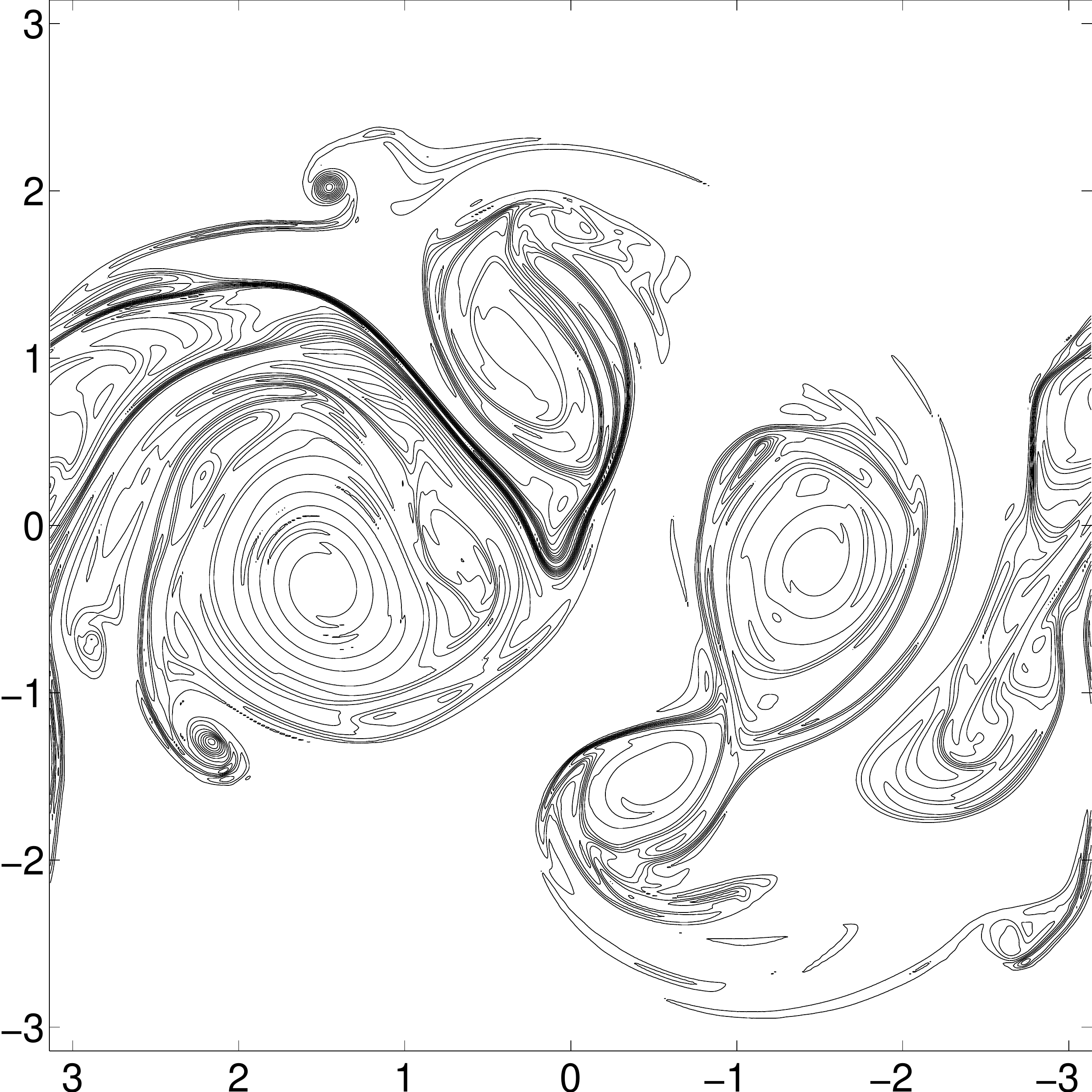}\\
  \end{tabular}
  \caption{\small\label{fig:vortices} Contour lines of the potential
    vorticity at day $4$ (up) and day $8$ (down) computed with the finite-volume
    scheme and finite-difference scheme, both on a $512\times 512$
    grid.  The vortical pattern closely resembles the results reported
    in \cite{Dritschel_etal}, although clearly with more numerical
    diffusion.}
\end{figure}

The time evolution of these seemingly simple initial data quickly
produces large complex vortical structures with many smaller vortex
filaments tearing off. In Figure~\ref{fig:vortices} we have plotted
the potential vorticity at two times for the finite-volume scheme and
the finite-difference scheme, respectively. Both schemes seem to
produce the same coarse scale vortices, but in addition the high order
finite-volume scheme also resolves several small scale vortices.  See
\cite{Dritschel_etal} for a comparison with a Semi-Lagrangian contour
advection algorithm.

%
%===============================================================================
\subsection {Convergence Test for a Barotropic Jet Problem}
\label{sec:jetConvTest}
%===============================================================================
%
%\marginlabel{This is new}
In this and the following two examples we study barotropic jets.
Here we show that the no-slip boundary condition described in
Section~\ref{sec:FV_InflowBC} yields the expected loss of convergence
rates. We also show that the free-slip boundary condition
gives high-order convergence rates.

As in \cite{Gjevik_etal_2002,Thiem_etal_2006} the water is initially
at rest.  Then the jet is started smoothly across the southern
boundary (see Figure~\ref{fig:domain}) with velocity
\begin{align}\label{eq:vjet1}
  v_{\textnormal jet}(x, t) = V_{\textnormal max}\exp\left(- \left(\frac{2(x-L_B)}{B}\right)^2\right)
  \gamma \left(\frac{t}{2000}\right),
\end{align}
where the growth function $\gamma$ is given by 
\begin{align}\label{eq:growth_function}
\gamma(\tau):=
\left \{
\begin{array}{ll}
 70\;\tau^9 -315\;\tau^8 + 540\;\tau^7 - 420\;\tau^6 + 126\;\tau^5, &  \textnormal{ if } \tau \leq 1,\\
 1, & \textnormal{ else. } 
\end{array}
\right ..
\end{align}
The centre of the jet $L_B = 100 km$ and the width is $B = 50 km$.
The maximum velocity is $V_{\textnormal max} = 0.04 \frac{m}{s}$.
The full strength of the jet is reached after $2000 s$.  We compute on
the domain $\Omega = [0,300 km]\times[0,300 km]$ with smooth
bottom topography given by
\begin{align*} 
  Q_1 &= 0.5 \; (D_O - D_S),\\
  Q_2 &= 0.5 \; (D_O + D_S),\\ 
  z(x,y) &= -Q_1 \; \textnormal{tanh}\left(\frac{x-X_O}{X_S} \right)- Q_2,
\end{align*}
where $D_S = 400 m$, $D_O = 1000 m$, $X_S = 40 km$, and $X_O = 120
km$. Initially, the water in $\Omega$ is at rest, so
$\eta(x,y) \equiv 0$. The boundary conditions in the $x$-direction are
reflective (east-west), while the boundary condition at $y=0$ (south)
is an inflow condition.

The northern boundary condition (at $y=300 km$) is transparent. We use
the radiation condition \cite{Flather_1976} for the finite-differences
and the absorbing boundary condition \cite{Engquist_Majda_1977} for
the finite-volume scheme. For this particular example, the two
conditions are equivalent. The acceleration of gravity $g=9.81
\frac{m}{s^2}$ and the Coriolis parameter $f=1.2\times 10^{-4}
s^{-1}$.

In
Tables~\ref{table:conv2D_nonsmooth_FD} and \ref{table:conv2D_nonsmooth_FV}
we have computed the rates of convergence of the two schemes at the
final time $T=3000 s$. At this time the jet has flooded a large part
of the domain and the west-going wave is partially reflected at the
boundary.  The reference solution was computed using the finite-volume
scheme on a $1600\times 1600$ grid. As predicted in
Section~\ref{sec:FV_InflowBC} both schemes do not converge with the
high rates obtained in Example~\ref{sec:TestingOrder}.  This is in
accordance with the discussion in
\cite{Oliger_Sundstroem_1978,Bennett_Kloeden_1981,Edwards_etal_1983},
which predicts that a no-slip inflow boundary condition will result in
a loss of smoothness in the whole domain, for computations on very fine grids.

To obtain a smoother solution we apply the free-slip boundary condition
developed in Section~\ref{sec:FV_InflowBC} to the finite-volume scheme.
Our boundary condition for the finite-difference scheme is that the
tangential volume-flux should be continuous, which is realised by
\begin{align}\label{eq:FD_smoothInflowBC}
U_{i,\textnormal{inflow}}=U_{i,\textnormal{inflow}+1}.
\end{align}
As shown in Tables~\ref{table:conv2D_smooth_FD} and
\ref{table:conv2D_smooth_FV} these boundary conditions recover the
expected higher orders of convergence, especially for the
finite-volume scheme. Once more the reference solution was computed
using the finite-volume scheme on a $1600\times 1600$ grid.
%
% Smooth Jet FD no slip BC %
%
\begin{table}  
\begin{tabular}{c|c c|c c|c c}\hline 
 $N$    & \multicolumn{2}{c|}{$\eta$}&\multicolumn{2}{c|}{$U$}&\multicolumn{2}{c}{$V$} \\
        & $L^1$-error & rate & $L^1$-error & rate & $L^1$-error & rate \\ 
\hline
  50 & 6.10E07 &      & 7.89E09 &      & 4.84E09 &      \\
 100 & 2.90E07 & 1.07 & 4.26E09 & 0.89 & 2.42E09 & 1.00 \\
 200 & 1.38E07 & 1.07 & 2.30E09 & 0.95 & 1.20E09 & 1.02 \\
 400 & 6.31E06 & 1.13 & 1.11E09 & 0.99 & 5.85E08 & 1.04 \\
 800 & 2.76E06 & 1.20 & 5.51E08 & 1.01 & 2.83E08 & 1.05 \\
\hline
\end{tabular}\\[0.3cm]
 \caption{\small\label{table:conv2D_nonsmooth_FD} The $L^1$-errors and
convergence rates for each of the components in the convergence test of
Section~\ref{sec:jetConvTest}. The solutions were computed with the second-order
finite-difference scheme with no-slip inflow conditions.  The
reference solution was computed with the high-order finite-volume
scheme on a $1600 \times 1600$ grid.}
\end{table}
%
% Smooth Jet FV no slip BC %
%
\begin{table}  
\begin{tabular}{c|c c|c c|c c}\hline 
 $N$    & \multicolumn{2}{c|}{$\eta$}&\multicolumn{2}{c|}{$U$}&\multicolumn{2}{c}{$V$} \\
        & $L^1$-error & rate & $L^1$-error & rate & $L^1$-error & rate \\ 
\hline
  50& 3.19E06 &      & 4.57E08 &      & 4.47E08 &      \\
 100& 2.11E05 & 3.92 & 8.21E07 & 2.48 & 5.70E07 & 2.97 \\
 200& 1.51E04 & 3.80 & 6.97E06 & 3.56 & 4.36E06 & 3.71 \\
 400& 5.10E03 & 1.57 & 4.19E06 & 0.73 & 1.88E06 & 1.21 \\
 800& 2.67E03 & 0.93 & 3.72E06 & 0.17 & 1.02E06 & 0.89 \\
\hline
\end{tabular}\\[0.3cm]
 \caption{\small\label{table:conv2D_nonsmooth_FV} The $L^1$-errors and
convergence rates for each of the components in the convergence test
of Section~\ref{sec:jetConvTest}.  The solutions were computed with
the finite-volume scheme with no-slip boundary conditions.}
\end{table}
%
% Smooth Jet FD free slip BC % 
%
\begin{table}  
\begin{tabular}{c|c c|c c|c c}\hline 
 $N$ & \multicolumn{2}{c|}{$\eta$}&\multicolumn{2}{c|}{$U$}&\multicolumn{2}{c}{$V$}\\ 
     & $L^1$-error & rate & $L^1$-error & rate & $L^1$-error & rate \\ 
  \hline 
   50 & 2.14E07 &      & 2.34E09 &      & 2.91E09 &      \\
  100 & 1.07E07 & 0.99 & 1.10E09 & 1.09 & 1.38E09 & 1.07 \\
  200 & 5.20E06 & 1.05 & 4.91E08 & 1.16 & 6.35E08 & 1.13 \\
  400 & 2.25E06 & 1.21 & 2.02E08 & 1.28 & 2.67E08 & 1.25 \\
  800 & 7.25E05 & 1.64 & 6.42E07 & 1.65 & 8.53E07 & 1.65 \\
  \hline
\end{tabular}\\[0.3cm]  
\caption{\small\label{table:conv2D_smooth_FD} The $L^1$-errors and
convergence rates for each of the components in the convergence test
of Section~\ref{sec:jetConvTest}.  The solutions were computed with
the finite-difference scheme with free-slip inflow conditions. The
reference solution was computed with the high-order finite-volume
scheme on a $1600 \times 1600$ grid.}
\end{table}  
%
% Smooth Jet FV free slip BC %
%
\begin{table}  
\begin{tabular}{c|c c|c c|c c}\hline 
 $N$    & \multicolumn{2}{c|}{$\eta$}&\multicolumn{2}{c|}{$U$}&\multicolumn{2}{c}{$V$} \\
        & $L^1$-error & rate & $L^1$-error & rate & $L^1$-error & rate \\ 
\hline
  50 & 3.20E06 &      & 4.64E08 &      & 4.57E08 &      \\ 
 100 & 2.20E05 & 3.86 & 9.08E07 & 2.35 & 6.43E07 & 2.83 \\ 
 200 & 1.66E04 & 3.72 & 1.37E07 & 2.73 & 7.85E06 & 3.04 \\ 
 400 & 1.32E03 & 3.65 & 1.59E06 & 3.11 & 8.23E05 & 3.25 \\ 
 800 & 1.02E02 & 3.69 & 1.45E05 & 3.45 & 7.16E04 & 3.52 \\ 
\hline
\end{tabular}\\[0.3cm]  
\caption{\small\label{table:conv2D_smooth_FV} The $L^1$-errors and
convergence rates for each of the components in the convergence test
of Section~\ref{sec:jetConvTest}.  The solutions were computed with
the finite-volume scheme with free-slip inflow boundary conditions.}
\end{table}  
%
%
%===============================================================================
\subsection {Development of Eddies in Shelf Slope Area due to a Barotropic Jet}
\label{sec:jetexample}
%===============================================================================
%
%===============================================================================
\subsubsection{Ormen Lange Shelf Experiment~I} \label{sec:BarostrophicJetSetupE1}
\begin{figure}
\centering
\includegraphics{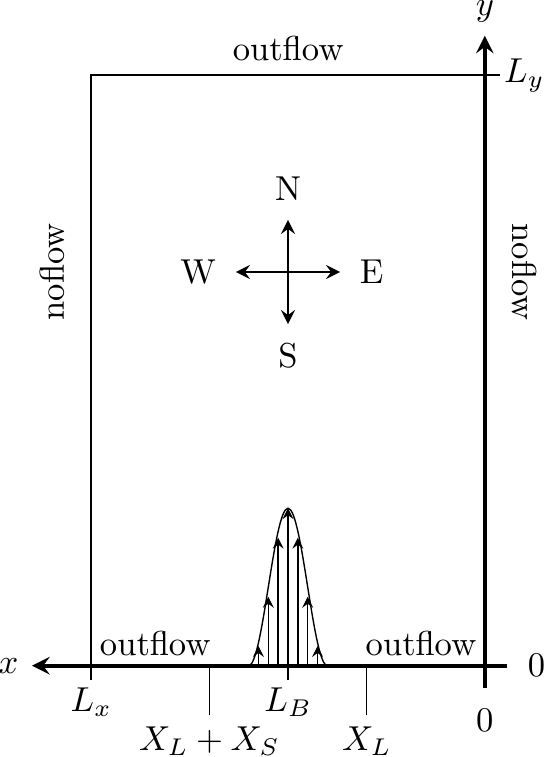}
\caption{\label{fig:domain} The  computational domain of the  example in Section
\ref{sec:jetexample}.  The types of boundary conditions used are indicated.}
\end{figure}

%===============================================================================
%
In \cite{Thiem_etal_2006}, Thiem et al. used a numerical model based
on the first-order finite-difference scheme of
Section~\ref{sec:TheBgridScheme} to study the impact of the shelf
geometry upon along-shelf currents. The setup is taken from the Ormen
Lange gas field off the western Norwegian coast. The shelf width in
this model is constant with a depth profile given by
\begin{align*}
 z(x,y) = \begin{cases}
   - D_S,                                      & x\leq X_L,\\
   - D_O + (D_O-D_S)\left(\frac{X_L+X_S-x}{X_S}\right)^2,& X_L \leq x \leq X_L+X_S, \\
   - D_O,&\mbox{otherwise},
   \end{cases}
\end{align*}
where $D_S=250 m$, $D_O=1600 m$, $X_L=100 km$, and $X_S=90 km$.  The
domain is $[0,L_x]\times[0,L_y]$, where $L_x=300km$ and $L_y=600km$.
Initially, the surface elevation $\eta=0 m$ and the water is initially
at rest.  The boundary conditions in the $x$-direction (east, west)
are reflective. On part of the southern boundary ($y=0$, $|x-L_B|\leq
B$ with $B=10 km$, $L_B = 115 km$) we prescribe an in-flowing jet with
velocity
\begin{align}\label{eq:vjet2}  
  v_{jet}(x,t) = V_{\textnormal max}\exp\left(-\left(\frac{2(x-L_B)}{B}\right)^2\right)(1-\exp(-\sigma t))
\end{align}
where $V_{\textnormal max} = 0.4 m/s$, the jet growth factor
$\sigma=2.3148\times 10^{-5}$, acceleration of gravity $g=9.81 m/s^2$
and Coriolis parameter $f=1.2\times 10^{-4} s^{-1}$.  For the rest of
the southern as well as for the northern boundary we prescribe an
absorbing radiation condition, see Figure~\ref{fig:domain}.

In Subsection~\ref{sec:jetConvTest} we have compared the no-slip
inflow boundary condition with the more accurate free-slip inflow
boundary condition for a smooth jet. Now we will study these boundary
conditions for the more realistic Ormen~Lange~setup described
above.

The first computation uses the finite-difference scheme with no-slip
inflow boundary condition as described in
Section~\ref{sec:BGridBoundaryConditions}. The second computation is
done by the finite-volume scheme. The no-slip inflow boundary
condition is the same as the free-slip boundary condition
\eqref{eq:FV_freeslipBC}, except that we set the tangential velocity
$u_A$ to zero. At the outflow boundary, we use conditions
\eqref{eq:absorbing2a} and \eqref{eq:absorbing2b}.

The final computation, again by the finite-volume scheme, uses the
free-slip inflow boundary condition \eqref{eq:FV_freeslipBC}.

%Various other absorbing boundary conditions were also tried, but the
%authors did not observe significant differences in the solution.  The
%boundary condition that has the greatest impact on the solution is the
%inflow condition.

\begin{figure}
  \centering
  \includegraphics[width=0.32\linewidth]{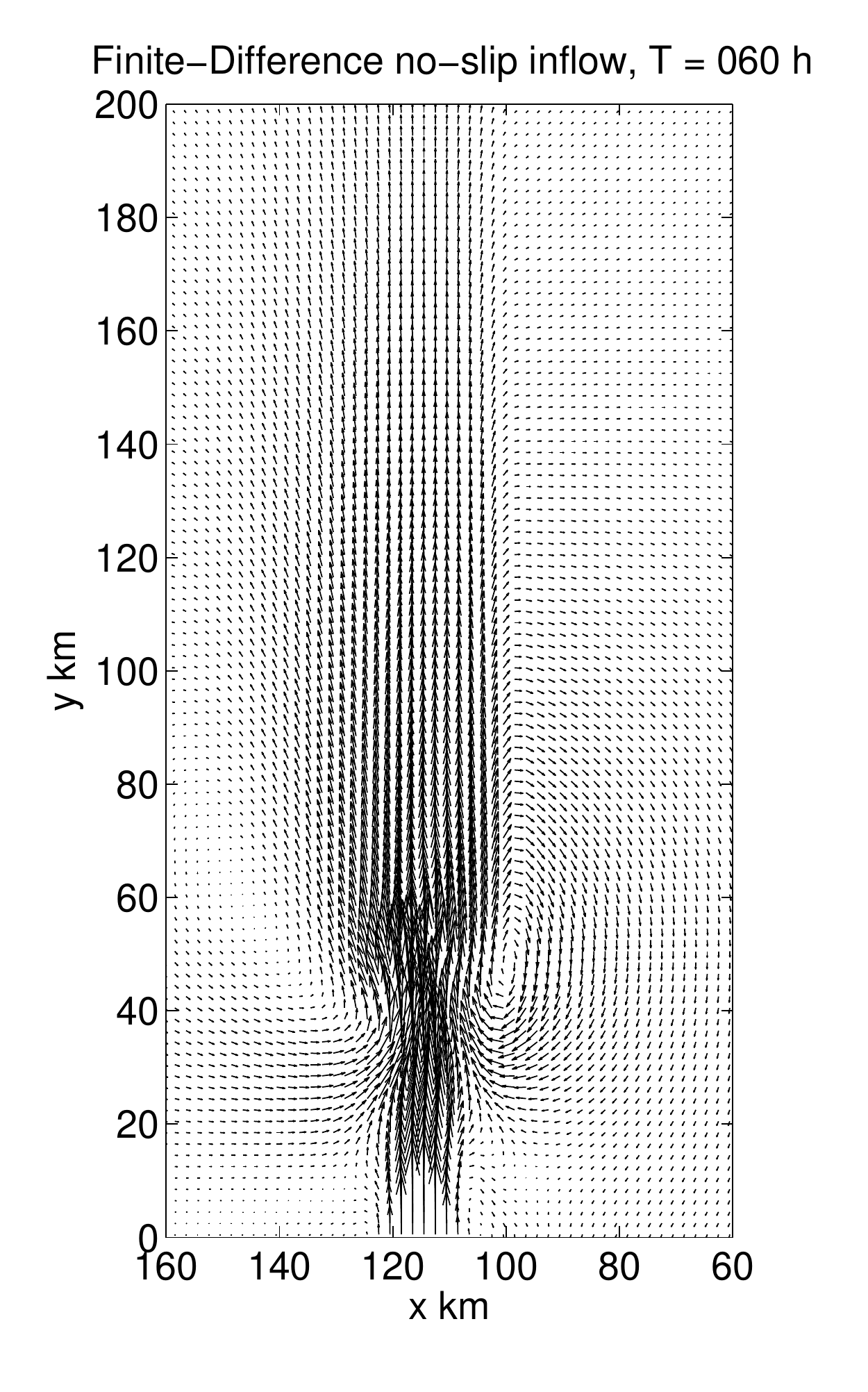}
  \includegraphics[width=0.32\linewidth]{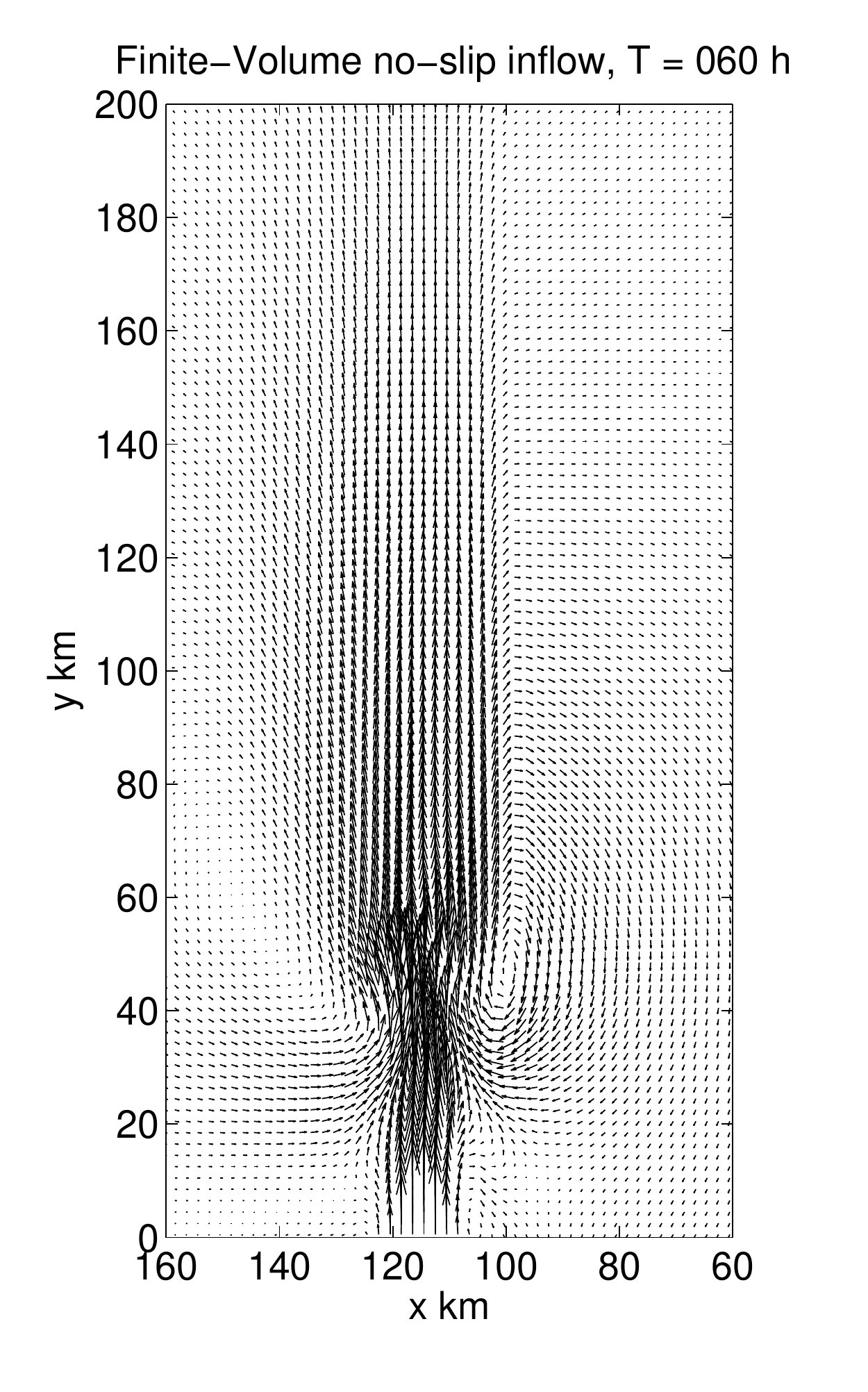}
  \includegraphics[width=0.32\linewidth]{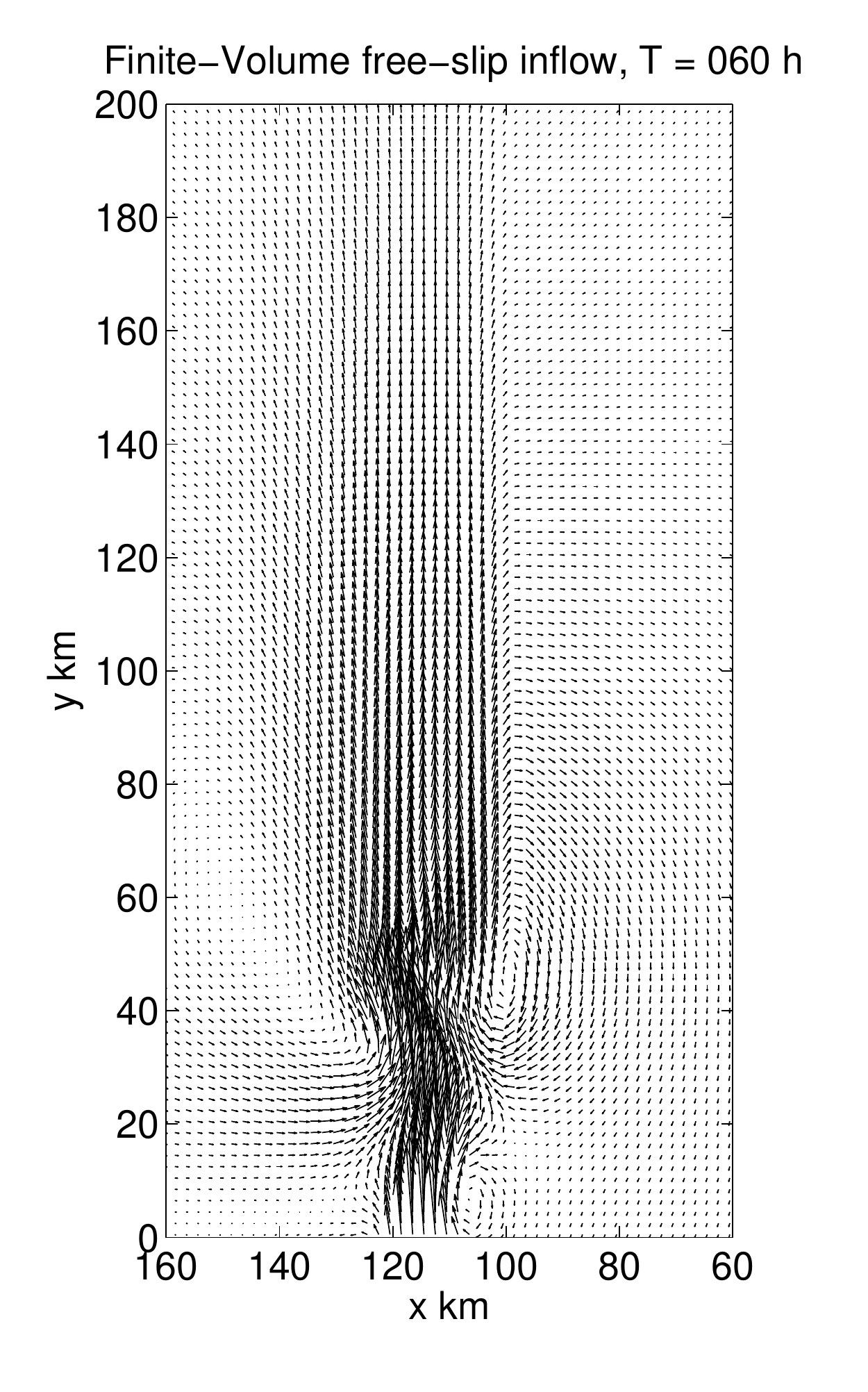}\\
  \vspace{-0.5cm}
  \includegraphics[width=0.32\linewidth]{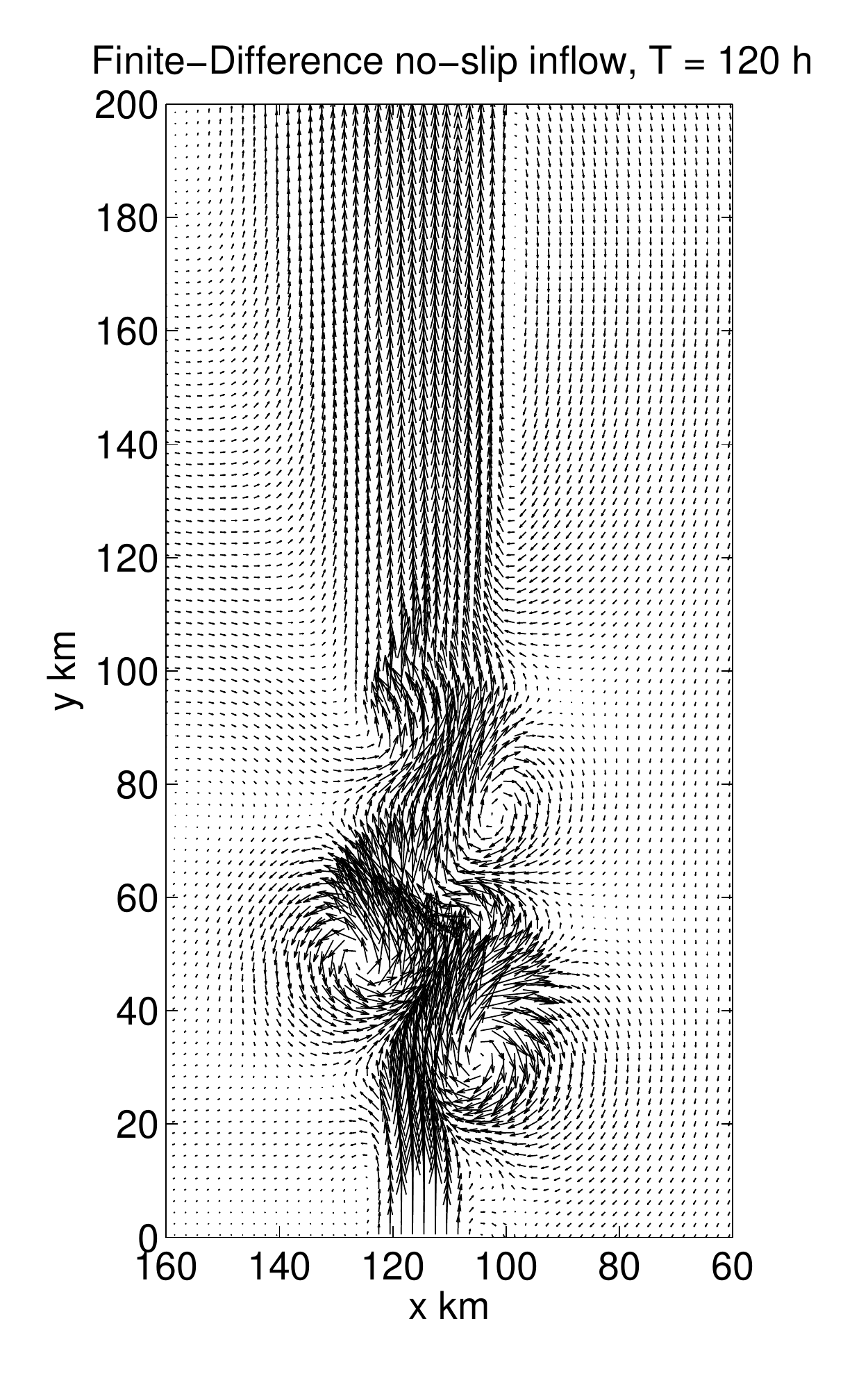}
  \includegraphics[width=0.32\linewidth]{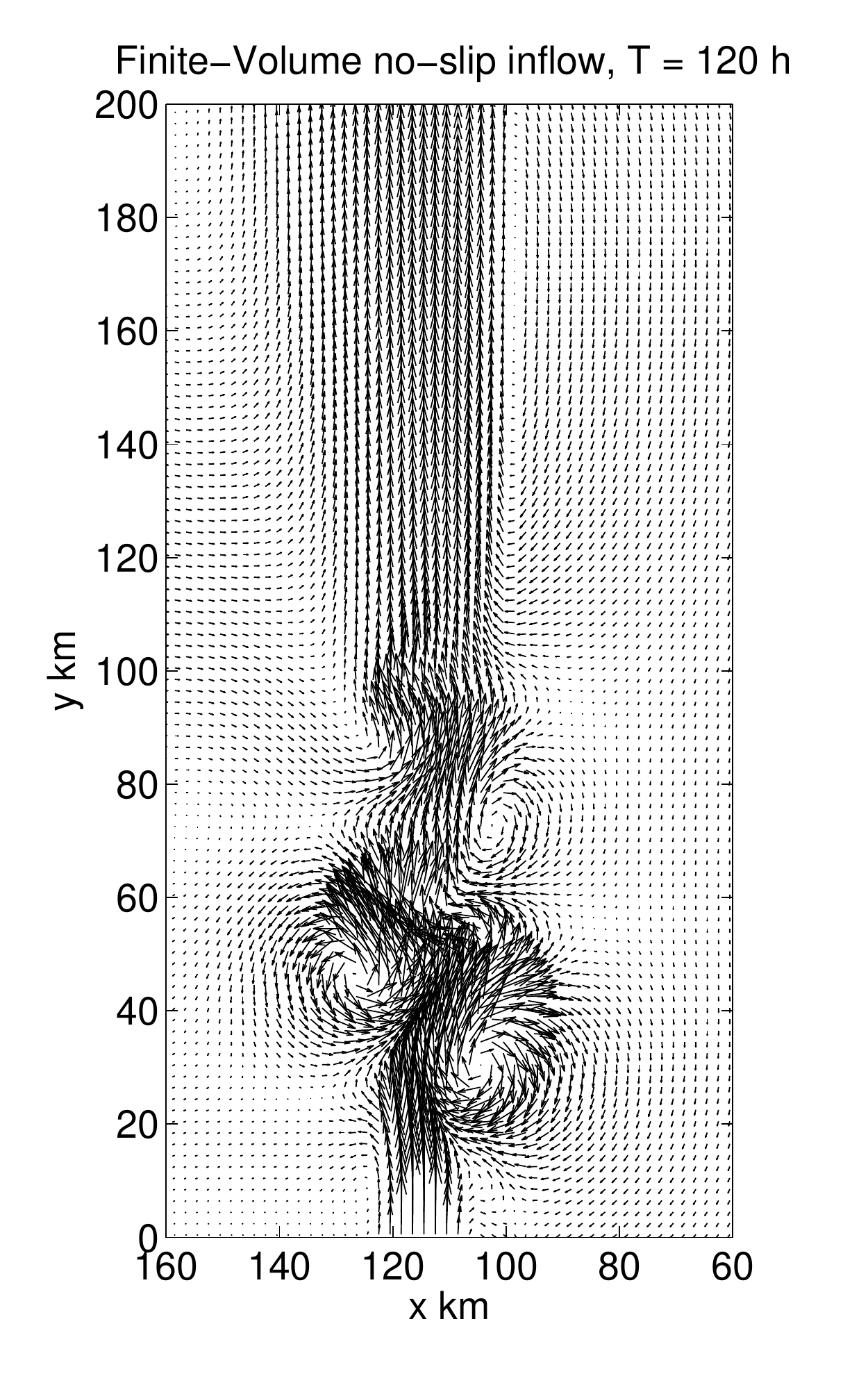}
  \includegraphics[width=0.32\linewidth]{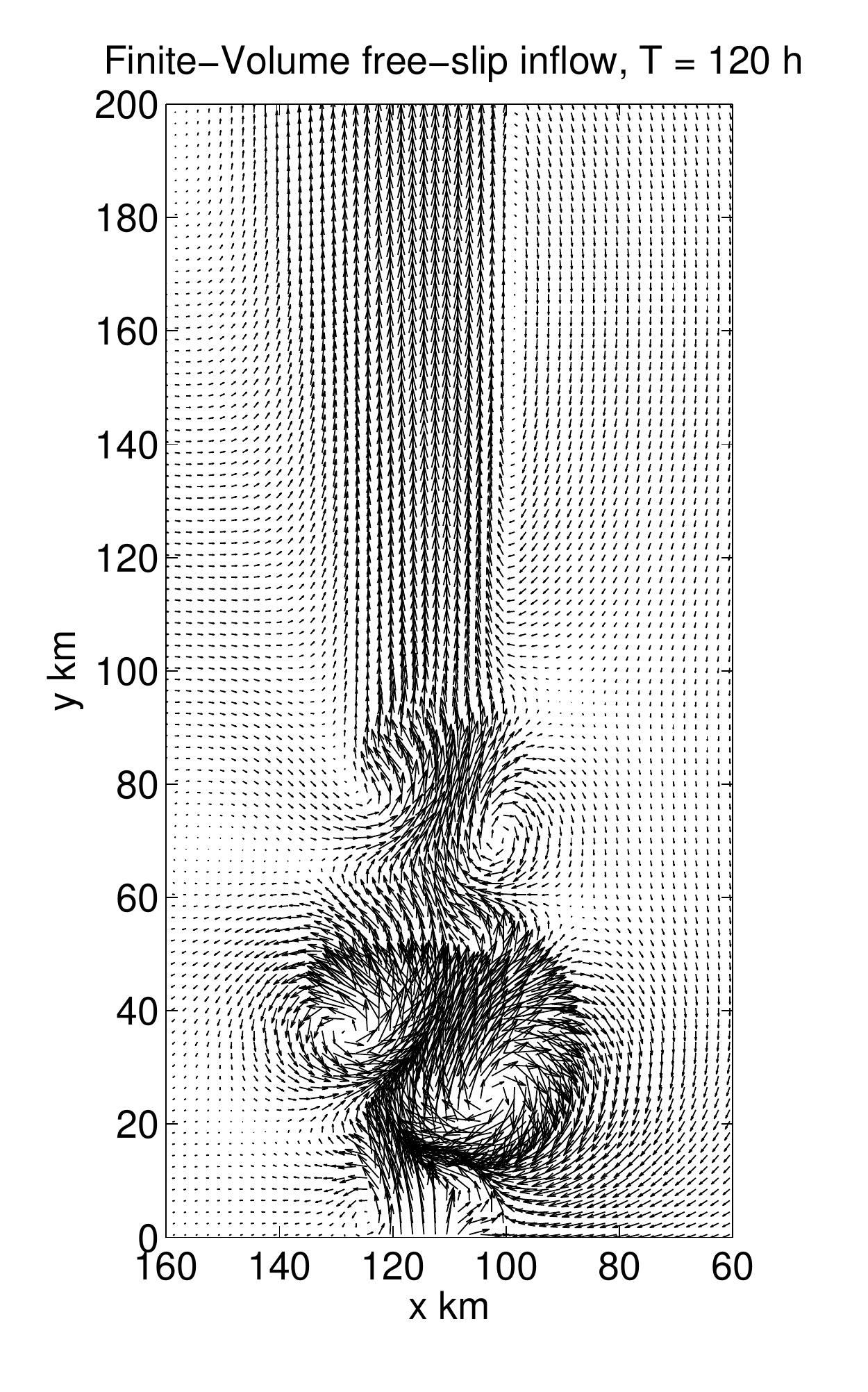}\\
  \vspace{-0.5cm}
  \includegraphics[width=0.32\linewidth]{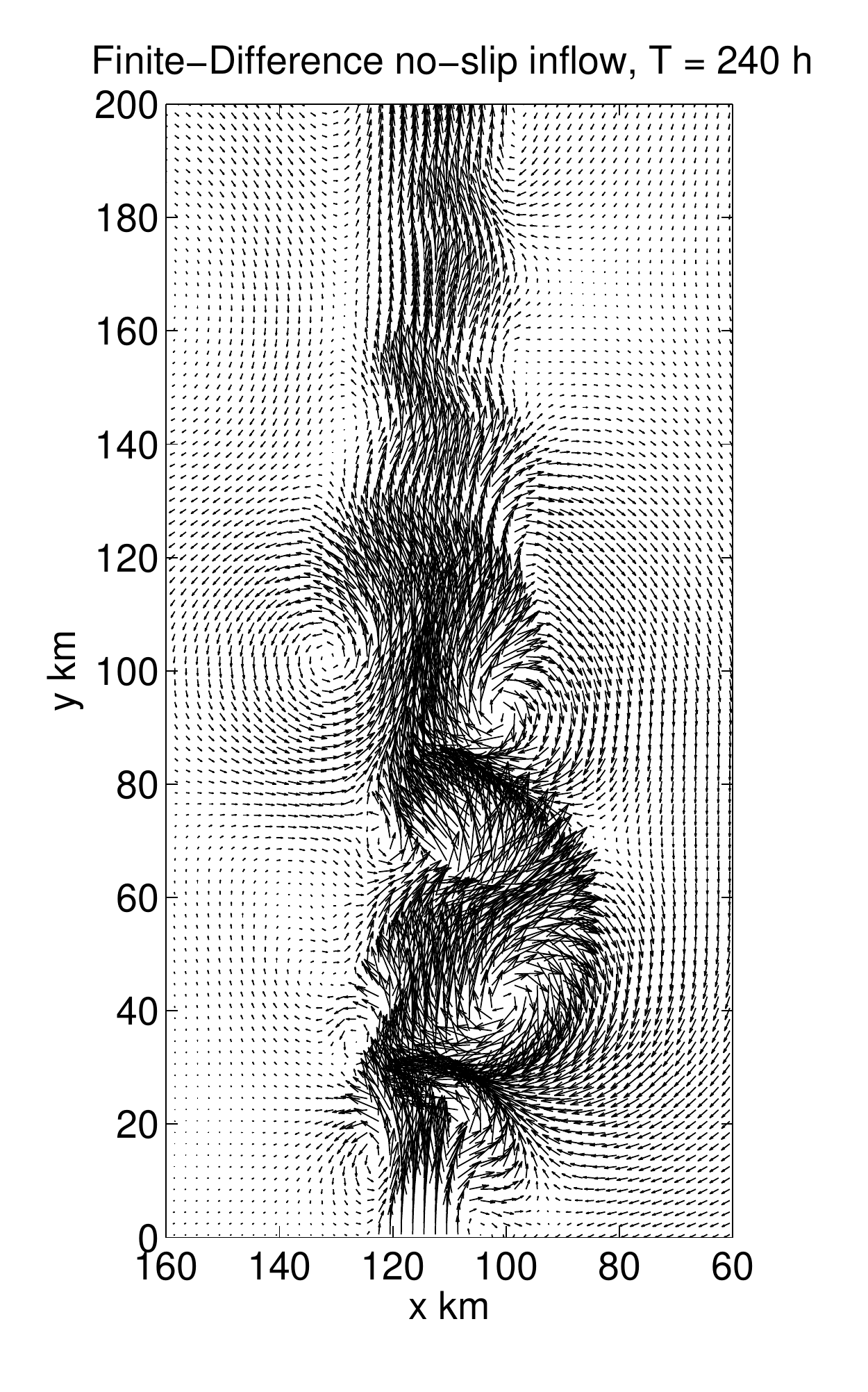}
  \includegraphics[width=0.32\linewidth]{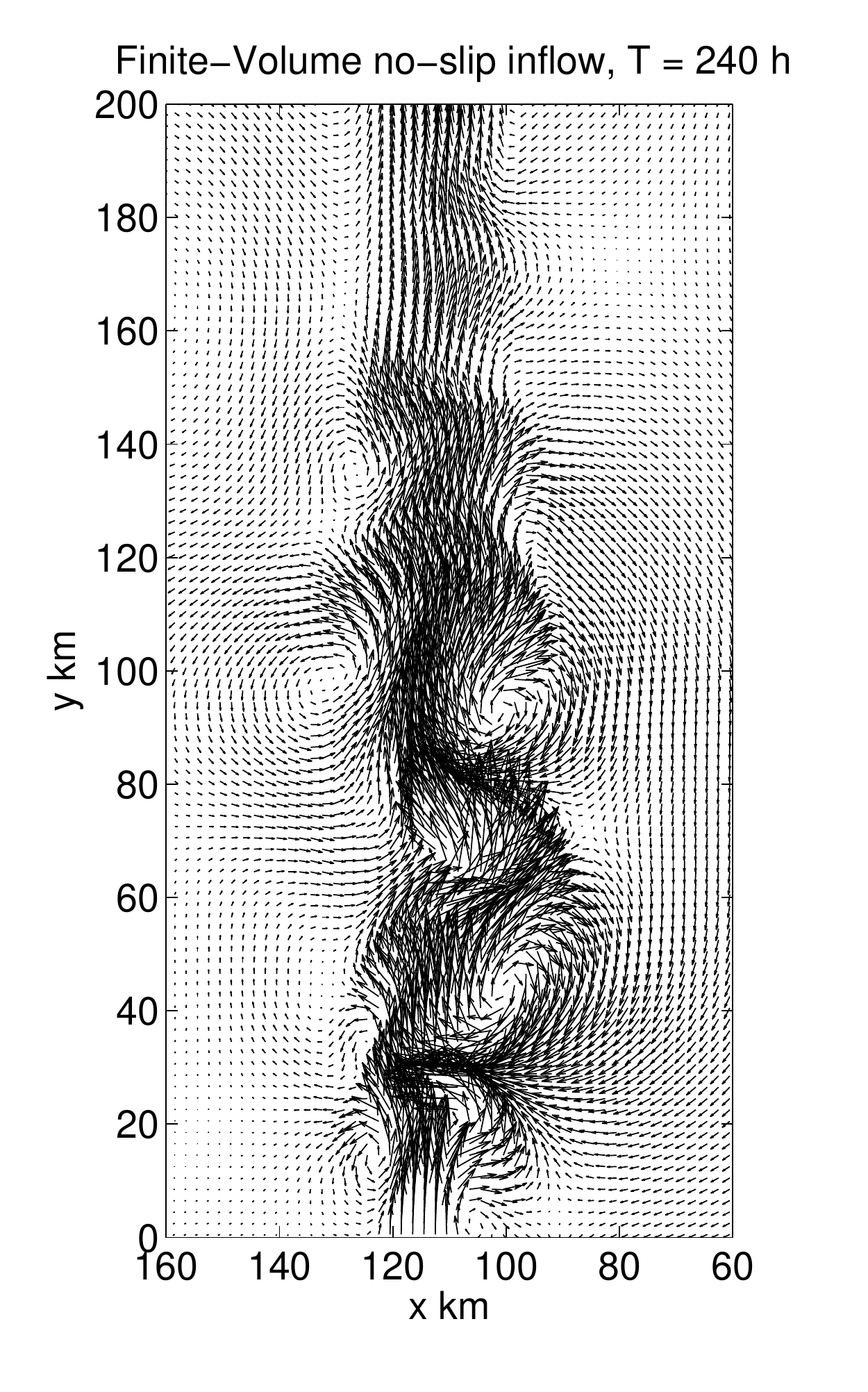}
  \includegraphics[width=0.32\linewidth]{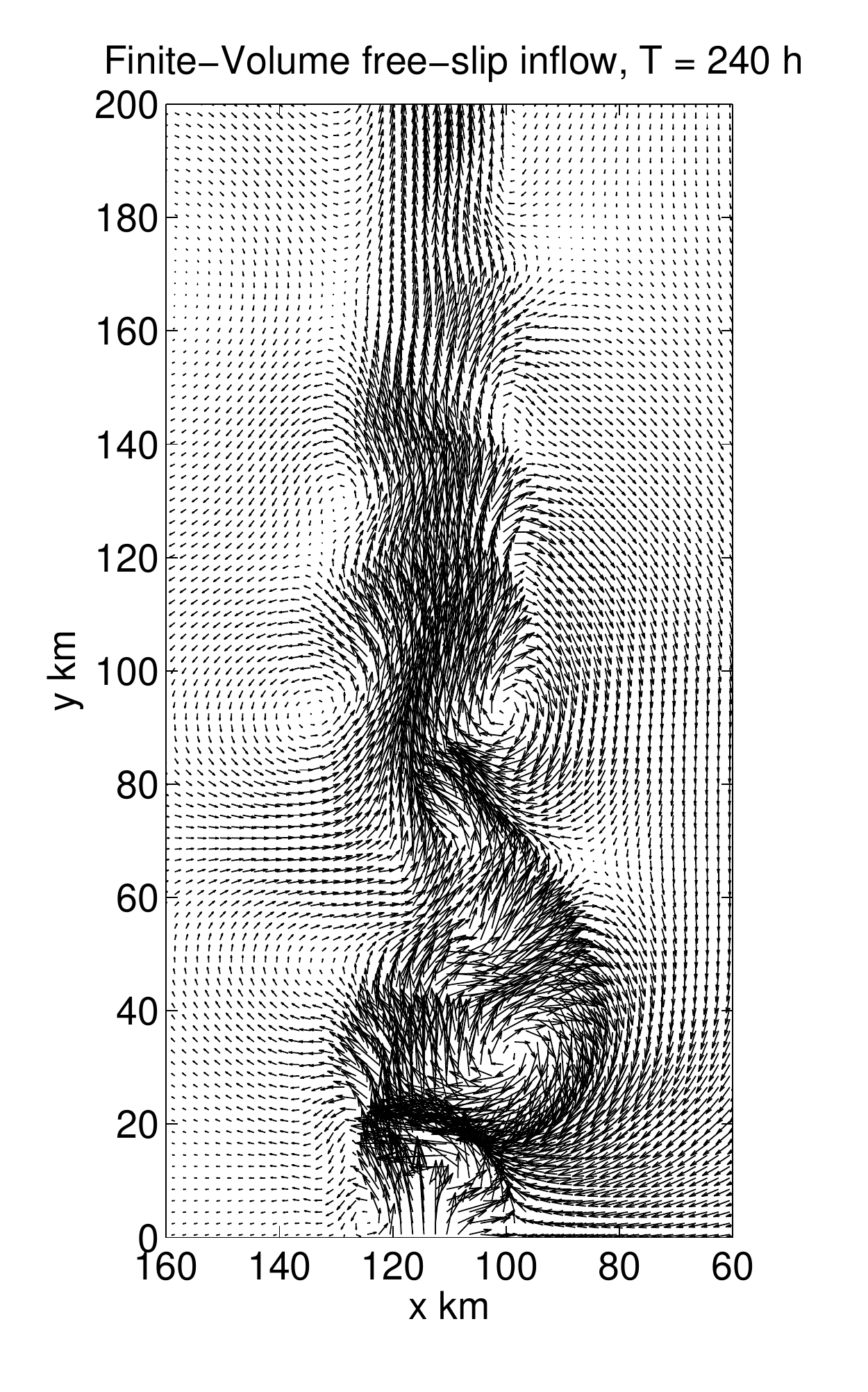}
  \caption{\small\label{fig:norwegian_shelfSurface_E1}Ormen Lange
    Experiment~I. Velocity plots at $60$ (top), $120$, and $240$
    (bottom) hours, computed with the Finite-Difference scheme no-slip
    boundary condition (left), Finite-Volume no-slip boundary
    condition and Finite-Volume free-slip boundary condition (right).}
\end{figure}

\begin{figure}
  \centering
  \includegraphics[width=0.32\linewidth]{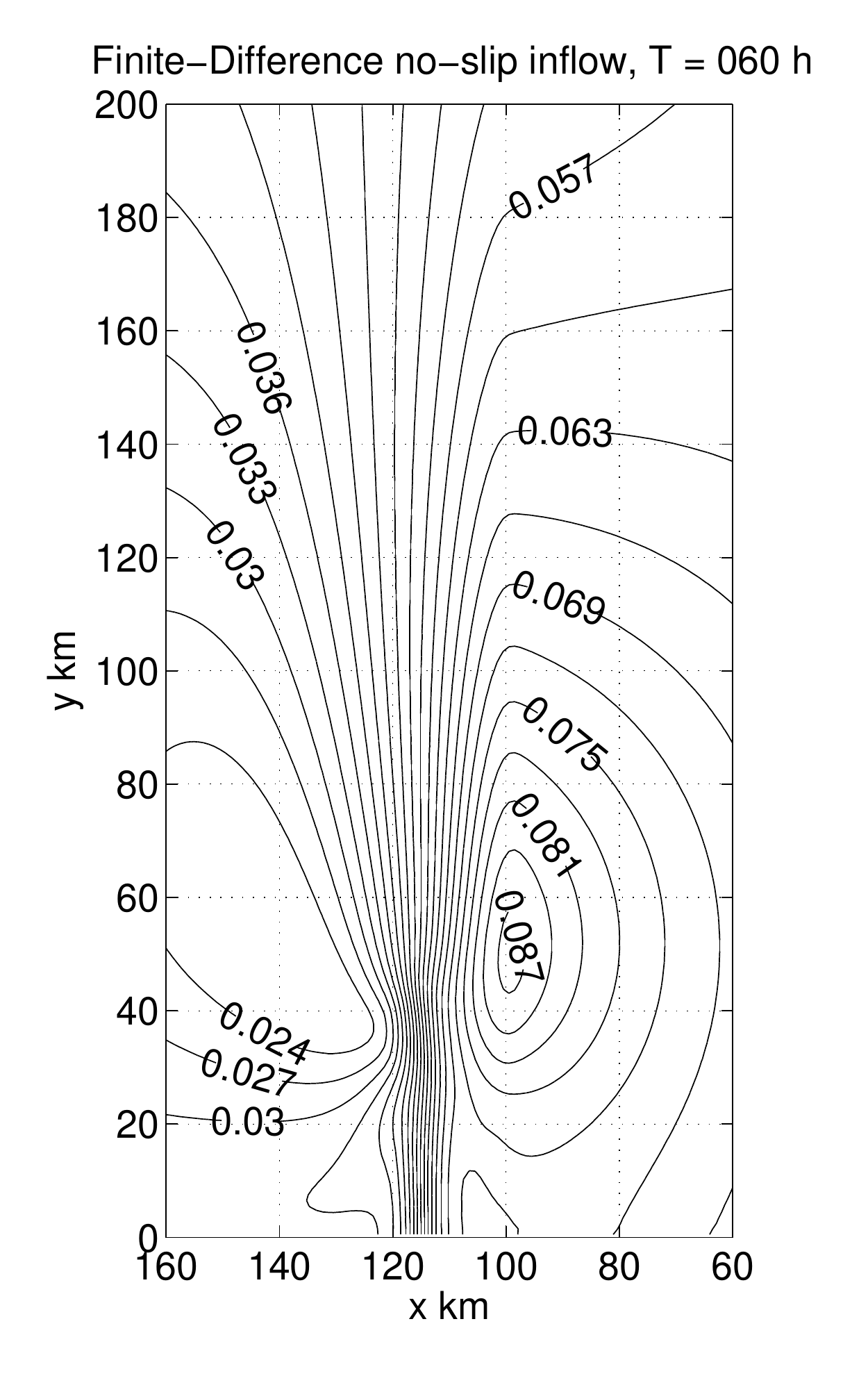}
  \includegraphics[width=0.32\linewidth]{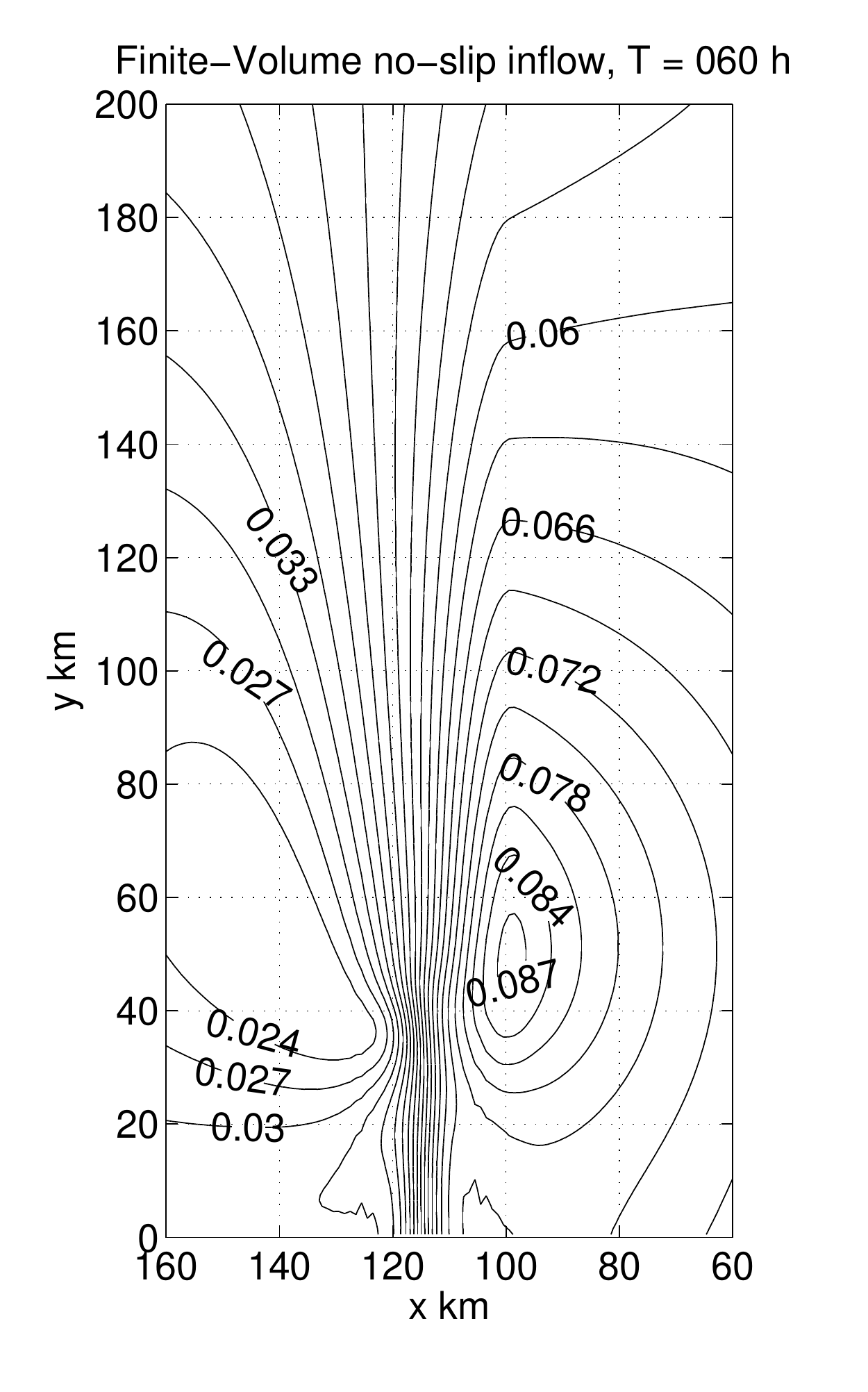}
  \includegraphics[width=0.32\linewidth]{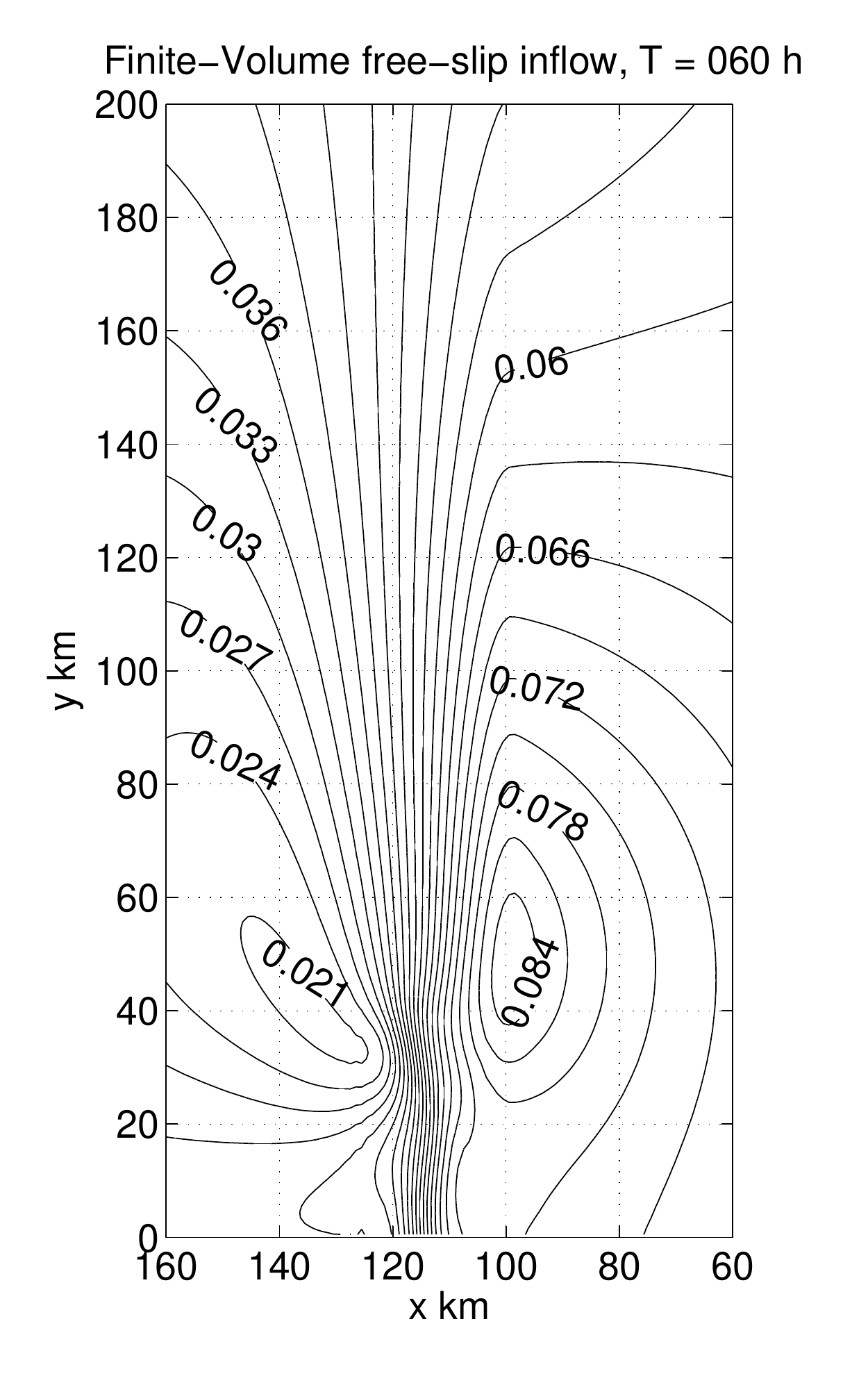}\\
  \vspace{-0.5cm}  
  \includegraphics[width=0.32\linewidth]{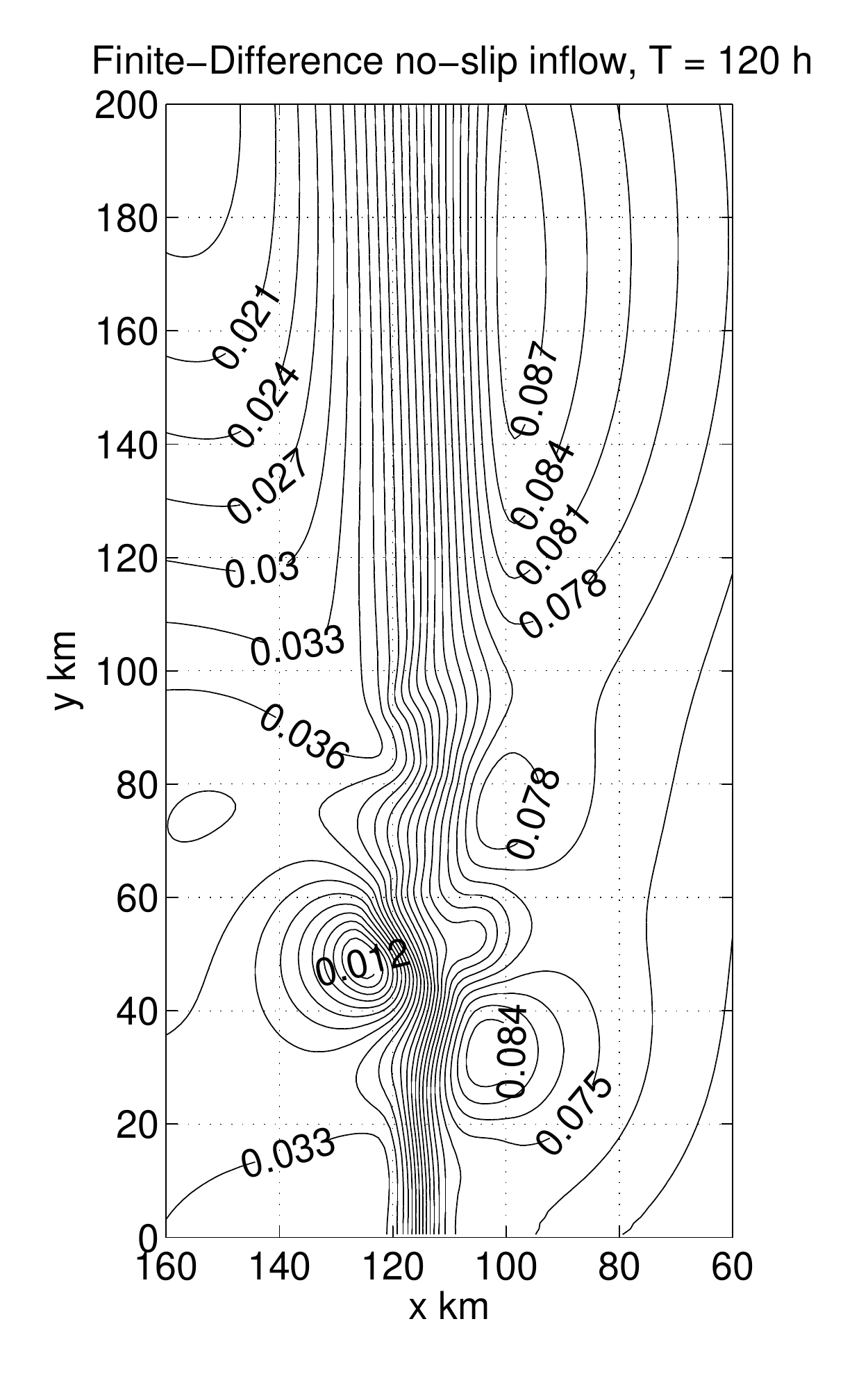}
  \includegraphics[width=0.32\linewidth]{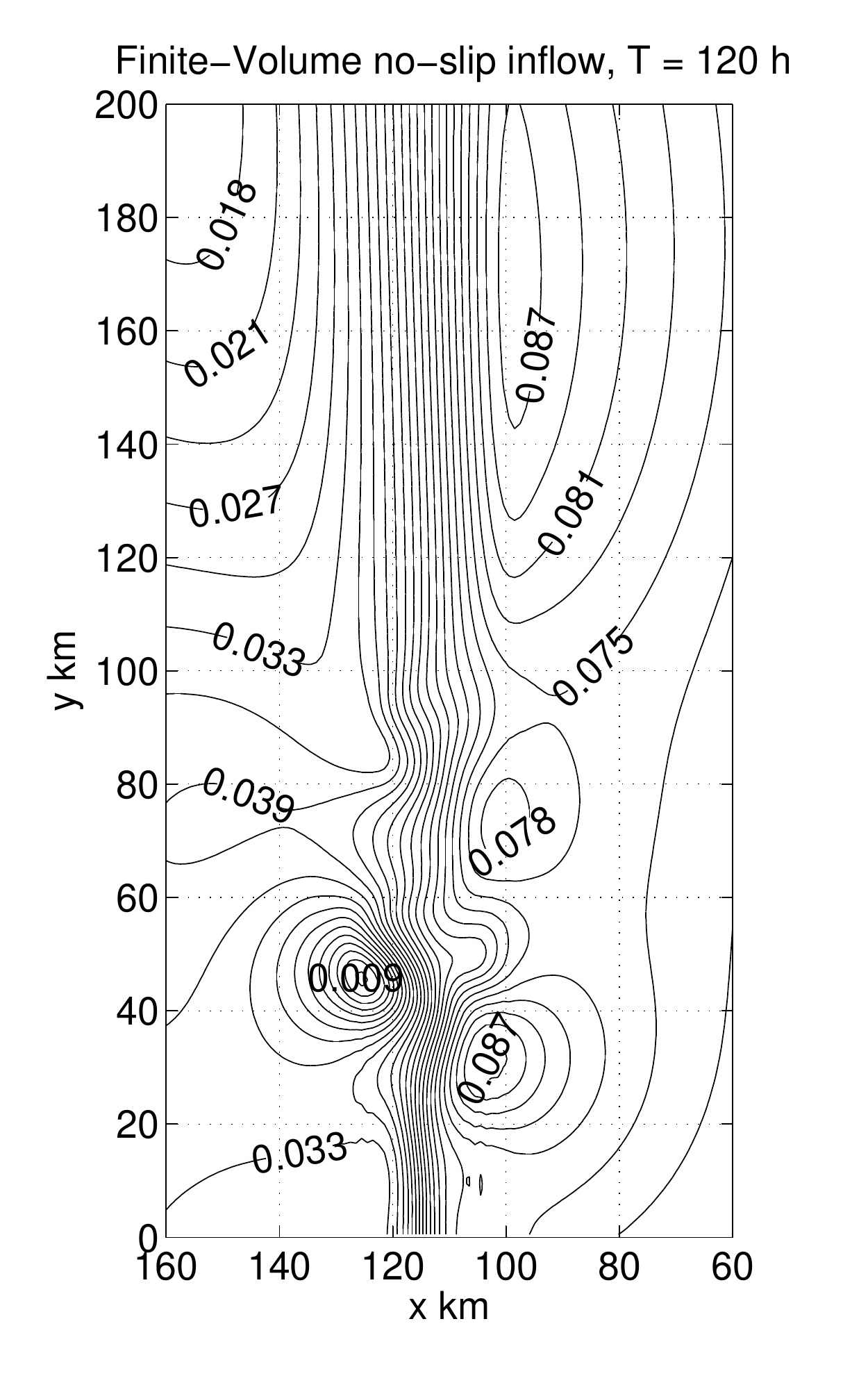}
  \includegraphics[width=0.32\linewidth]{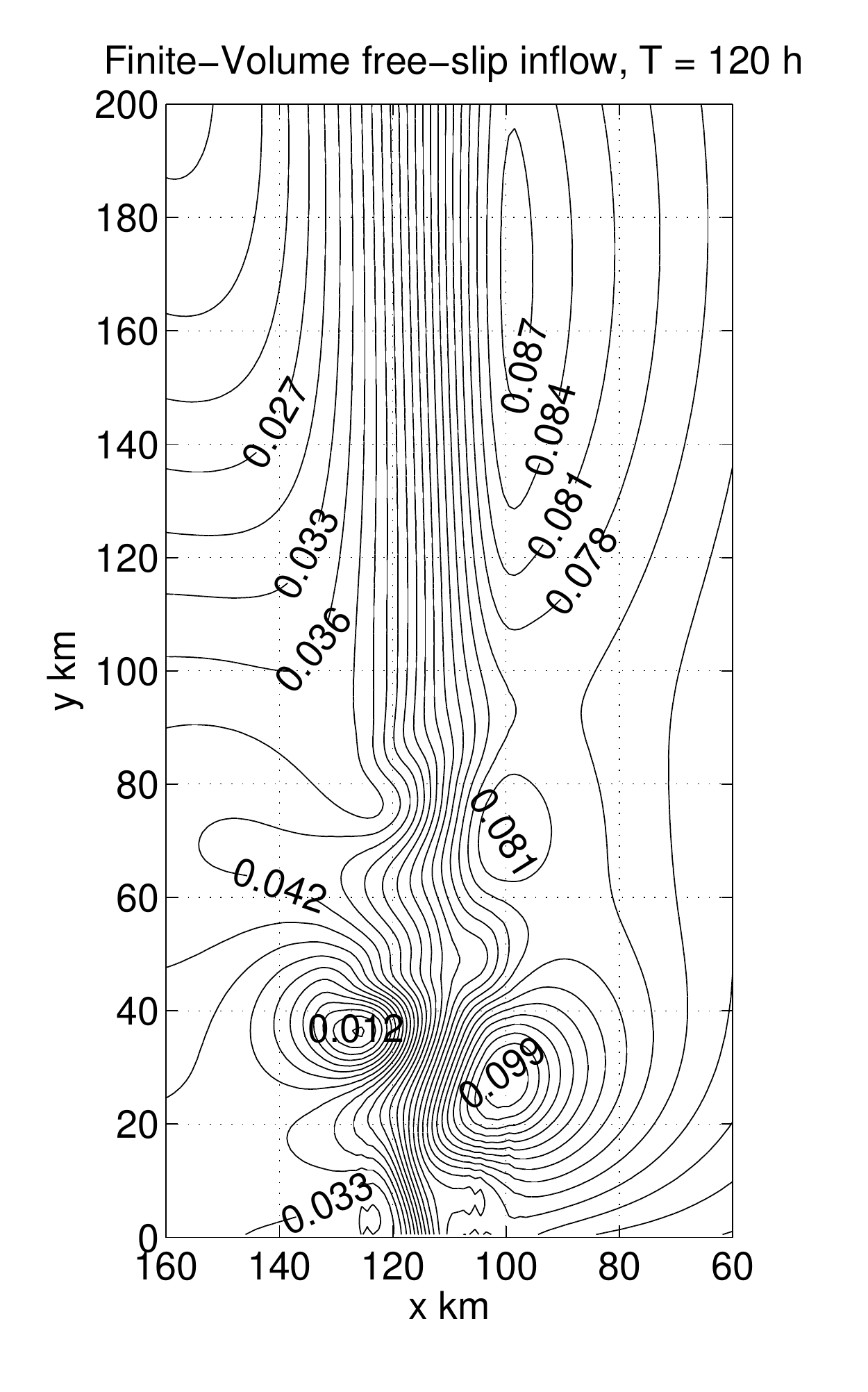}\\
  \vspace{-0.5cm}  
  \includegraphics[width=0.32\linewidth]{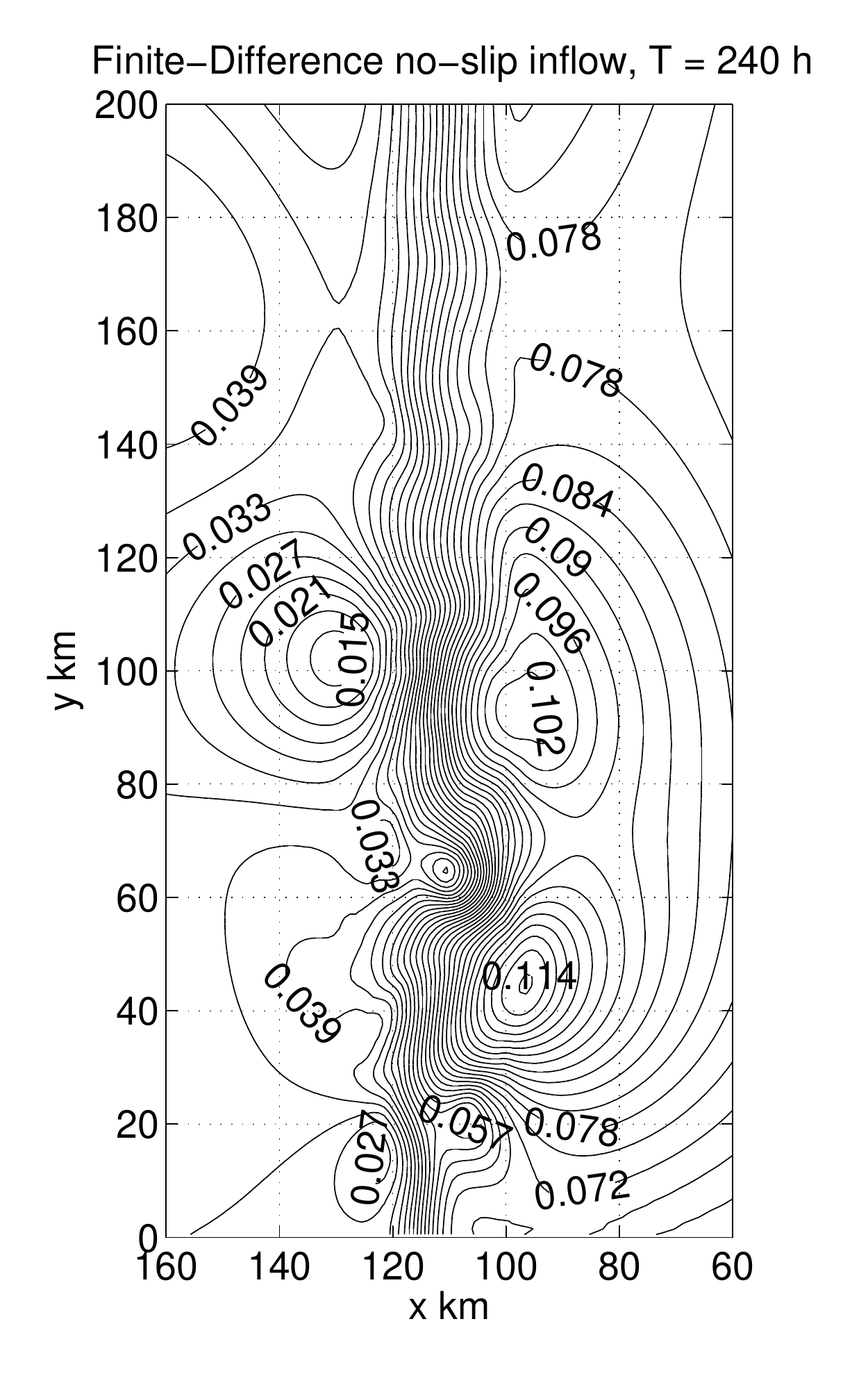}
  \includegraphics[width=0.32\linewidth]{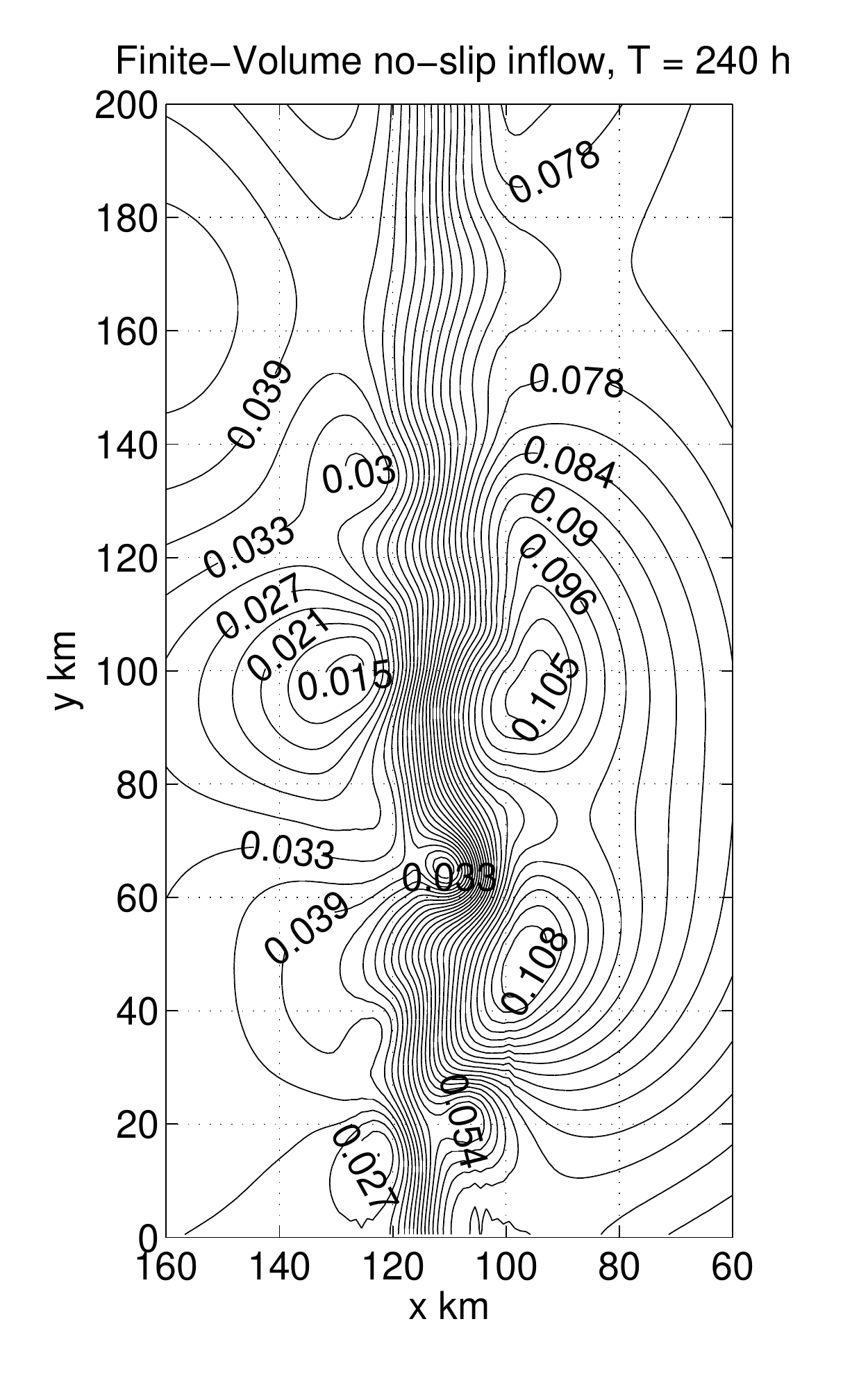}
  \includegraphics[width=0.32\linewidth]{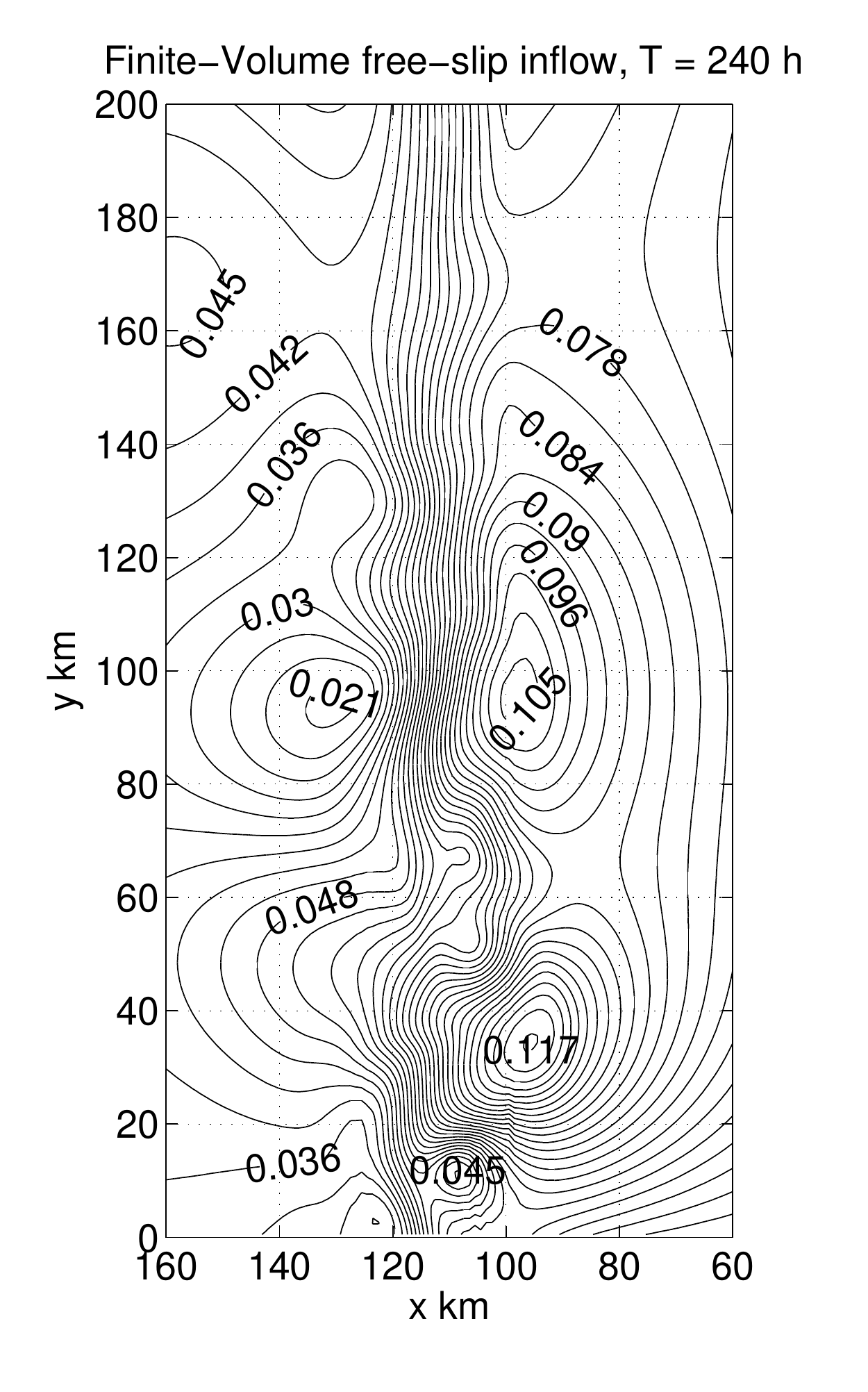}
  \caption{\small\label{fig:norwegian_shelfVelocity_E1}Ormen Lange
  Experiment~I. Contour plots of surface elevation at $60$ (top),
  $120$, and $240$ (bottom) hours, computed with the Finite-Difference
  scheme no-slip boundary condition (left), Finite-Volume no-slip
  boundary condition and Finite-Volume free-slip boundary condition
  (right).}
\end{figure}

The results of the three computations after 60, 120, and 240 hours are
shown in Figures~\ref{fig:norwegian_shelfSurface_E1} and
\ref{fig:norwegian_shelfVelocity_E1}.  The plots of the
finite-difference and finite-volume solutions with no-slip inflow
boundary condition look quite similar. After a short time, the narrow
current starts to oscillate and large eddies are generated. However,
we would like to point out that in addition to these physical
oscillations the finite-difference develops large numerical
oscillations, which we damp by adding artificial diffusion as in
\cite{Gjevik_etal_2002} Equation~(4) by adding eddy viscosity $\nu$,
given by
\begin{align}
\nu = q l^2\left[\left(\frac{\partial \bar{u}}{\partial x} \right)^2 +   
     \left(\frac{\partial \bar{u}}{\partial x} + \frac{\partial \bar{v}}{\partial y} \right)^2  +
     \left(\frac{\partial \bar{v}}{\partial y} \right)^2 \right]^{\oh},
\end{align}
according to Smagorinsky~\cite{Smagorinsky_1963}. Where $l$ denotes the grid size and the depth mean
current velocity defined to first order by,
\begin{align}
\bar{u}=\frac{U}{H}, \quad \bar{v}=\frac{V}{H}.
\end{align} 
The diffusion parameter $q$ is set to $q=0.1$ in all finite-difference
computation.  The finite-volume solution with the free-slip inflow
boundary condition looks different, eddies are close to the inflow.  
%\marginlabel{\color{red} Bj\o rn, can you possibly tell which solution is more physical?}

%
%
%===========================================================================================
\subsubsection{Setup for Ormen Lange Shelf Experiment~II} \label{sec:BarostrophicJetSetupE2}
%===========================================================================================
%
In the setup of Section~\ref{sec:BarostrophicJetSetupE1}, the
in-flowing jet was cut off at $x=L_B\pm B/2$. After some time, these
points become transition points with a noticeable discontinuous shear
layer. For the next experiment we avoid such a discontinuous shear
layer and change the boundary condition by assigning the in-flowing
jet profile $v_{jet}$ defined in equation \eqref{eq:vjet2} on the
whole southern boundary.

\begin{figure}
  \centering
  \includegraphics[width=0.32\linewidth]{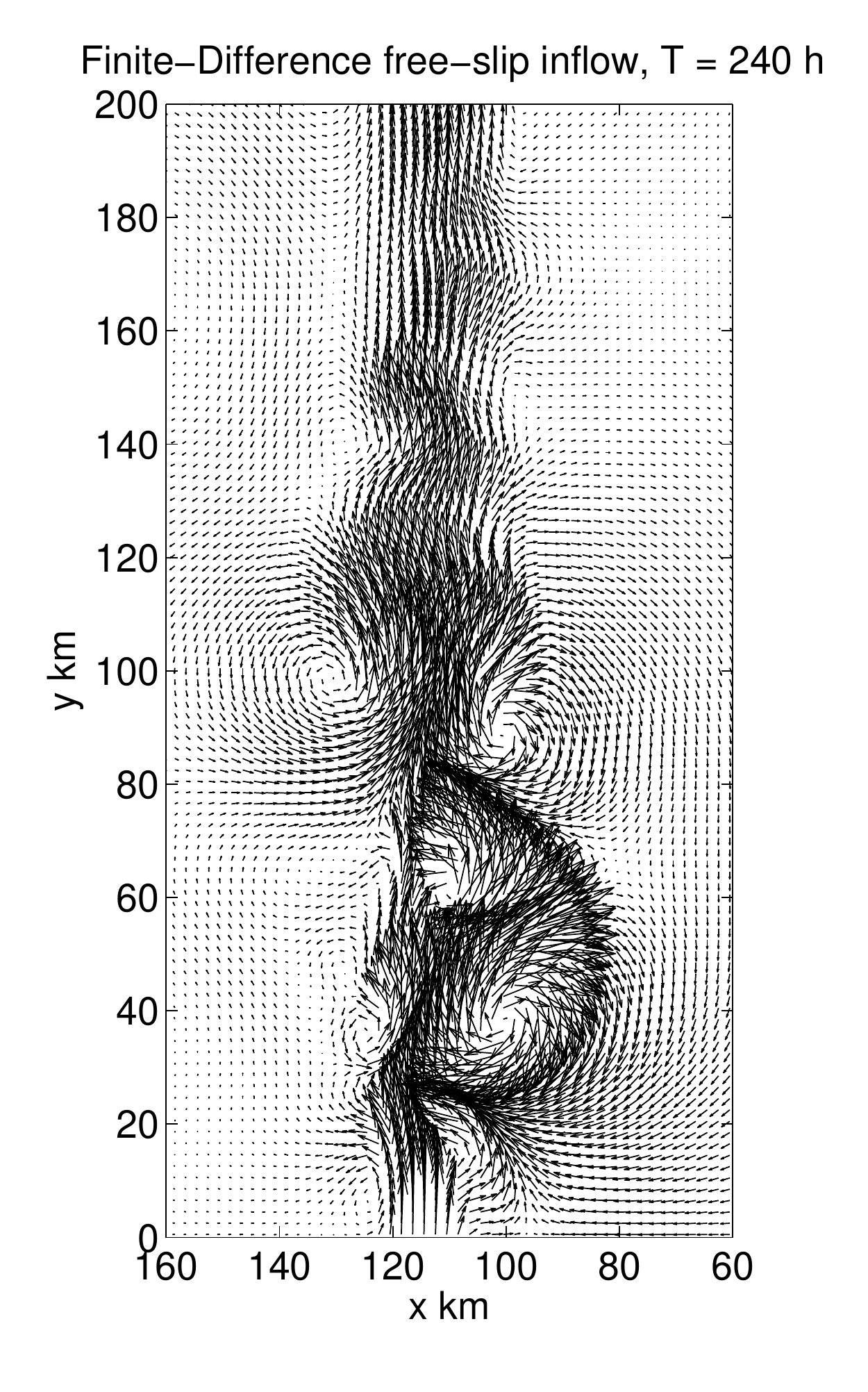}
  \includegraphics[width=0.32\linewidth]{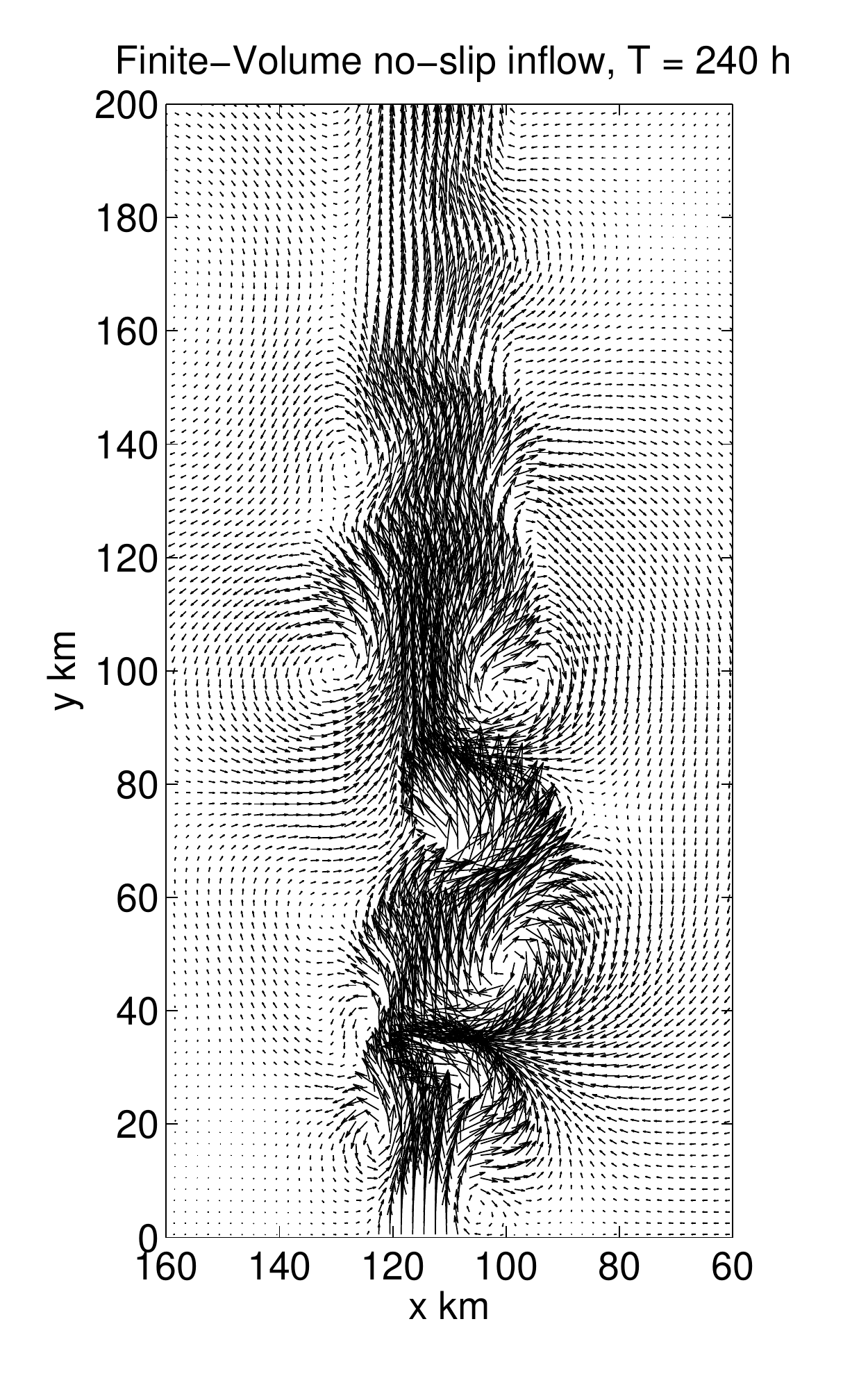}
  \includegraphics[width=0.32\linewidth]{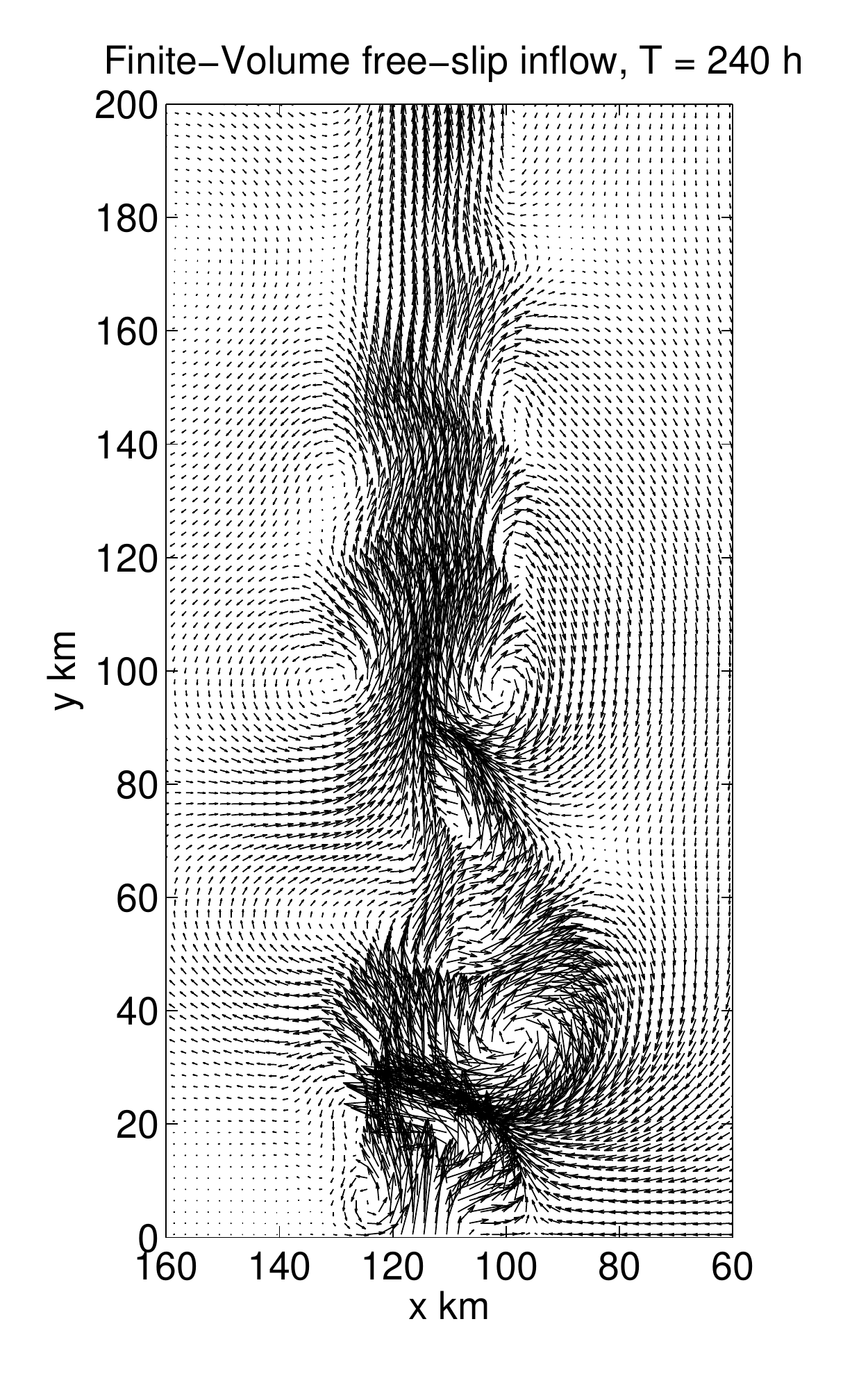}\\
  \vspace{-0.5cm}
  \includegraphics[width=0.32\linewidth]{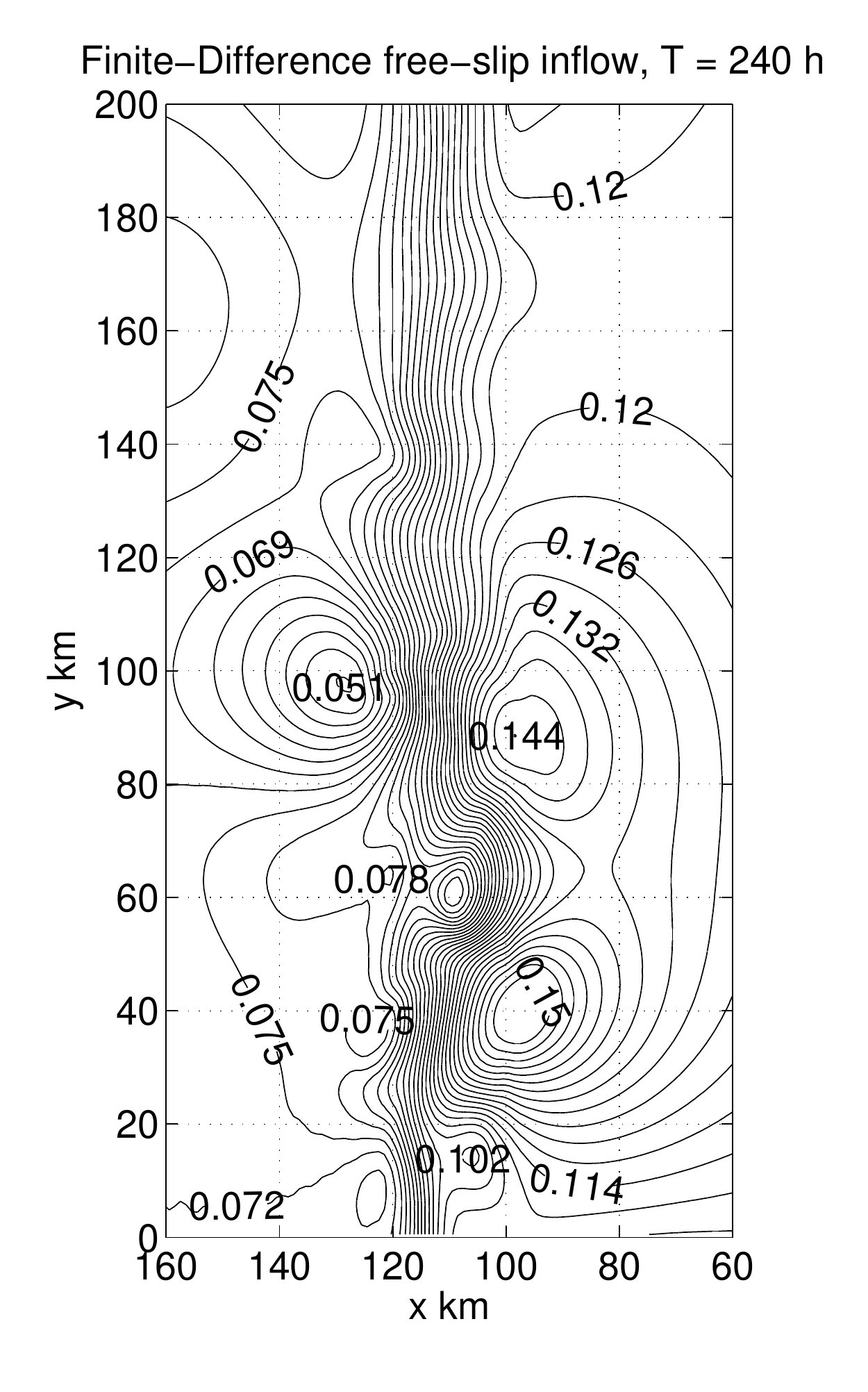}  
  \includegraphics[width=0.32\linewidth]{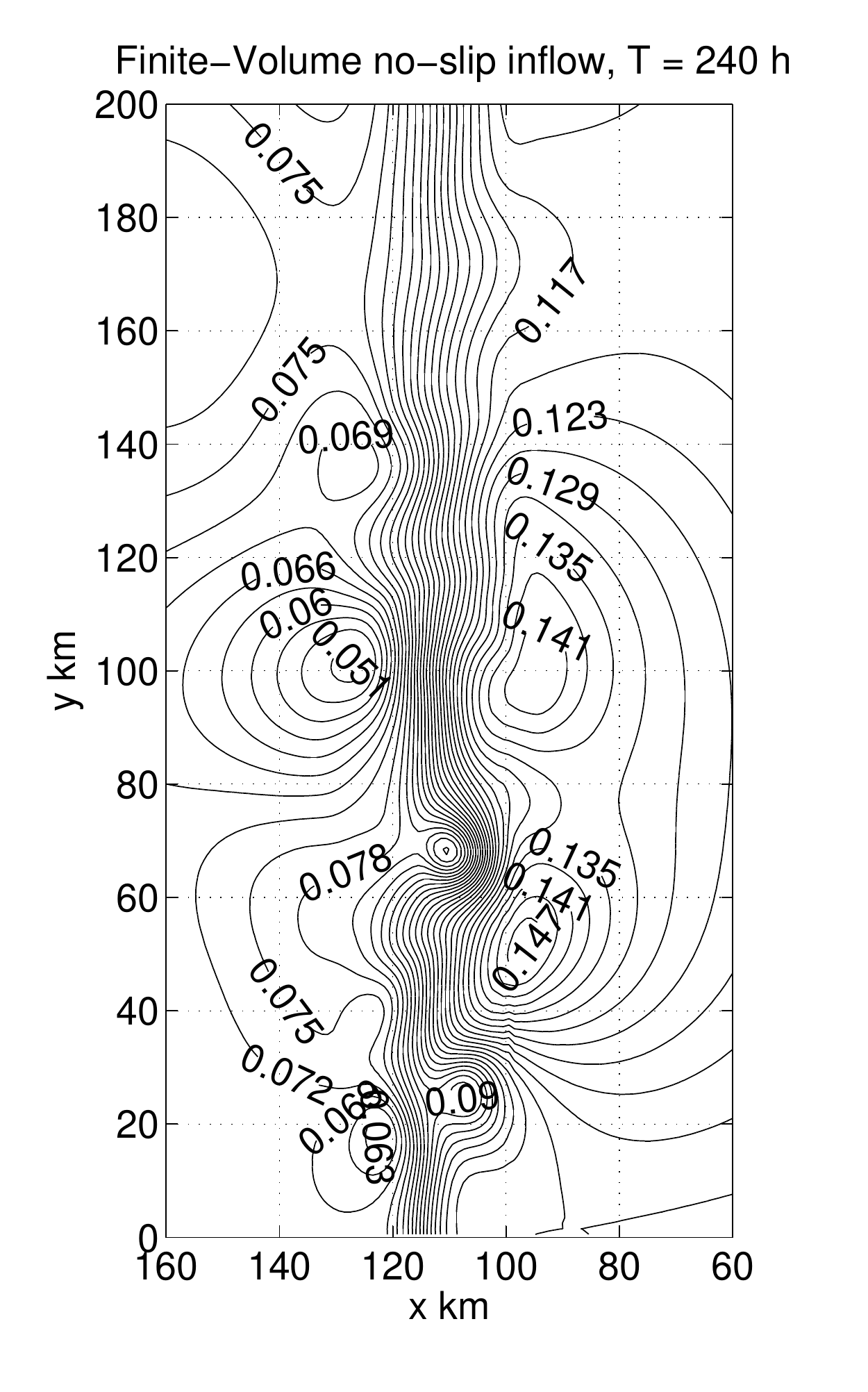}
  \includegraphics[width=0.32\linewidth]{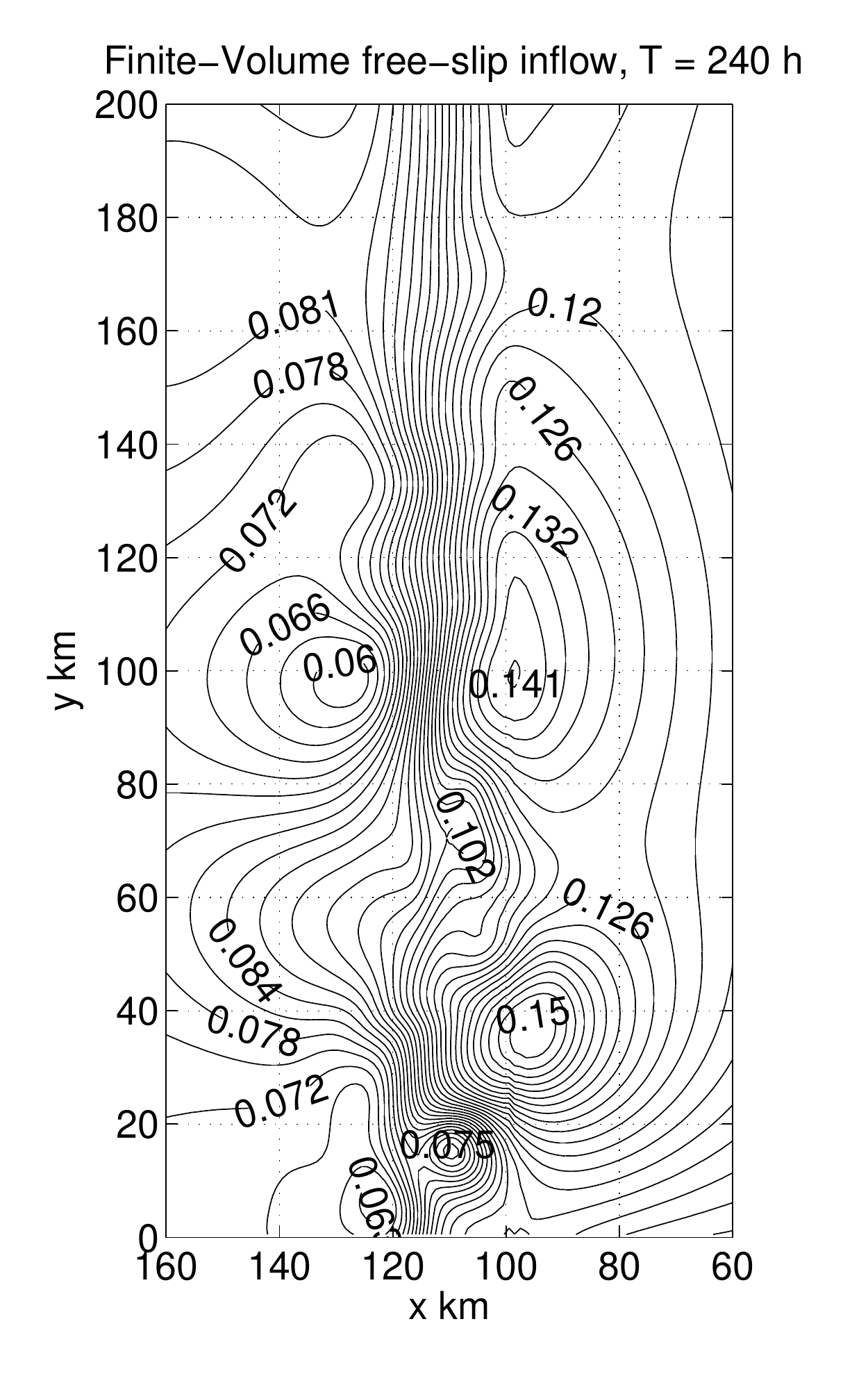}  
  \caption{\small\label{fig:norwegian_shelfVelocity_E2}Ormen Lange
  Experiment~II.. Contour plots of surface elevation (top) and
  velocity plots (bottom) at $240$ hours, computed with the Finite-Difference
  scheme free-slip boundary condition (left), Finite-Volume no-slip
  boundary condition and Finite-Volume free-slip boundary condition
  (right).}
\end{figure}

The results are displayed in
Figures~\ref{fig:norwegian_shelfVelocity_E2}. 
Eddies are still created and are of similar strength as in the previous section
Figures~\ref{fig:norwegian_shelfSurface_E1}. Note, however, that the
maximal water level is now about 15cm, which is 3.3cm higher than
before (11.7cm). This was to be expected because our new southern
boundary condition does not allow any outflow.

From this experiment we can conclude that the non-smooth patching of
the boundary condition at the southern boundary is not the mechanism
which creates the instability. In the next experiment we will
investigate if the instability is effected by the start-up procedure.

%
%===========================================================================================
\subsubsection{Setup for Ormen Lange Shelf Experiment~III} \label{sec:BarostrophicJetSetupE3}
\begin{figure}
  \centering
  \includegraphics[width=0.6\linewidth]{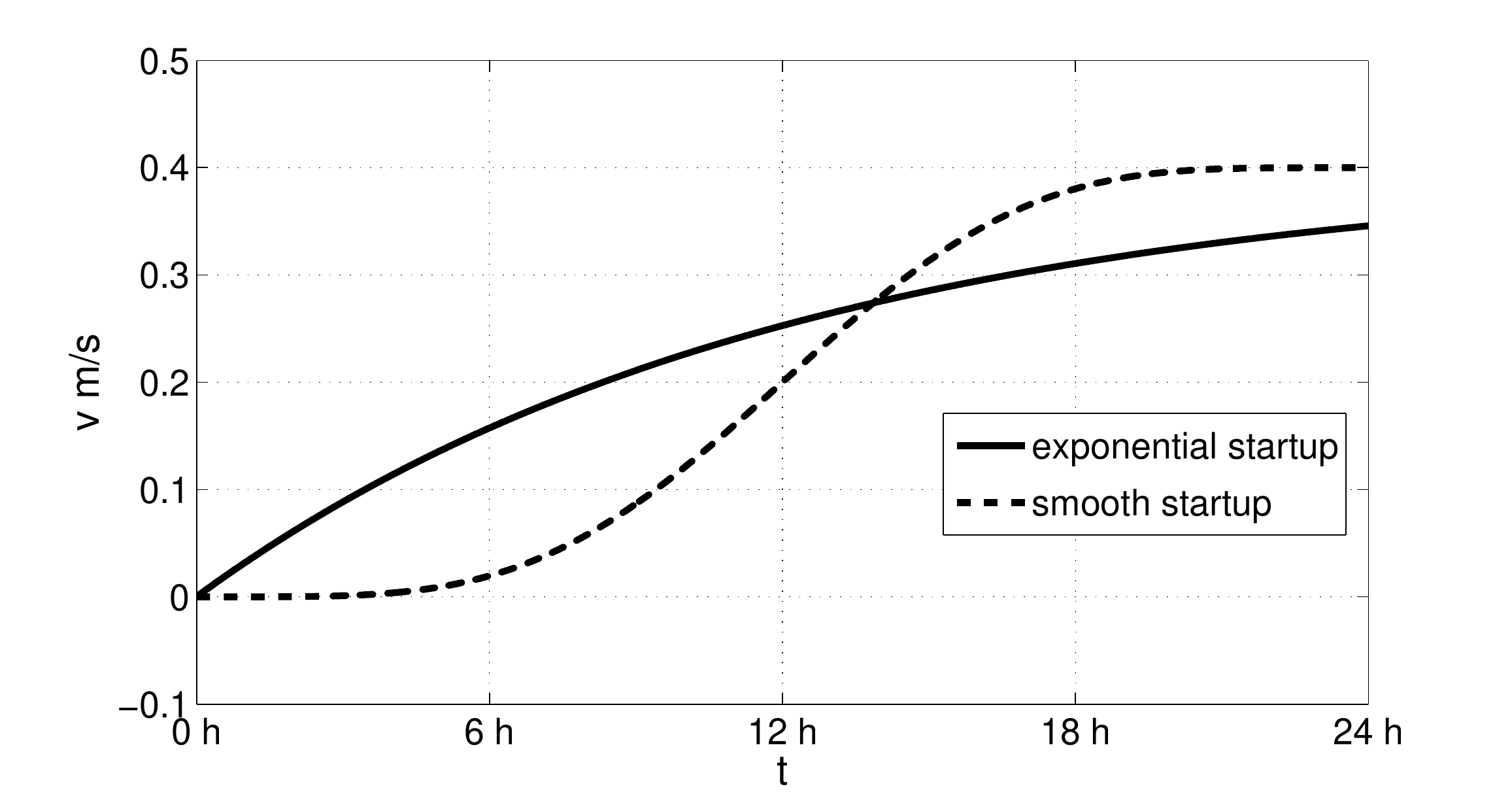}\\
  \caption{\small\label{fig:startup}Smooth startup function $\gamma(\tau)$
    compared with exponential startup.}
\end{figure}

%===========================================================================================
We use the same setup as in
Subsection~\ref{sec:BarostrophicJetSetupE2} but here we use the smooth
(four times continuously differentiable) growth function $\gamma(\tau)$ of \eqref{eq:growth_function}
with $\tau=t/24h$, see Figure~\ref{fig:startup}. 
Eddies are still being created and are qualitatively about
the same as before. This is also suggested
by the linear stability analysis in \cite{Gjevik_etal_2002}.
\begin{figure}
  \centering
  \includegraphics[width=0.32\linewidth]{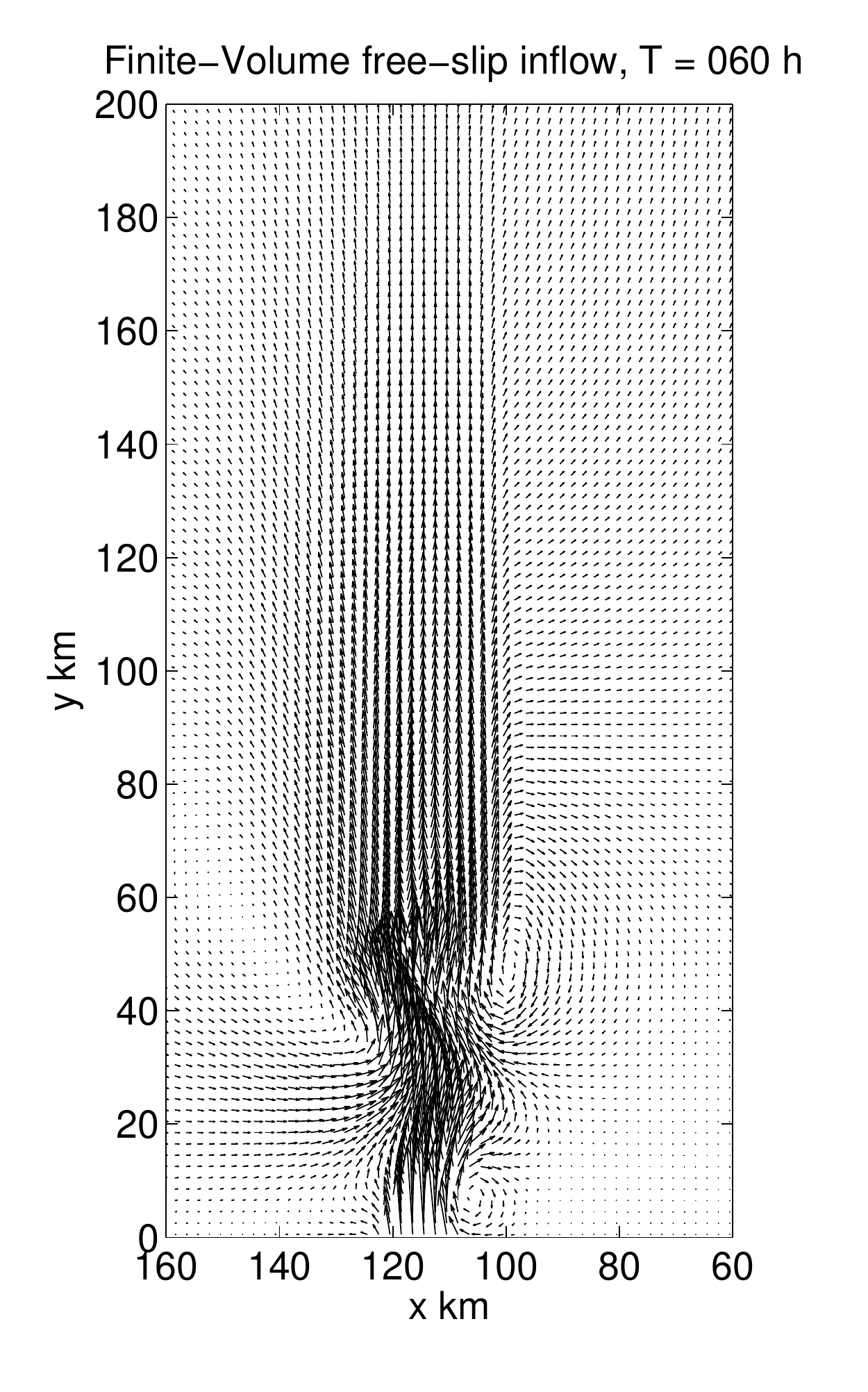}
  \includegraphics[width=0.32\linewidth]{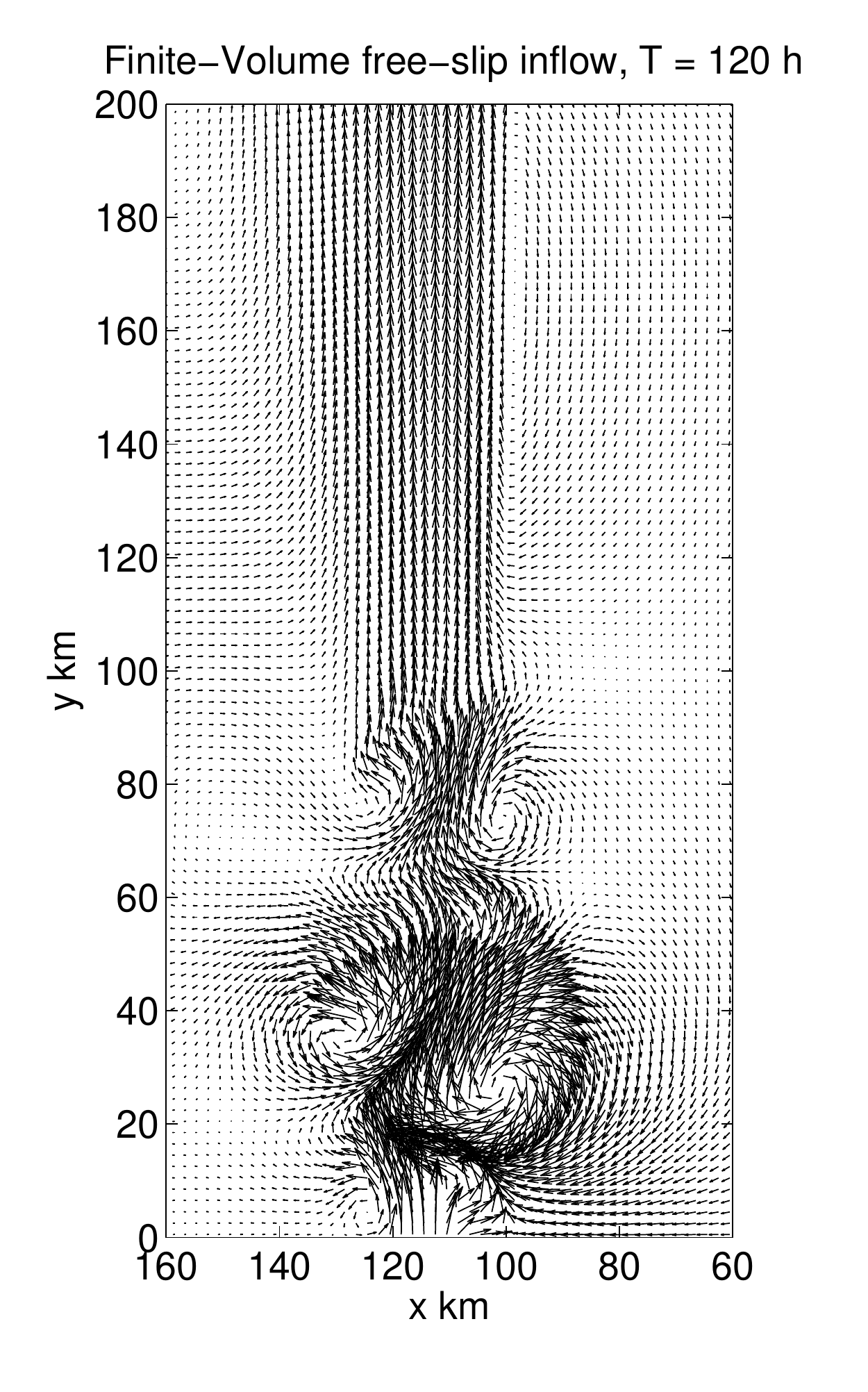}
  \includegraphics[width=0.32\linewidth]{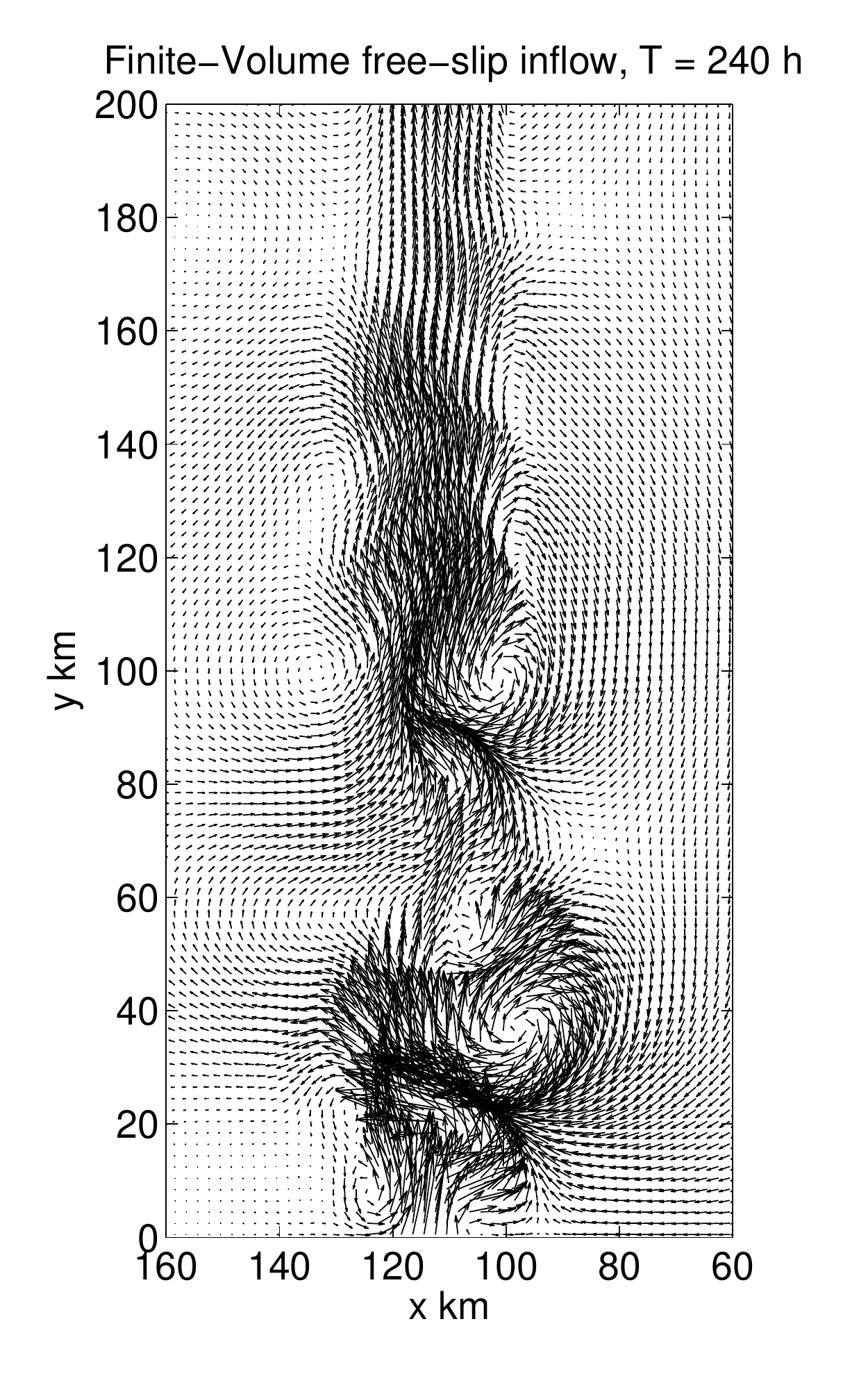}\\
  \vspace{-0.5cm}  
  \includegraphics[width=0.32\linewidth]{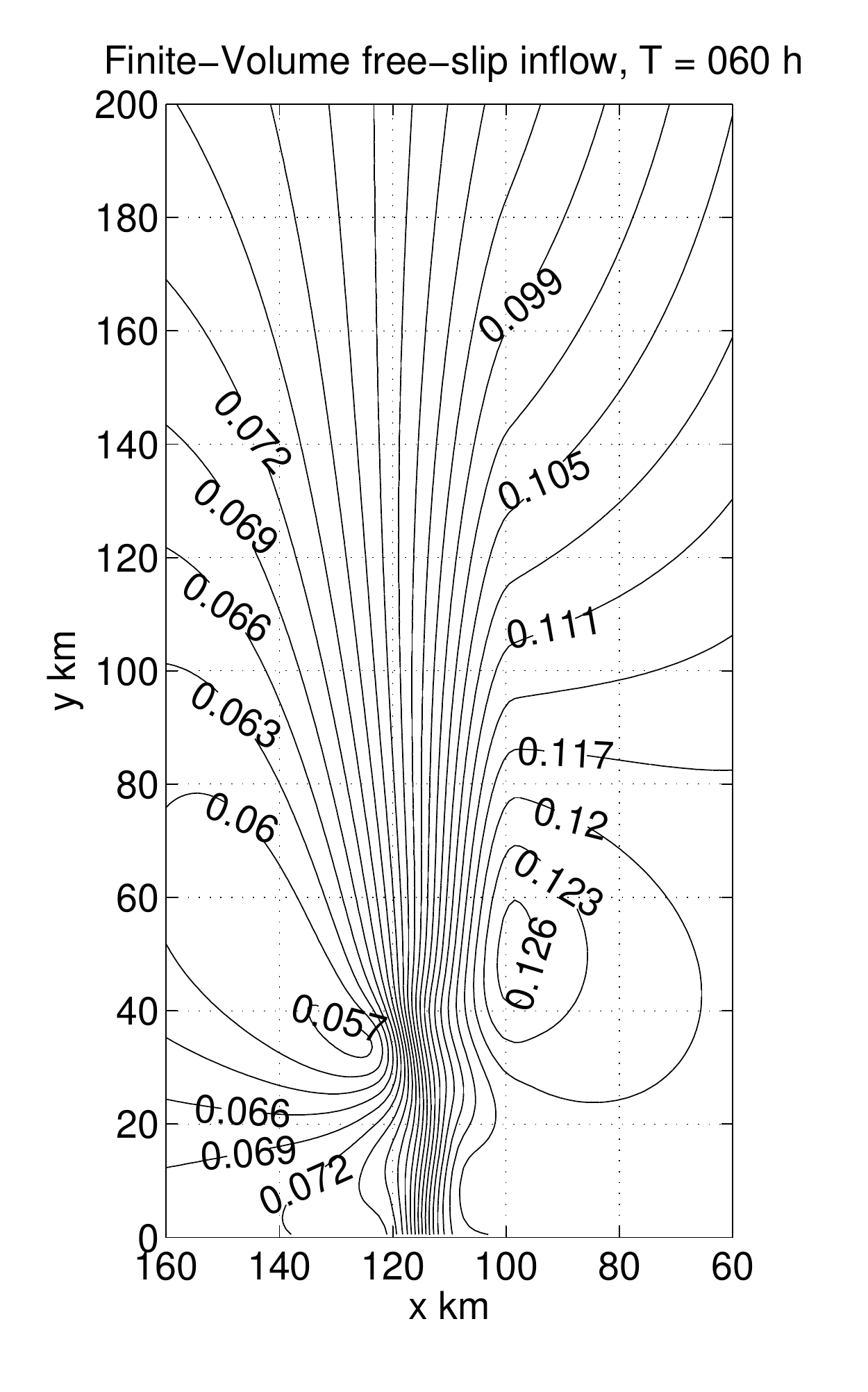}
  \includegraphics[width=0.32\linewidth]{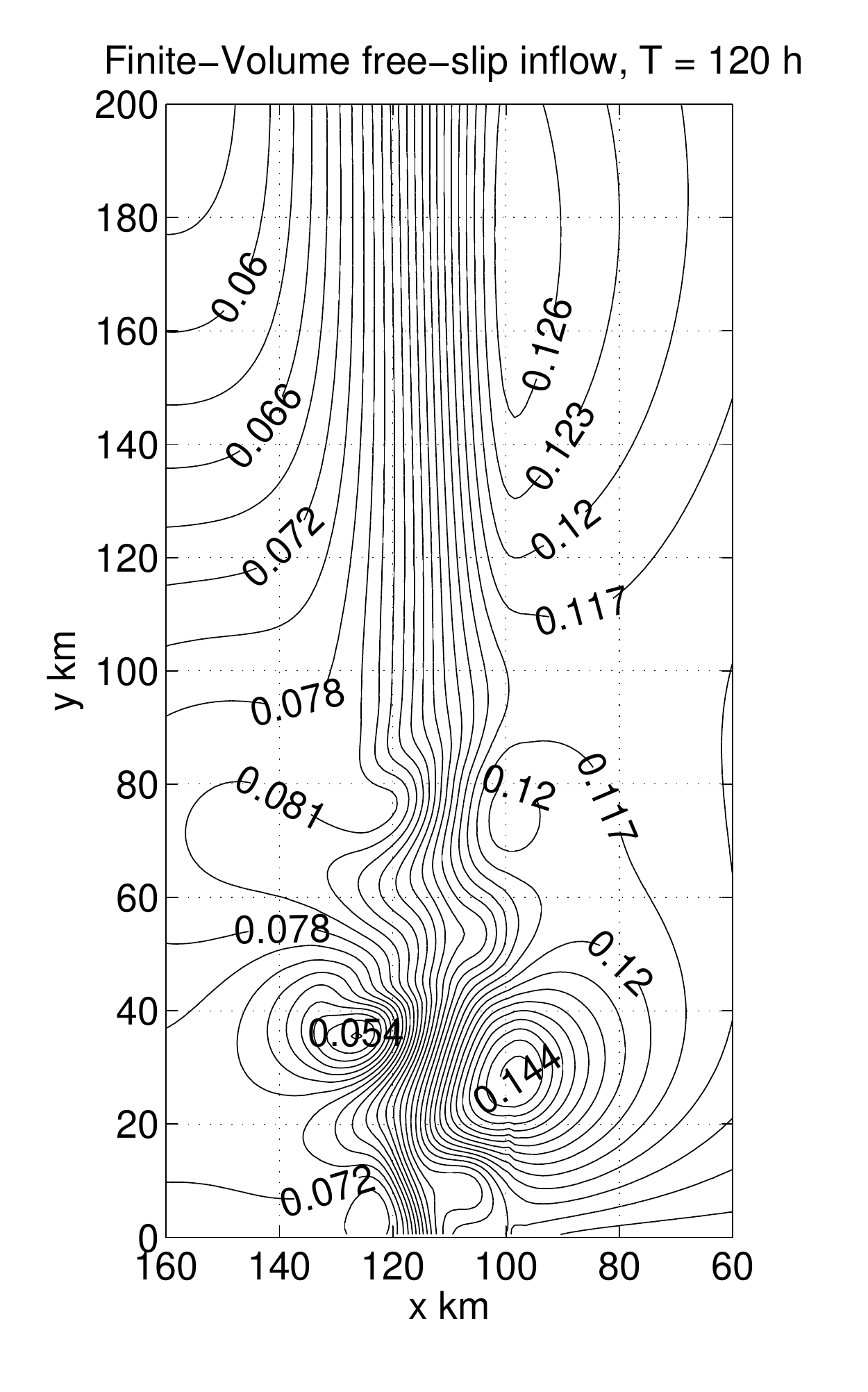}
  \includegraphics[width=0.32\linewidth]{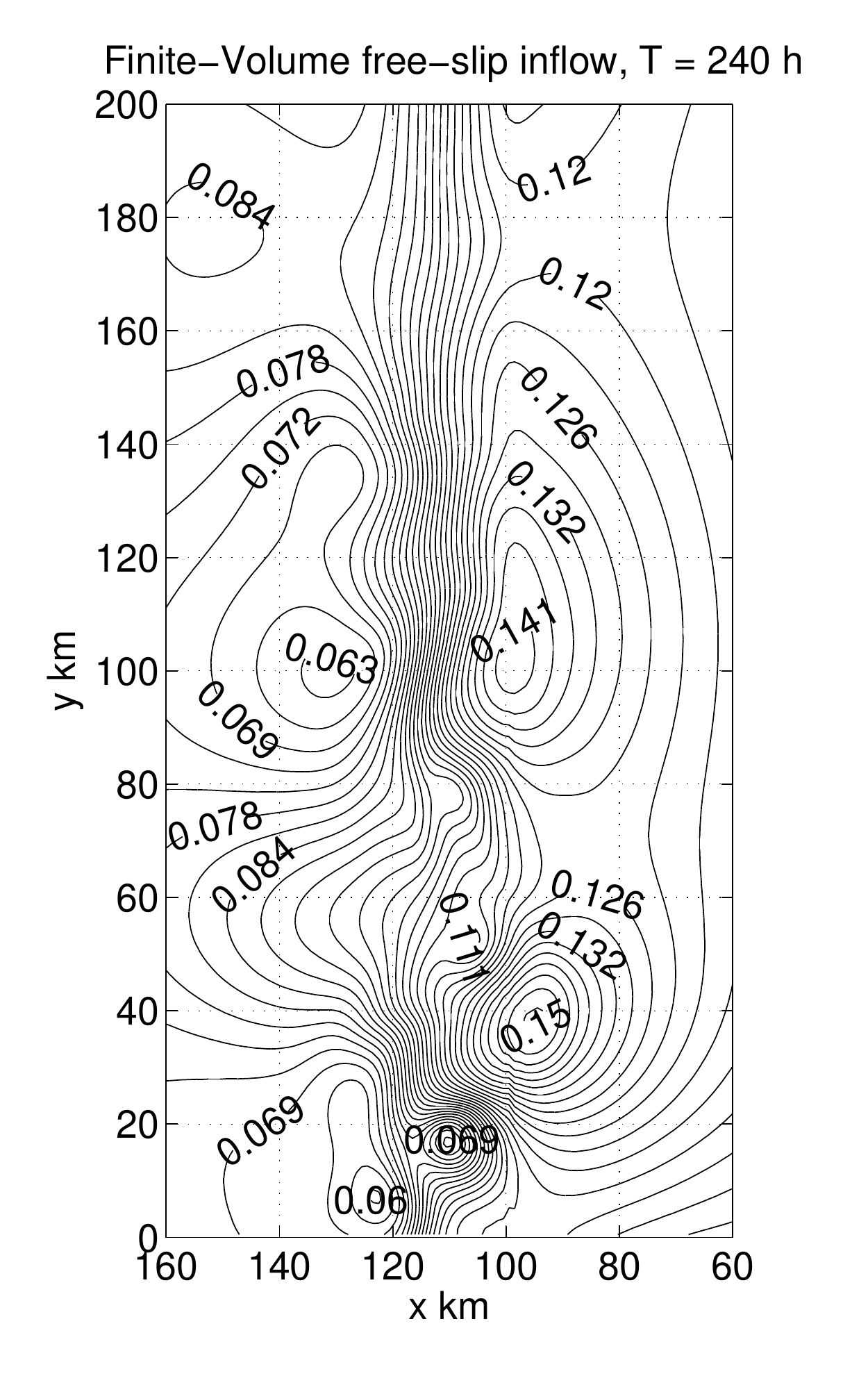}\\
  \caption{\small\label{fig:norwegian_shelfVelocityE3}Ormen Lange
  Experiment~III. Contour plots of surface elevation (top) and
  velocity plots (bottom) at $60$ (left), $120$, and $240$ (right)
  hours, computed with the Finite-Volume scheme free-slip boundary
  condition and smooth growth function
  Equation~\eqref{eq:growth_function}.}
\end{figure}
%===========================================================================================
\subsubsection{Setup for Ormen Lange Shelf Experiment~IV} \label{sec:BarostrophicJetSetupE4}
%===========================================================================================
Here we present another variant of the southern boundary condition. In
Section~\ref{sec:BarostrophicJetSetupE1} we used a discontinuous patch
of an in-flowing jet in the centre and open outflow at the periphery.
In Section~\ref{sec:BarostrophicJetSetupE2} we prescribed inflow
everywhere. Now we joint the inflow- and the open outflow boundary
conditions smoothly: Let $F_{i+\oh,j}^{infl}$ be the flux determined
by the free-slip inflow boundary condition and $F_{i+\oh,j}^{absorb}$
the one given by the absorbing outflow boundary condition. Now we use
the following convex combination to obtain the effective boundary flux
\begin{align}
F_{i+\oh,j} := \chi(x) F_{i+\oh,j}^{infl} + (1-\chi(x)) F_{i+\oh,j}^{absorb}.
\end{align}
The function $\chi(x)$, which prescribes the transition from the open
outer region towards the jet in the centre of the domain, is given by
\begin{align*}
  \chi(x)=\left \{
  \begin{array}{ll}
    0  &,   x < T_W - R,\\
    \Theta\left(\frac{x- T_W + R)}{2 R}\right) &, T_W - R \leq x \leq T_W + R,  \\
    1 &, T_W + R < x < T_E - R, \\
    \Theta\left(\frac{2 R - (x- T_E + R)}{2 R}\right) &, T_E - R \leq x \leq T_E + R,  \\  
    0 &,  x > T_E + R.
  \end{array}
  \right .
\end{align*}
Here 
$\Theta(x):= x^5\left[126 + x(-420 + x(540 + x(-315 + 70x)))\right]$
and the smoothing radius is $R = 5000$.  The transition
points are $T_W:= L_B-\oh B$ and $T_E:= L_B+\oh B$.
\begin{figure}
  \centering
  \includegraphics[width=0.8\linewidth]{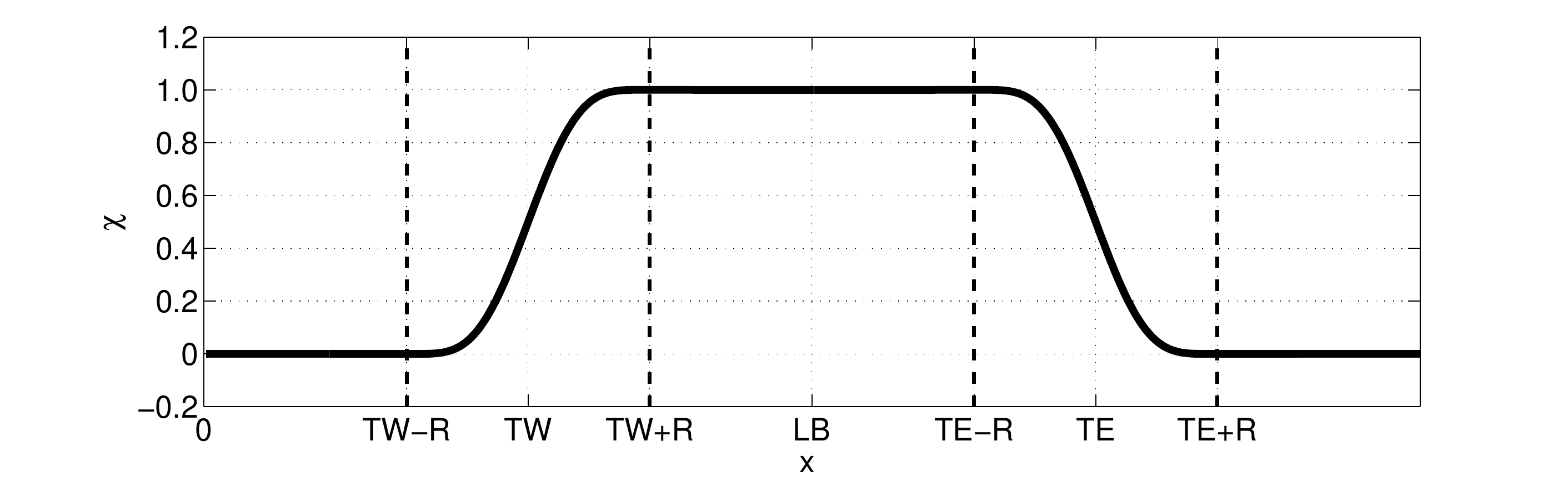}
  \caption{\small\label{fig:Theta} The function $\chi(x)$ prescribes
    the convex combination of the absorbing condition and the jet}
\end{figure}
\begin{figure}
  \centering
  \includegraphics[width=0.32\linewidth]{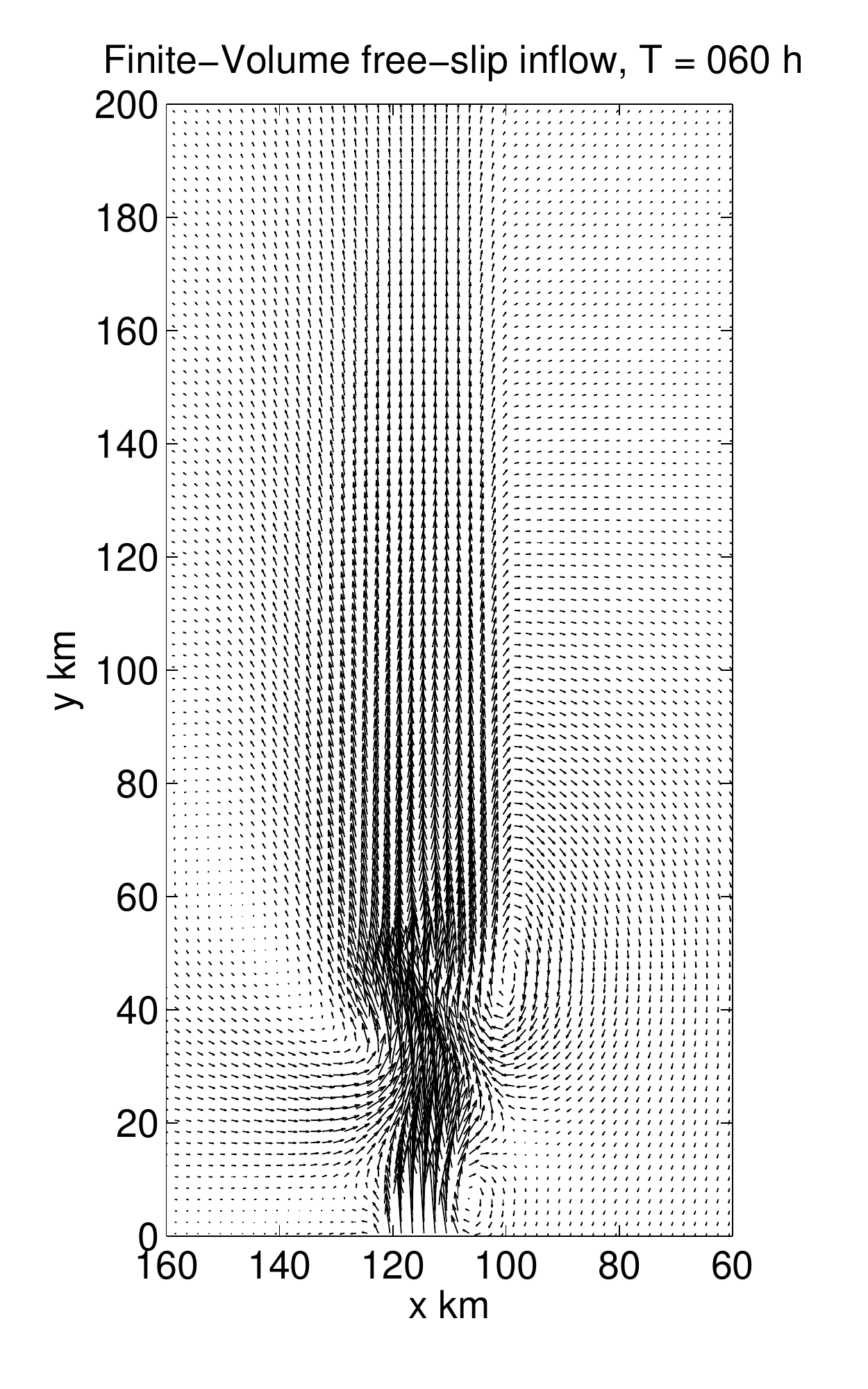}
  \includegraphics[width=0.32\linewidth]{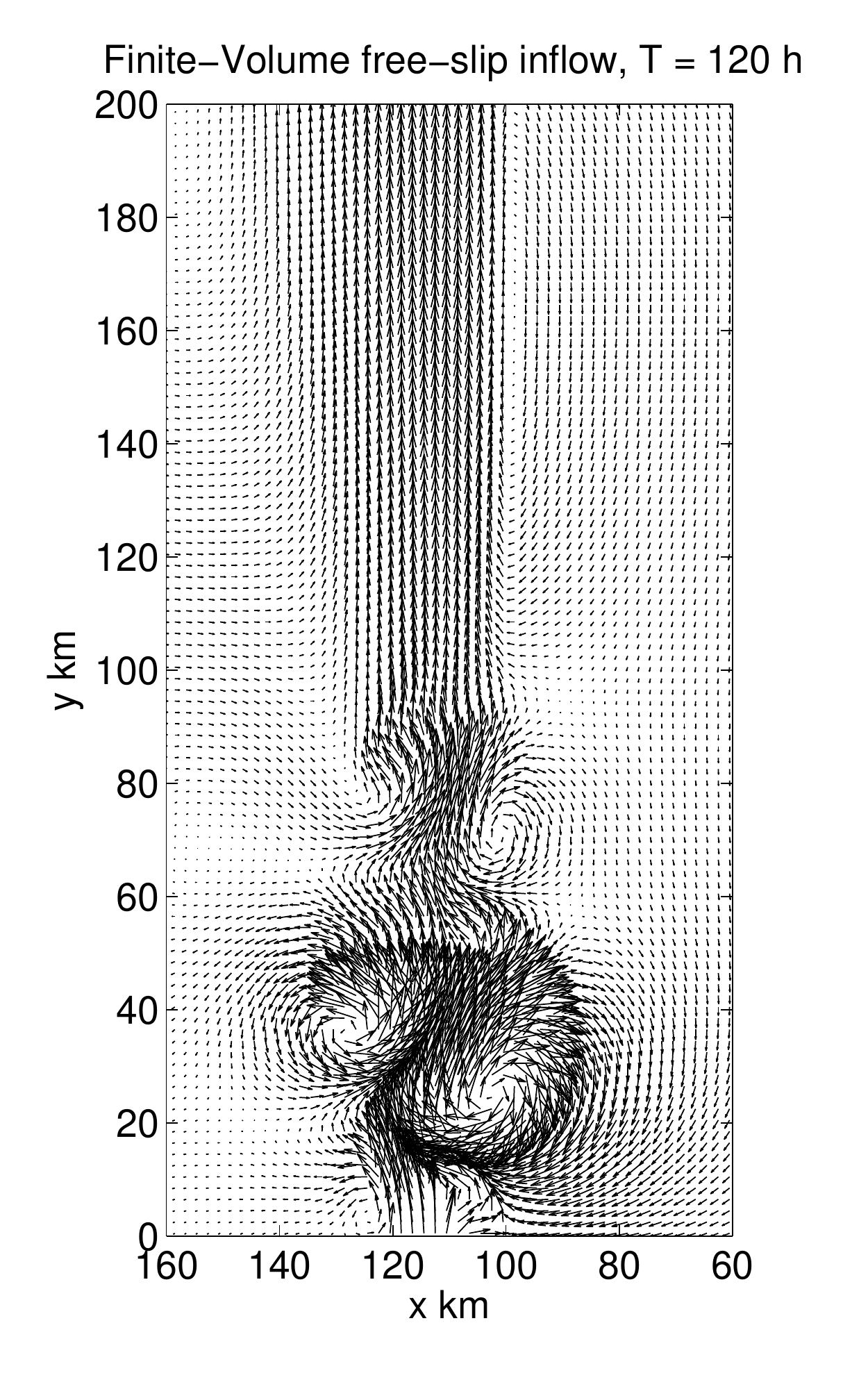}
  \includegraphics[width=0.32\linewidth]{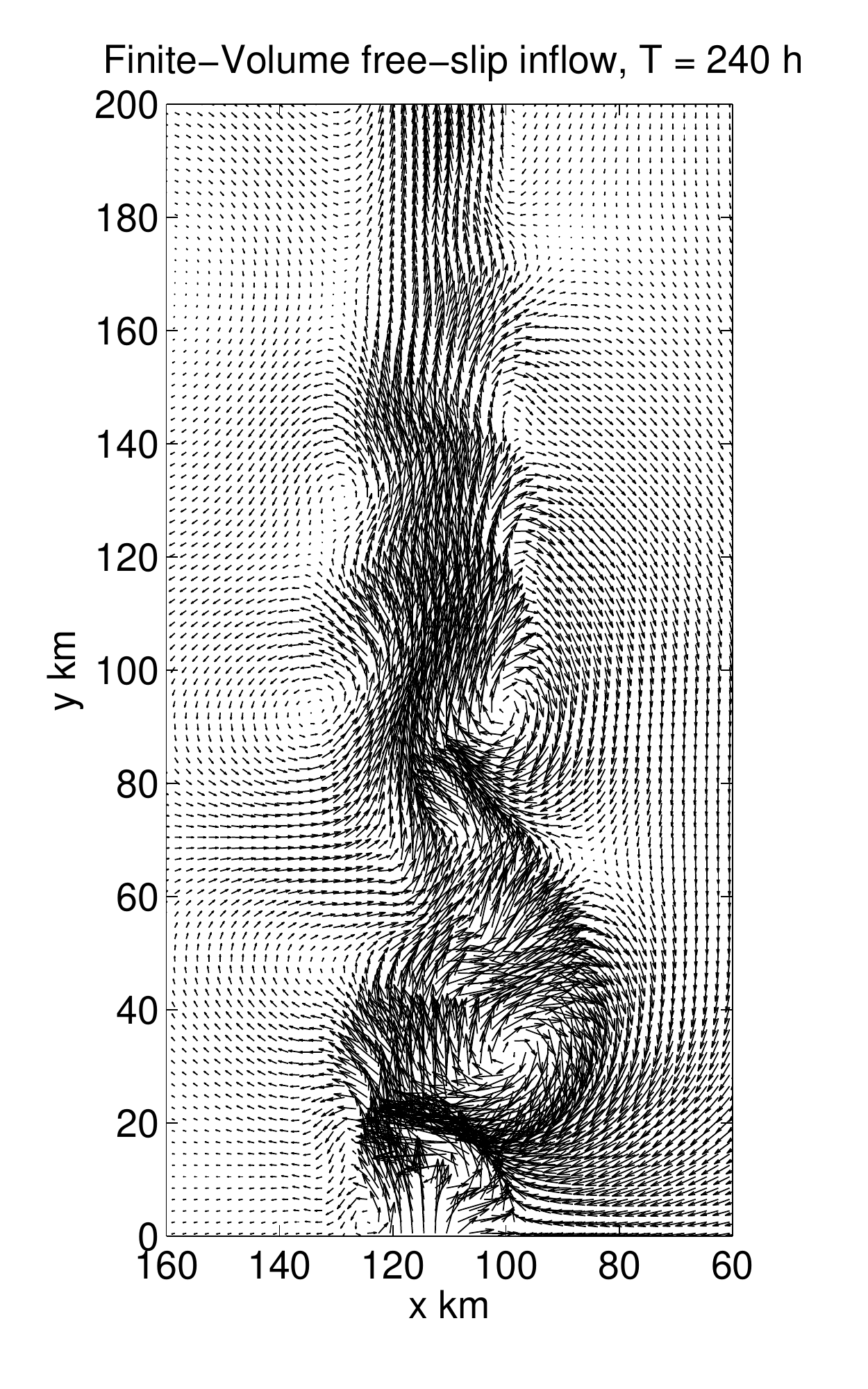}\\
  \vspace{-0.5cm}  
  \includegraphics[width=0.32\linewidth]{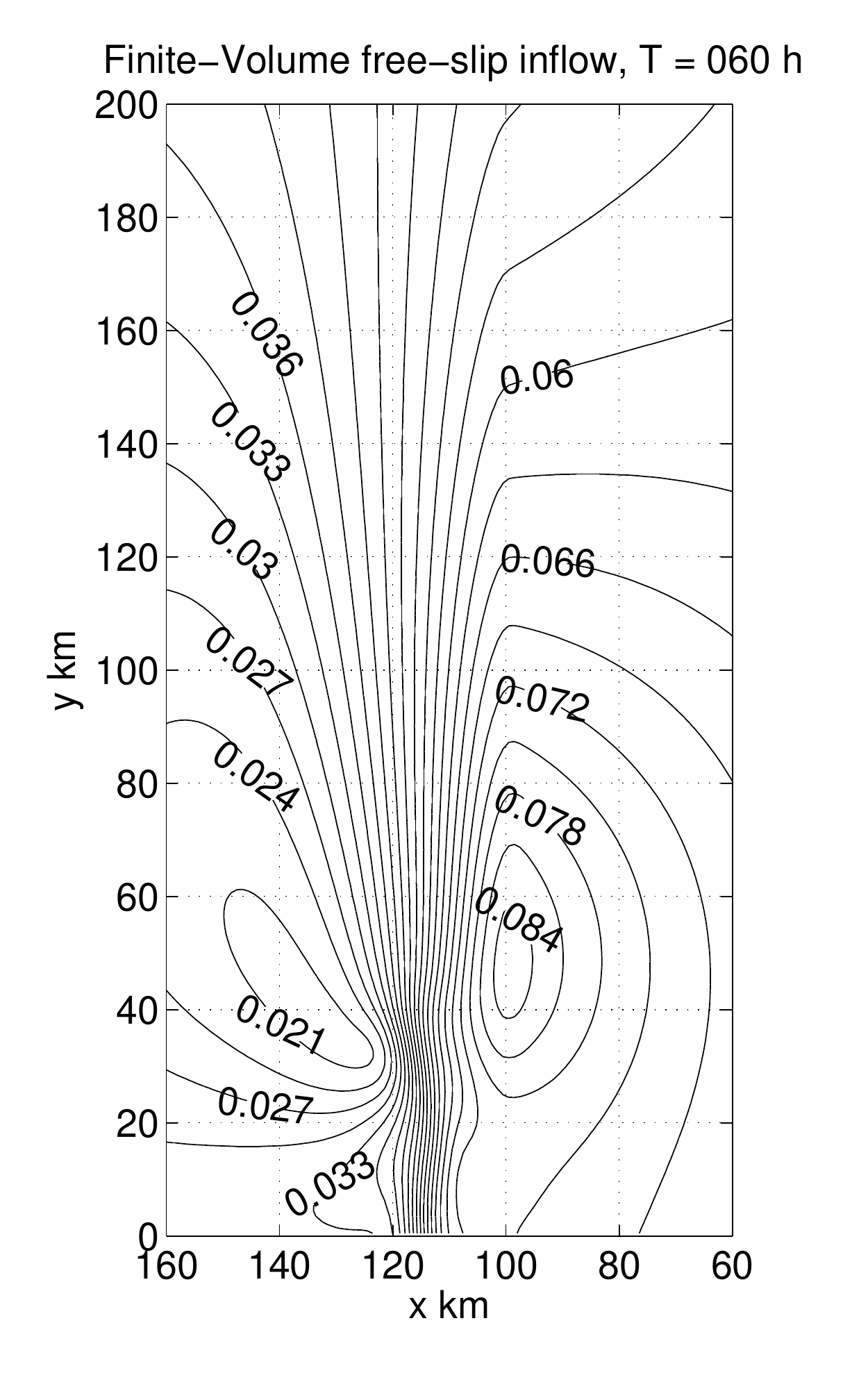}
  \includegraphics[width=0.32\linewidth]{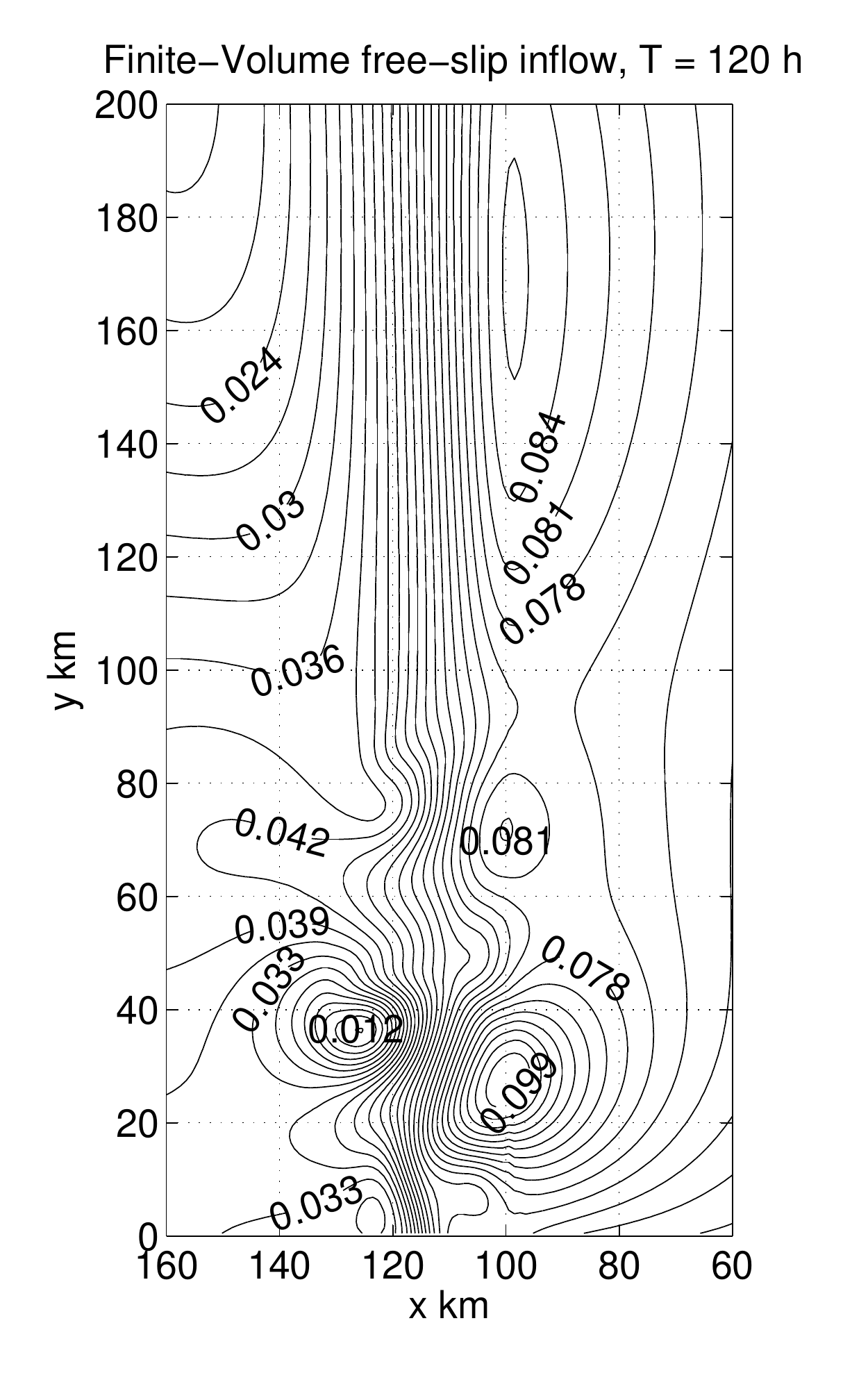}
  \includegraphics[width=0.32\linewidth]{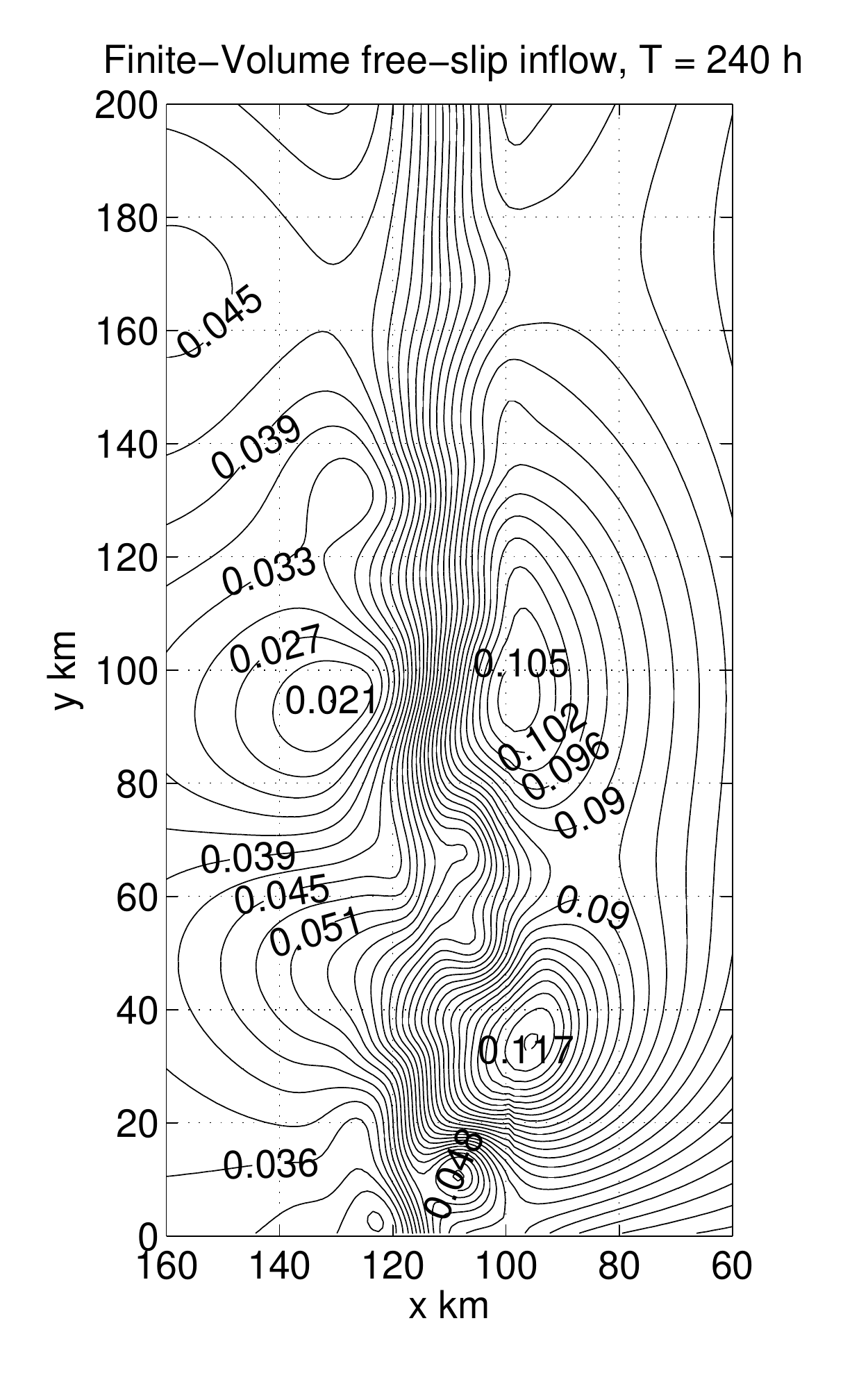}\\
  \caption{\small\label{fig:norwegian_shelfVelocityE4}Ormen Lange
  Experiment~IV.  Contour plots of surface elevation (top) and
  velocity plots (bottom) at $60$ (left), $120$, and $240$ (right)
  hours, computed with the Finite-Volume scheme free-slip boundary
  condition and exponential growth function
  Equation~\eqref{eq:growth_function} and transition smoothing.}
\end{figure}%
The numerical results are shown in
Figure~\ref{fig:norwegian_shelfVelocityE4}.  They agree in
considerable detail with the previous computations and hence confirm
the development of eddies, without introducing any discontinuity via
the numerical boundary treatment.
%===========================================================================================
\subsubsection{Balanced inflow boundary conditions: Ormen Lange Shelf Experiment~V} \label{sec:BarostrophicJetSetupE5}
%===========================================================================================
There remains one technical issue concerning the previous boundary condition:
the transition points $T_E$ and $T_W$ have to be chosen by hand. This
is not necessary for the volume-flux balanced boundary condition derived
in Section ~\ref{sec:FV_MomBalInf_BC} . There the decision of
outflow/inflow is taken automatically by the Riemann solver.

\begin{figure}
  \centering
  \includegraphics[width=0.32\linewidth]{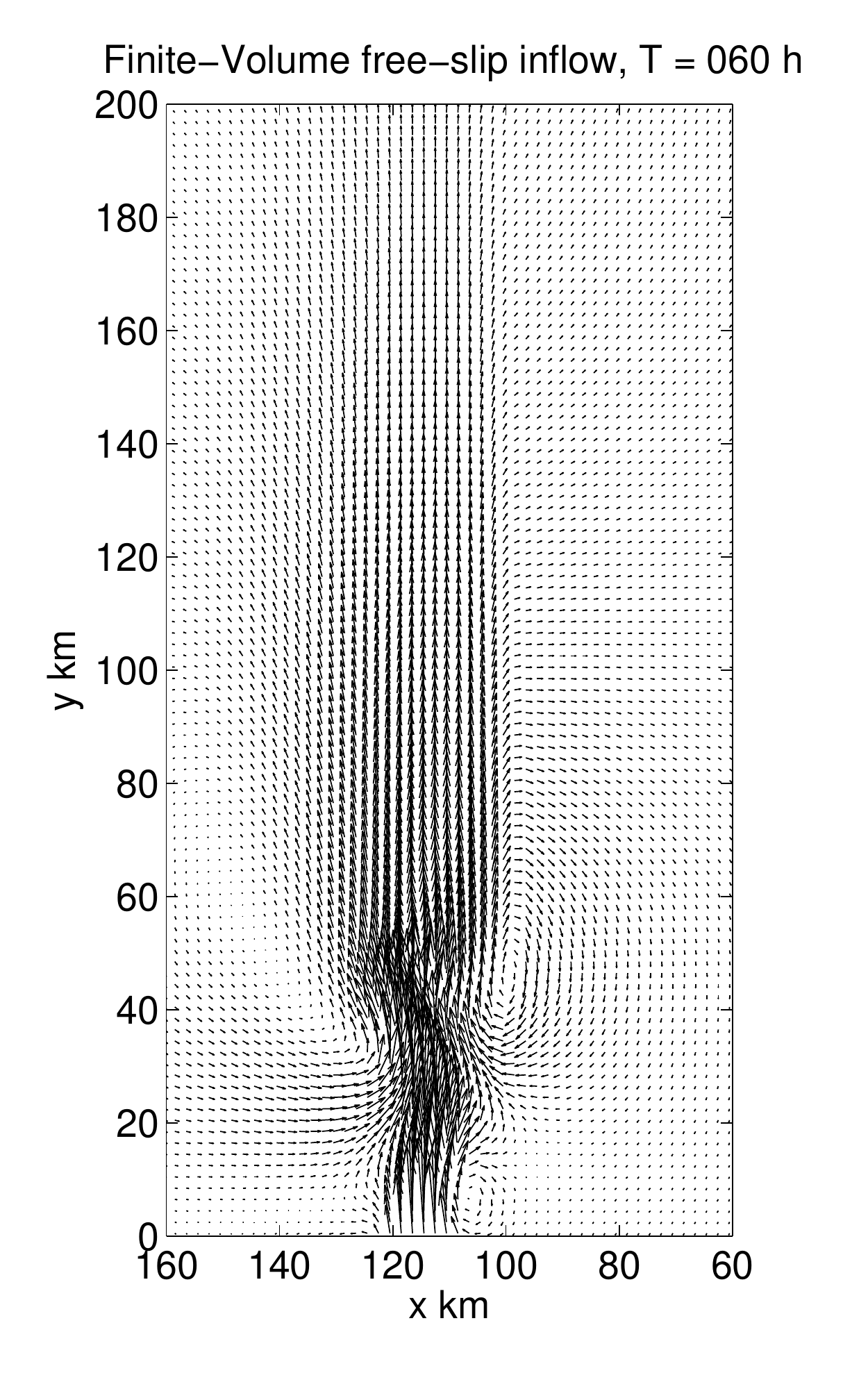}
  \includegraphics[width=0.32\linewidth]{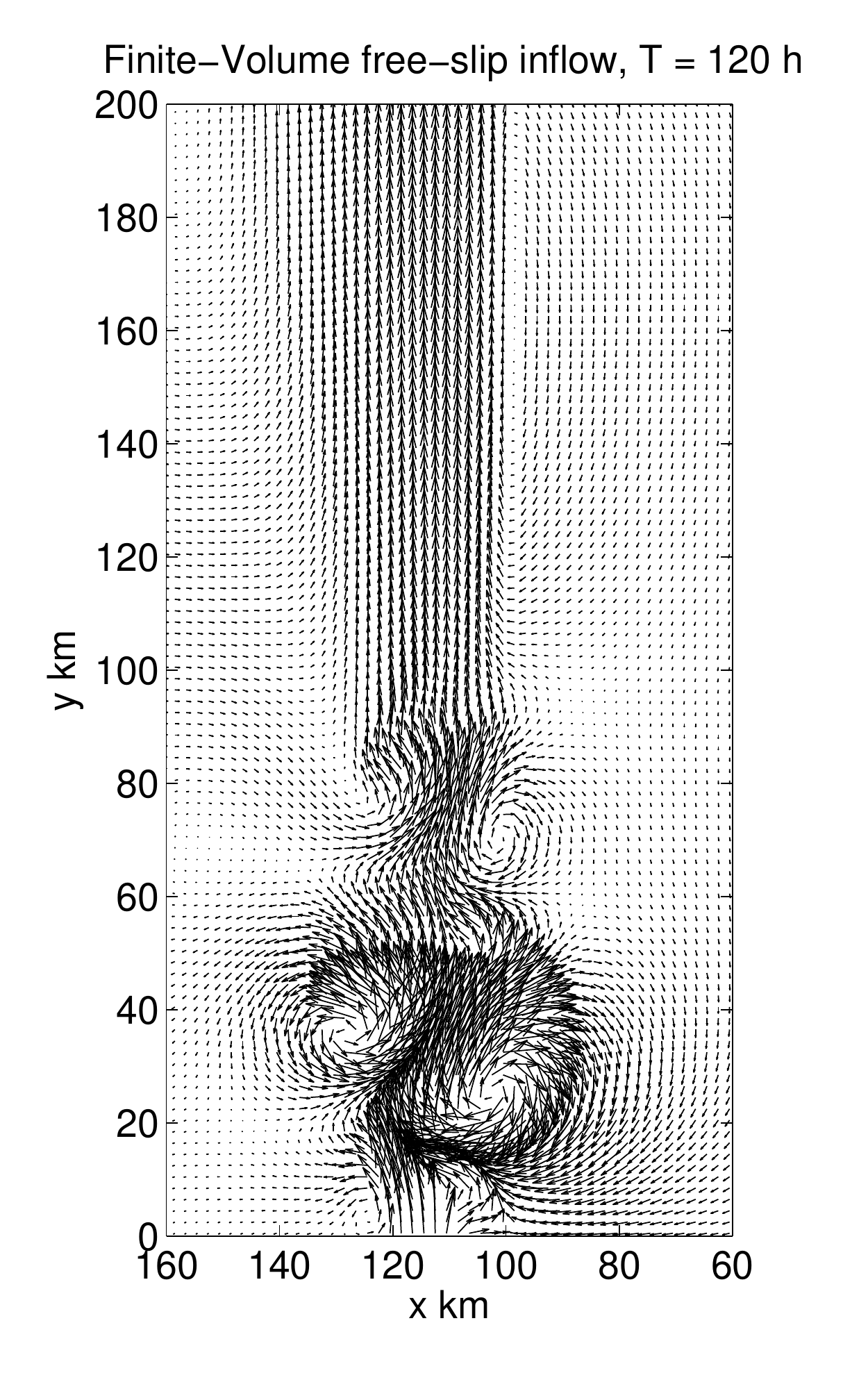}
  \includegraphics[width=0.32\linewidth]{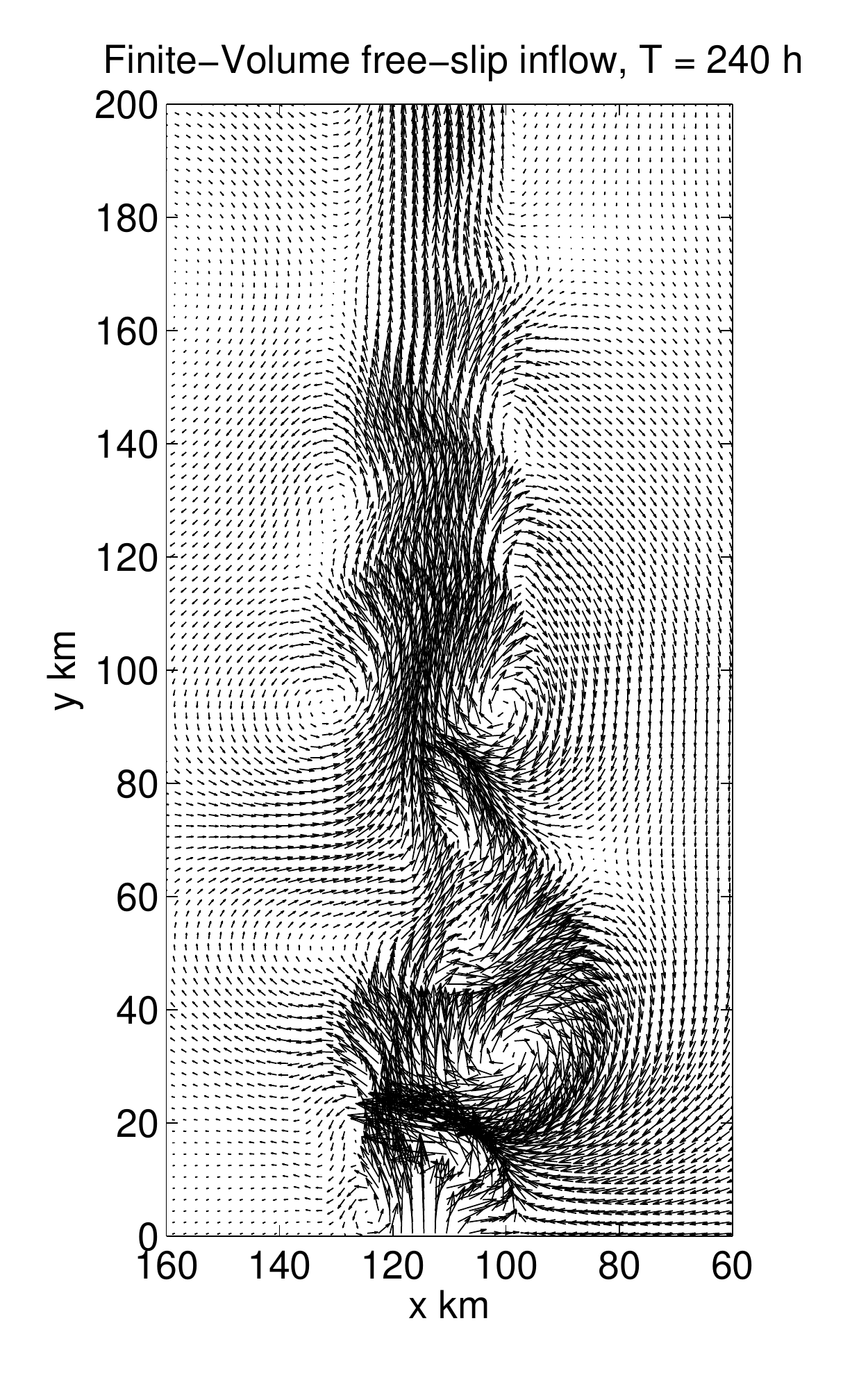}\\
  \vspace{-0.5cm}  
  \includegraphics[width=0.32\linewidth]{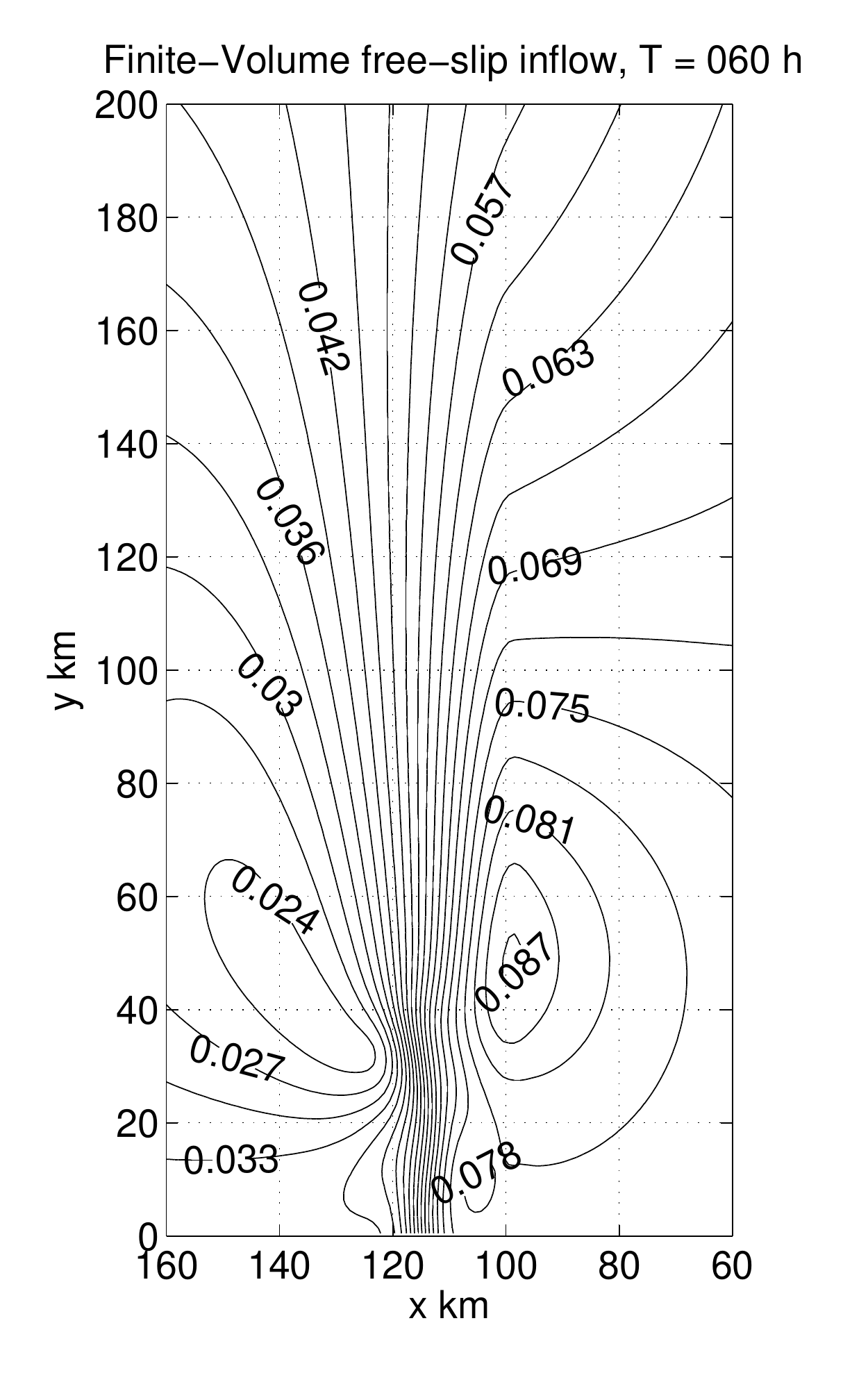}
  \includegraphics[width=0.32\linewidth]{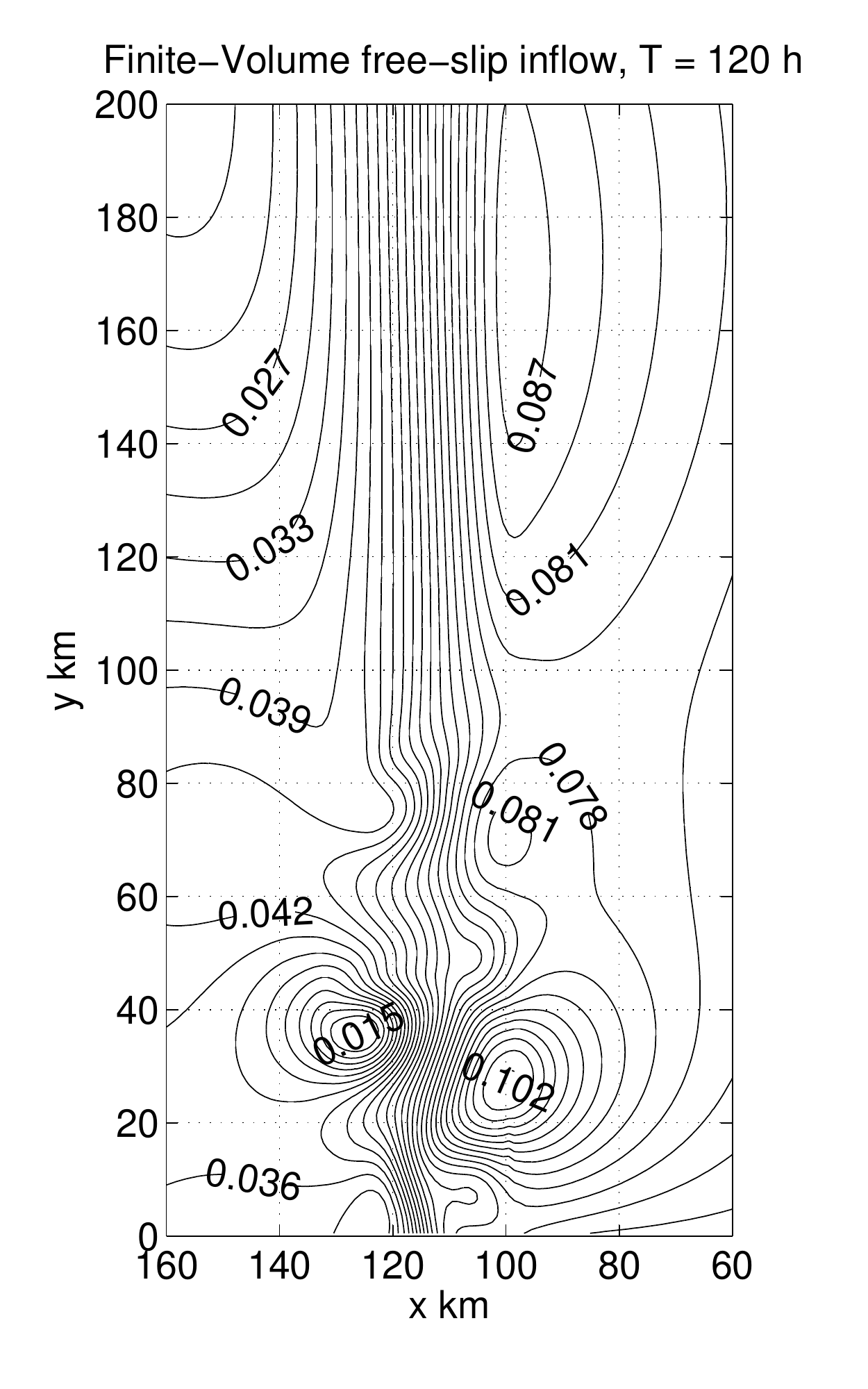}
  \includegraphics[width=0.32\linewidth]{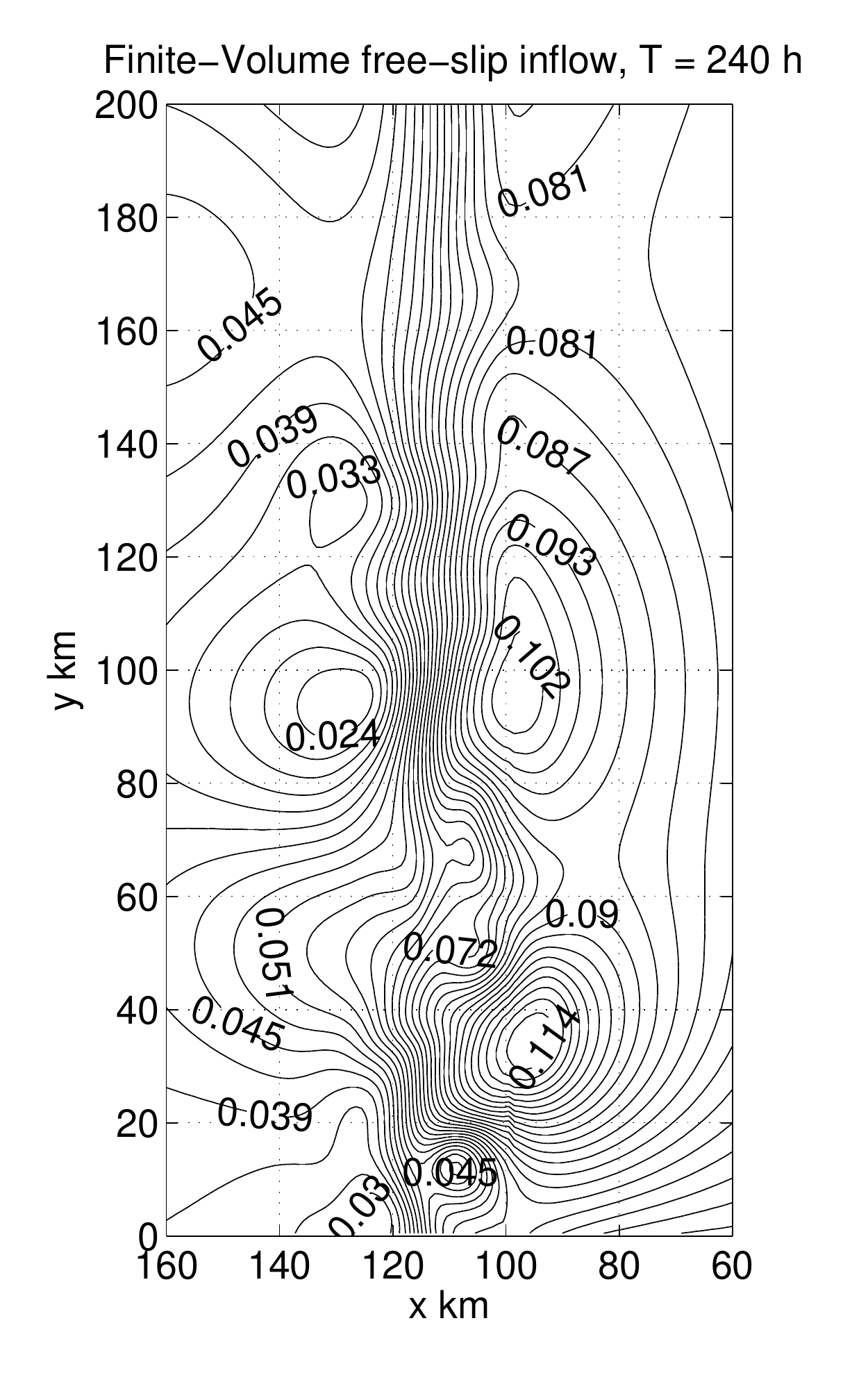}\\
  \caption{\small\label{fig:norwegian_shelfVelocityE6}Ormen Lange
  Experiment~VI. Contour plots of surface elevation (top) and
  velocity plots (bottom) at $60$ (left), $120$, and $240$ (right)
  hours, computed with the Finite-Volume scheme free-slip boundary
  condition, exponential growth function and balanced
  Equation~\eqref{eq:MomentumBalancePV2} at inflow.}
\end{figure}

The results for the balanced boundary condition are shown in
Figure~\ref{fig:norwegian_shelfVelocityE6}. They are in excellent
agreement with the results in
Figure~\ref{fig:norwegian_shelfVelocityE4}. This shows that the
volume-flux boundary condition is an interesting alternative to the
previous treatments, if we know the far-field values $\eta(x_0,y_0)$
and $U(x_0,y_0)$. The results for the geostrophically balanced boundary
condition are almost identical, and hence we do not display them here.

%\marginlabel{\color{red} Please check this carefully!}
%In all the computations, the eddy formation near the inflow boundary
%seems to generate a shelf wave with length of about $1000$km and
%period that equals the $35.6$ hour period of the eddies.  These long
%shelf waves appear with a strong signal in the surface elevation that
%grows considerably after ten days. The waves do not seem to affect the
%eddy formation at the inflow boundary, which is on a much shorter
%length scale.
%
%===========================================================================================
\subsubsection{Comparison with linear stability analysis} \label{sec:ComparisonLinearAnalysis}
%===========================================================================================
Linear stability analysis described in \cite{Gjevik_2002} and
\cite{Thiem_etal_2006} shows that the along shelf jet as defined for
the Ormen Lange case, section~\ref{sec:BarostrophicJetSetupE1}, is
unstable with respect to along shelf wave perturbations. The maximum
predicted exponential growth rate of $0.44$ day$^{-1}$ occurs for a
wave length of $44$~km.  The corresponding wave period is
$34.2$~hours. A second unstable mode has a maximum growth rate of
$0.28$ day$^{-1}$, a wave length of $54$~km and a period of
$41.1$~hours.  There are also steady, neutrally stable, shelf wave
oscillations in the band of wave lengths around $1000-1200$~km with
corresponding period $35.6-40.5$~hours.

\begin{figure}
  \centering
  \includegraphics[width=0.495\linewidth]{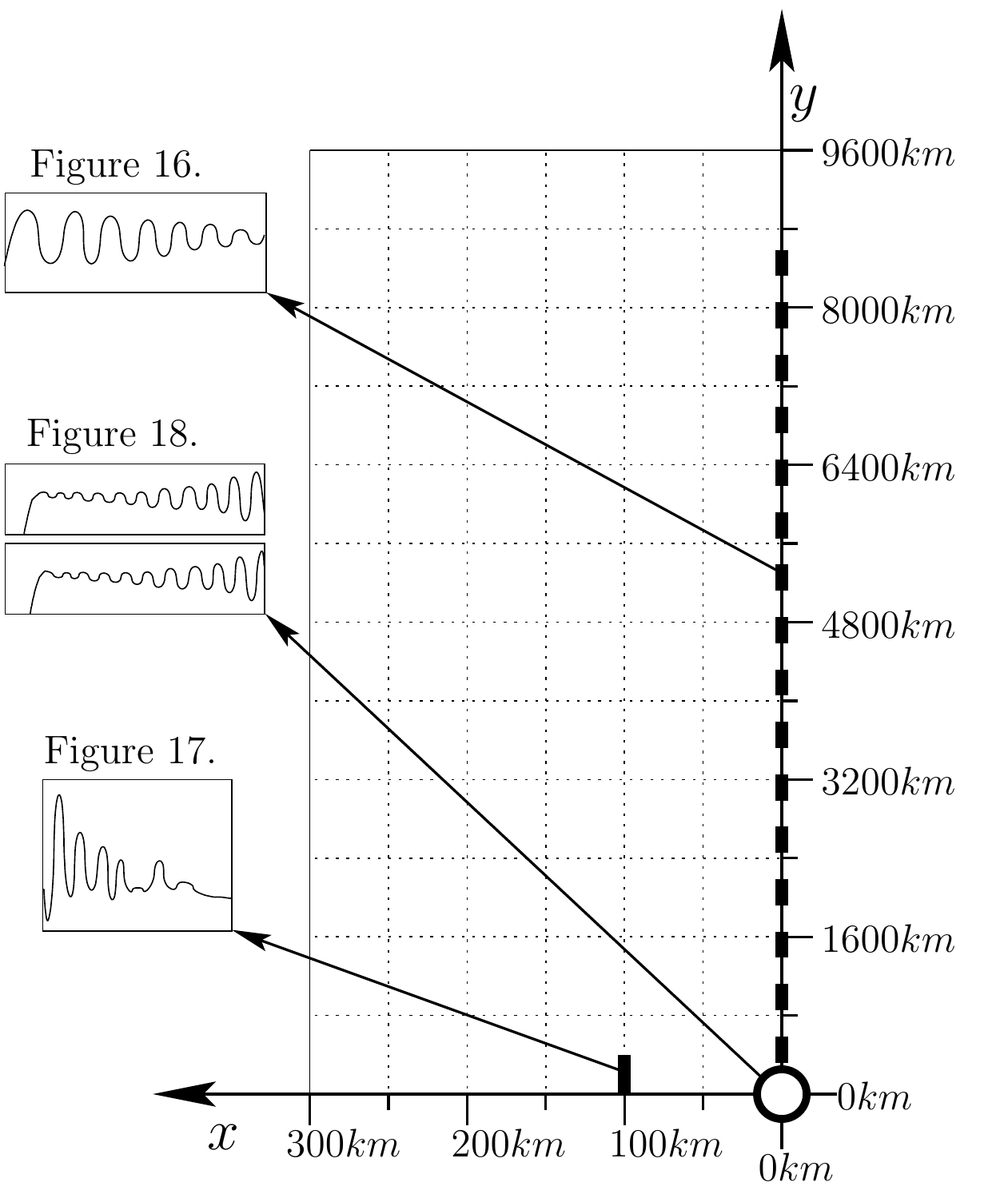}
   \caption{\small\label{fig:domainLarge} Large domain of $300 \times
     9600$~km$^2$, where the two cross-sections from
     Figure~\ref{fig:EddyAmpPlot_x100km} and
     Figure~\ref{fig:AmpPlot_x0km_T480h} are marked, the position from
     which the frequency plots Figure~\ref{fig:Frequency_x0km_y0km}
     are taken is marked with the black circle.  }
\end{figure}

\begin{figure}
  \centering
  \includegraphics[width=0.495\linewidth]{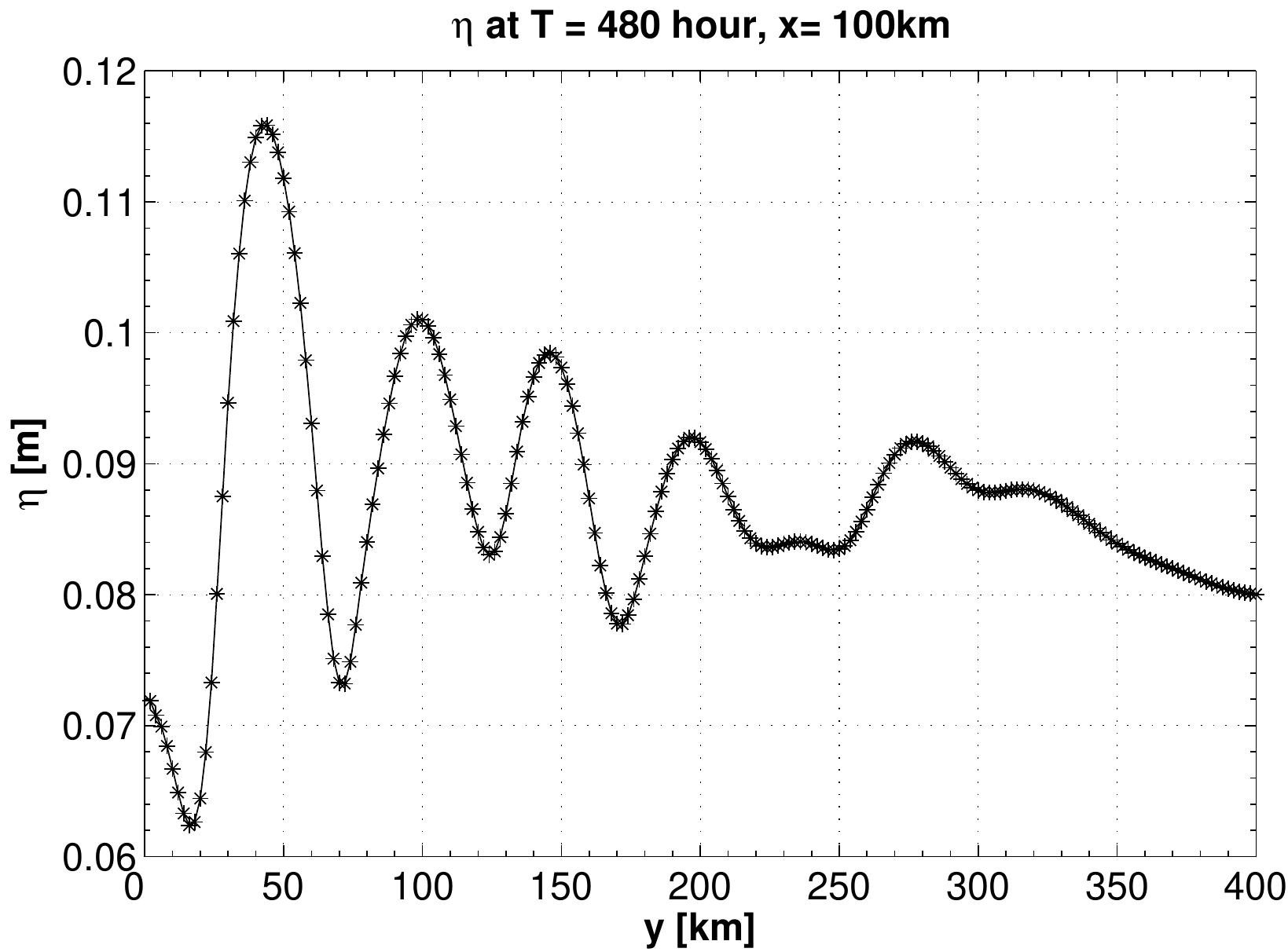}
   \caption{\small\label{fig:EddyAmpPlot_x100km}} Section of surface
   displacement at $x=100$~km, computed on a domain of $300 \times
   9600$~km$^2$. The maximum exponential growth rate is observed for a
   wave length of $49$~km.
\end{figure}

\begin{figure}
  \centering
  \includegraphics[width=\linewidth]{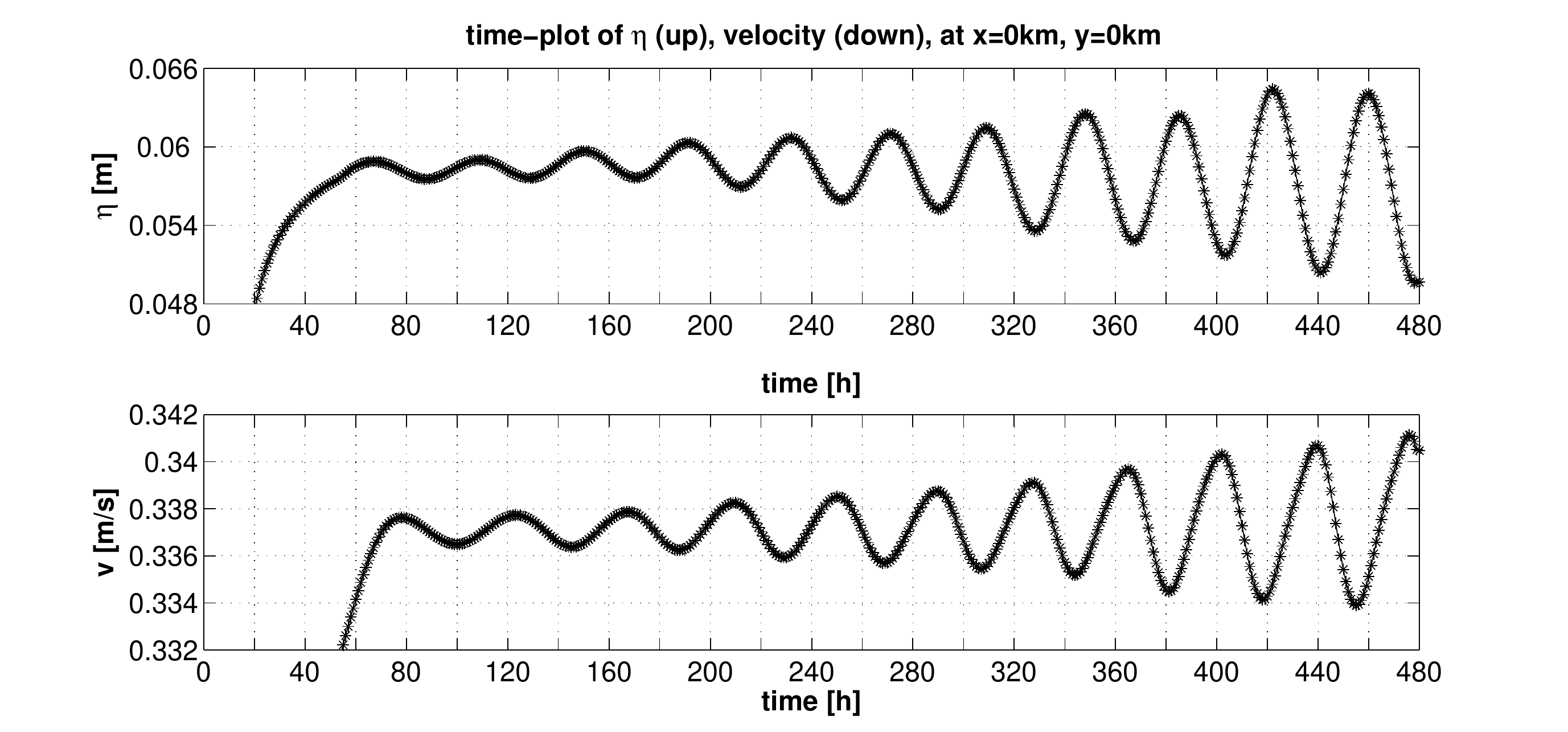}
   \caption{\small\label{fig:Frequency_x0km_y0km} The frequency of the
     maximum unstable waves is approximately $39.5$ hour.  The upper
     plot shows the frequency in the surface displacement and the
     lower plot shows the the frequency in normal velocity $v$.  }
\end{figure}

To compare the results of the linear stability analysis with the
solution of the finite volume scheme in more detail we did the same
computation as in section~\ref{sec:BarostrophicJetSetupE1} on an
enlarged domain of $300 \times 9600 \textnormal{km}^2$, grid-width
$2$~km and final time $480$~hours (see Figure~\ref{fig:domainLarge}).
Since several periods of the long waves (wavelength $1000 - 1200$~km)
fit into this domain, it is possible to measure the wavelength very
accurately.
Figure~\ref{fig:EddyAmpPlot_x100km} shows the surface elevation for
the section $x=100$~km, $0$~km $\leq y \leq 400$~km, which is the
upper shelf-edge. It is here that we observe the strongest wave
amplitudes. The peaks indicate a wave length of $49$~km, which is the
distance of the eddies.
\begin{figure}
  \centering
  \includegraphics[width=0.9\linewidth]{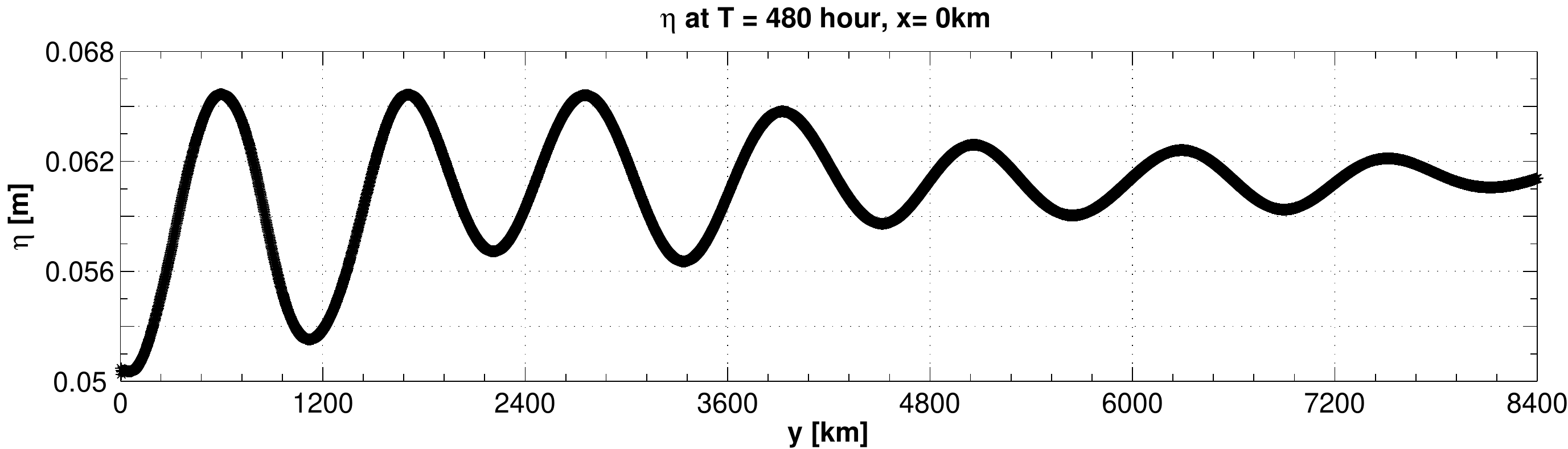}
   \caption{\small\label{fig:AmpPlot_x0km_T480h}.The wave length of
     the neutrally stable long- wave at the eastern boundary $x=0$~km
     is $1152$~km.}
\end{figure}
Figure~\ref{fig:AmpPlot_x0km_T480h} shows the surface elevation for 
the section along the coast ($x=0$~km, $0$~km $\leq y \leq 8400$~km). Here we 
observe the second strongest wave amplitudes. The peaks indicate a  
wave length of $1152$~km. 
Both wave, the one in Figure~\ref{fig:EddyAmpPlot_x100km} with wave
length $49$~km, and the one in Figure~\ref{fig:AmpPlot_x0km_T480h}
with $1152$~km, have the same wave period of $39.5$~hours. For the
second wave, time-plots of surface displacement, and velocity $v$ are
shown in Figure~\ref{fig:Frequency_x0km_y0km}. Time plots for the
first wave are similar, and not shown here.

The computational results of Figure~\ref{fig:EddyAmpPlot_x100km},
showing the development of eddies with period $35-40$~hours and an
along shelf separation of $40-60$~km, are in close, but not complete,
agreement with the predictions of the linear stability analysis.

The oscillation with wave length of about $1200$~km and period about
$40$~hours (see Figure~\ref{fig:AmpPlot_x0km_T480h}) is most
pronounced in the sea level $\eta$ and its amplitude grows
considerably over a time span of $10$~days. Clearly, this oscillation
corresponds to the steady long shelf wave oscillations found by the
stability analysis.  In the numerical simulations the oscillation
seems to be excited by the periodic eddy formation near the inflow
boundary and propagates subsequently downstream with a speed of about
$30$~km/hours.

A perfect correspondence between the linear stability analysis and the
numerical simulations of the inflow jet cannot be expected due to
nonlinear effects and the downstream development of the eddies in the
model.

Note that for the situations computed above, the finite difference
scheme yields almost equal results as the finite volume scheme.
\iffalse
%
\begin{figure}
  \centering
  \includegraphics[width=0.32\linewidth]{figures/FV_freeslip_vel_GeoBal_T060h}
  \includegraphics[width=0.32\linewidth]{figures/FV_freeslip_vel_GeoBal_T120h}
  \includegraphics[width=0.32\linewidth]{figures/FV_freeslip_vel_GeoBal_T240h}\\
  \vspace{-0.5cm}  
  \includegraphics[width=0.32\linewidth]{figures/FV_freeslip_eta_GeoBal_T060h}
  \includegraphics[width=0.32\linewidth]{figures/FV_freeslip_eta_GeoBal_T120h}
  \includegraphics[width=0.32\linewidth]{figures/FV_freeslip_eta_GeoBal_T240h}\\
  \caption{\small\label{fig:norwegian_shelfVelocityE5}Ormen Lange
  Experiment~V.  Contour plots of surface elevation (top) and
  velocity plots (bottom) at $60$ (left), $120$, and $240$ (right)
  hours, computed with the Finite-Volume scheme free-slip boundary
  condition, exponential growth function and geostrophic balanced
  Equation~\eqref{eq:GeostrophicBalance} inflow.}
\end{figure}
\fi
%
%
%
%===============================================================================
\section{Conclusion}\label{sec:Conclusion}
%===============================================================================
%
%
In this paper we have presented a comparison of a finite-difference
and a high-order finite-volume scheme for geophysical flow problems.
To conclude our report we discuss 
\begin{enumerate}
\item[1.]
the efficiency and stability of the FD and FV solvers,
\item[2.]
the numerical inflow boundary conditions,
\item[3.]
the geophysical implication of the computational results.
\end{enumerate}
\subsection{Efficiency and stability of the FD and FV solvers}
The results indicate that the two schemes compute qualitatively and
quantitatively similar solutions.  The rates of convergence are as
expected, i.e., first order for the original finite-difference scheme
used by Gjevik et al.  in \cite{Gjevik_etal_2002}, second order for
the modified finite-difference scheme, and fourth order (almost fifth)
for the finite-volume scheme.

An exact quantitative comparison of run-times is not possible, since
the two codes are research codes written in different languages and by
different programmers. However, it is fair to say that for the very
smooth test problems \ref{sec:TestingOrder} (see
Table~\ref{table:convergence2D_improved_enlag} and
Table~\ref{table:convergence2D}) and
\ref{sec:LargeEddiesDoublePeriodicDomain} the higher order finite
volume scheme is asymptotically more efficient. For the more
realistic, and less smooth, test problems
Section~\ref{sec:jetexample}, both the FD and the FV code give
qualitatively the same results on the same grid, but the FD scheme is
much faster than the FV scheme. It would be desirable to run the FV
scheme on a coarser grid to reduce the runtime. But then the inflow
data for the jet would be resolved by less then 10 cells, and the flow
is not sufficiently resolved any more. This might be different for
broader currents. Note that the FV scheme could resolve small gravity
waves which were completely smeared by the FD scheme.  However, these
waves quickly leave the computational domain and do not seem to have a
noticeable impact on the major currents.

But the FV scheme has an important advantage over the FD scheme: it is
much more stable in cases with strong gradient. We can run it with CFL numbers of 0.5 (in all our
computations) and sometimes up to 1, without adding any artificial
viscosity. This includes solutions with shock-like discontinuities,
e.g. hydraulic jumps. If we run the FD scheme without artificial viscosity,
and for smooth solutions, it may already produce instabilities for CFL
numbers of 0.5. This happened for example when we implemented the
free-slip boundary condition into the FD scheme. For hydraulic jumps,
a lot of artificial viscosity has to be added to stabilise the FD
scheme, and this reduces the accuracy of the scheme.

\subsection{Numerical inflow boundary conditions}
%\marginlabel{\color{red} Please check the implementation of\\ transparent inflow conditions.}
Using Riemann decompositions, we could successfully translate the FD
boundary conditions to the FV solver. We could also improve the
no-slip inflow boundary conditions by the free-slip condition, which
yields smoother solutions
(Tables~\ref{table:conv2D_nonsmooth_FD}-\ref{table:conv2D_smooth_FV}).
Moreover, we developed a transparent inflow condition, which allows 
waves to leave the domain through the inflow boundary.

\subsection{Geophysical implication of the computational results}
Various numerical experiments for the Ormen Lange cases presented in
Section~\ref{sec:jetexample}, with the FD and the FV schemes and
different implementations of the boundary conditions led to almost
identical results for two different startup profile configurations.
These results are also in close agreement with linear stability
analysis (see Section~\ref{sec:ComparisonLinearAnalysis}). Therefore
the computations presented here fully confirm the results of
\cite{Thiem_etal_2006} about instabilities of the shelf slope jet and
the formation of eddies.

\subsection{Further perspectives}
The finite-volume scheme is much more expensive with respect with
computer time than the traditional finite-difference scheme, but one
benefit from a lot higher accuracy.  To obtain a similar accuracy with
the finite-difference scheme one would have to refine the grid several
times.
%???
This may make the finite-volume scheme attractive for studies of high
frequency oscillations associated with strong current shears or small
scale bathymetric features on the shelf edge.

%We believe to have obtained a solid bases for further studies of 
%barostrophic flows, including more complex topography, bottom
%friction, turbulence modelling, tidal currents, multi layer
%models, and wind forces.

%
%
%
%===============================================================================

\end{document}